\documentclass[12 pt]{amsart}

\usepackage{amscd,amsfonts,amssymb,amsmath,amsthm,enumerate,latexsym,floatflt,longtable}
\usepackage{xcolor}
\usepackage{placeins}
\usepackage{hyperref}

\usepackage{adjustbox}

\usepackage[utf8]{inputenc}

\title[The Table of the Structure Constants]{The Table of the Structure Constants for the Complex Simple Lie Algebra of Type $F_4$
and its Application to the Calculation of Commutators in the Chevalley Group of Type $F_4$ over Fields and Rings}

\author[Sergey G.~Kolesnikov,  Anna I.~Polovinkina]{Sergey G.~Kolesnikov,  Anna I.~Polovinkina}

\date{\today}

\begin{document}
\sloppy



\maketitle
\begin{abstract}


This work is the first in a series of papers devoted to constructing tables of structure constants for the complex simple Lie algebras and to finding an explicit form of Chevalley commutator formulas.

The work consists of three parts. In the first part, expressions are found for the structure constants of the complex
simple Lie algebra of type $F_4$ in the form of functions of structure constants corresponding to
extraspecial pairs of roots. As a consequence, all Chevalley commuta\-tor formulas $[x_r(u),x_s(y)]$ 
are calculated when the sum $r+s$ is a root.

Further, in the second part, tables of structure constants and Chevalley commutator formulas are given in the special 
case when all constants corresponding to extraspecial pairs are equal to one. 

Finally, in the third part, directed and weighted graphs associated with root systems are constructed. 
It is shown that the elements of the exponent of the adjacency matrices of directed graphs are the numbers $P_{rs}$, where $P_{rs}$ is the number of representations of the root $r$ in the form of a sum of the root $s$ and fundamental roots such 
that any initial segment of the sum is a root. It is also shown that the elements of an exponent of the weight 
matrix of a weighted graph are the values of sums arising when calculating complex commutators in Chevalley groups.


\end{abstract}

\maketitle
\tableofcontents





\section{The Table of the Structure Constants}

\subsection{Information about the Root System of Type $F_4$}

Following \cite{Bur72}, we choose the orthonormal basis $e_1,e_2,e_3,e_4$ in the Euclidean space $V=\mathbb{R}^4$  with the scalar product $(\,,\,).$  Then the following set of vectors:
$$
\pm e_i\ (1\leqslant i\leqslant 4),
\quad
\pm e_i \pm e_j\ (1\leqslant i\leqslant 4),
\quad
(\pm e_1 \pm e_2 \pm e_3 \pm e_4)/2
$$
forms the reduced indecomposable root system of type $F_4.$

The following subsets
$$
\{ a=e_2-e_1,\quad b=e_3-e_2,\quad c=e_4,\quad d=(e_1-e_2-e_3-e_4)/2 \}
$$
and
$$
\{ e_i\ (1\leqslant i\leqslant 4),
\quad
e_i \pm e_j\ (1\leqslant i\leqslant 4),
\quad
(e_1 \pm e_2 \pm e_3 \pm e_4)/2 \}
$$
form, respectively, the system of fundamental roots $\Pi(F_4)$ and the system of positive roots
$F_4^+$ in $F_4.$ Next, for convenience, we represent the positive root
$$
r=\alpha a + \beta b + \gamma c + \delta d
$$
by the quadruple of numbers $\alpha\beta\gamma\delta.$ The sum $\alpha + \beta + \gamma + \delta$
is called the height of the root $r$ and is denoted by ${\rm ht}(r).$

We define the order $\prec$ on $V$, consistent with the choice of the fundamental root system, as follows.
For arbitrary vectors $x,y\in V$ we assume
$$
x = \tau_1 a+\tau_2 b+\tau_3 c+\tau_4 d \prec y = \rho_1 a+\rho_2 b+\rho_3 c+\rho_4 d
$$
if and only if the first coefficient in the expansion of the difference
$$
y-x = (\rho_1-\tau_1) a+(\rho_2-\tau_2) b+(\rho_3-\tau_3) c+(\rho_4-\tau_4) d
$$
is positive. With respect to the order $\prec$ the positive roots of the root system $F_4$ are ordered in the following way:
$$
0001 \prec 0010 \prec 0011 \prec 0100 \prec 0110 \prec 0111 \prec 0120 \prec 0121 \prec 0122 \prec 1000 \prec 1100 \prec 
$$
$$
1110 \prec 1111 \prec 1120 \prec 1121 \prec 1122 \prec 1220 \prec 1221 \prec 1222 \prec 1231 \prec 1232 \prec 1242 \prec
$$
$$
1342 \prec 2342.
$$

The following table shows the correspondence between the number of a positive root in the ordered list and its quadruple representation.
The table also shows root height and length.
\vskip0.5cm

\begin{center}
\begin{tabular}{|c|l|c|c|c|c|}
    \hline
        ¹  & $r$           & $\alpha\beta\gamma\delta$ & $\theta e_1+\rho e_2+\nu  e_3+\mu e_4$       & ${\rm ht}(r)$ & $(r,r)$ \\ \hline
        1  & $d$           & 0001                      & $(e_1-e_2-e_3-e_4)/2$                        & 1             & 1       \\ \hline
        2  & $c$           & 0010                      & $e_4$                                        & 1             & 1       \\ \hline
        3  & $c+d$         & 0011                      & $(e_1-e_2-e_3+e_4)/2$                        & 2             & 1       \\ \hline
        4  & $b$           & 0100                      & $e_3-e_4$                                    & 1             & 2       \\ \hline
        5  & $b+c$         & 0110                      & $e_3$                                        & 2             & 1       \\ \hline
        6  & $b+c+d$       & 0111                      & $(e_1-e_2+e_3-e_4)/2$                        & 3             & 1       \\ \hline
        7  & $b+2c$        & 0120                      & $e_3+e_4$                                    & 3             & 2       \\
\hline
\end{tabular}
\end{center}

\centerline{\textit{Table 1.}}

\subsection{Calculation of Structure Constants} In this section we proved the following
\medskip

\textbf{Theorem 1.}
\textit{Let the structure constants $N_{r,s},$ corresponding to extraspecial pairs of roots be chosen as indicated in Table~2.}
\vskip5mm

\begin{center}
\begin{tabular}{|c|c|c|c|c|c|c|c|c|c|c|c|}
\hline
 $N{r,s}$   & 0010         & 0110         & 0120         & 0121          & 1110         & 1120         & 1121          & 1220         & 1221          & 1231         & $N{r,s}$   \\\hline
 0001       & $\epsilon_1$ & $\epsilon_2$ & $\epsilon_3$ & $2\epsilon_4$ & $\epsilon_5$ & $\epsilon_6$ & $2\epsilon_7$ & $\epsilon_8$ & $2\epsilon_9$ & $\epsilon_0$ & 0010       \\\hline
 0100       & 0110         & 1100         & 1110         & 1221          & 1232         & $N{r,s}$     & 1000          & 1120         & 1242          & $N{r,s}$     & 1342       \\\hline
 $\delta_1$ & $2\delta_2$  & $\delta_3$   & $2\delta_4$  & $\delta_5$    & $2\delta_6$  & 0100         & $\gamma_1$    & $\gamma_2$   & $\gamma_3$    & 1000         & $\alpha_1$ \\\hline
\end{tabular}
\end{center}

\vskip2mm
\centerline{\textit{Table 2.}}
\vskip5mm

\textit{Then the values of non-zero structure constants $N_{r,s}$ and $N_{r,-s},$ where $0\prec r\prec s,$ will be as indicated in Tables~3 and 4 below. }
\vskip5mm

\begin{center}
\begin{tabular}{|c|c|c|c|c|c|}
\hline
$r$          & \boldmath$N_{c,r}$                     & form. & $r$          & \boldmath$N_{b+2c,r}$                              & form. \\ \hline
$b+c+d$      & $\epsilon_{23}\delta_2$                & (011)   & $a$          & $\delta_{1234}\gamma_1$                            & (073)   \\ \hline
$a+b+c+d$    & $\epsilon_{56}\delta_4$                & (033)   & $a+b$        & $\delta_{1234}\gamma_2$                            & (075)   \\ \hline
$a+2b+2c+2d$ & $\epsilon_{90}\delta_5$                & (115)   & $a+b+c+d$    & $-\epsilon_{58}\delta_{1245}\gamma_2$              & (077)   \\ \hline
$r$          & \boldmath$N_{c+d,r}$                   & form. & $a+b+2c+2d$  & $\epsilon_{6780}\delta_{1256}\gamma_2$             & (079)   \\ \hline
$b$          & $\epsilon_{12}\delta_1$                & (008)   & $a+2b+2c+2d$ & $\epsilon_{90}\delta_{1256}\gamma_3$               & (081)   \\ \hline
$b+c$        & $\epsilon_{13}\delta_2$                & (009)   & $r$          & \boldmath$N_{b+2c+d,r}$                            & form. \\ \hline
$b+c+d$      & $2\epsilon_{1234}\delta_2$             & (013)   & $a$          & $\epsilon_{36}\delta_{1234}\gamma_1$               & (083)   \\
\hline
\end{tabular}
\end{center}

\centerline{\textit{Table 3.}}

\begin{center}
\begin{tabular}{|c|c|c|c|c|c|}
\hline
$r$            & \boldmath$N_{d,r}$                      & form. & $r$            & \boldmath$N_{b+2c+2d,r}$                           & form. \\ \hline
$-c-d$         & $-\epsilon_1$                           & (001)   & $-a-b-2c-2d$   & $-\epsilon_{3467}\delta_{1234}\gamma_1$            & (096)   \\ \hline
$-b-c-d$       & $-\epsilon_2$                           & (003)   & $-a-2b-2c-2d$  & $-\epsilon_{3489}\delta_{1234}\gamma_2$            & (098)   \\ \hline
$-b-2c-d$      & $-2\epsilon_3$                          & (005)   & $-a-2b-3c-2d$  & $-\epsilon_{3480}\delta_{1245}\gamma_2$            & (100)   \\ \hline
$-b-2c-2d$     & $-\epsilon_4$                           & (006)   & $-a-2b-4c-2d$  & $-\epsilon_{3480}\delta_{1256}\gamma_2$            & (102)   \\ \hline
$-a-b-c-d$     & $-\epsilon_5$                           & (015)   & $-a-3b-4c-2d$  & $-\epsilon_{3480}\delta_{1256}\gamma_3$            & (104)   \\ \hline
\end{tabular}
\end{center}

\centerline{\textit{Table 4.}}

\begin{center}
\begin{tabular}{|c|c|c|c|c|c|}
\hline
$-a-b-2c-d$    & $-2\epsilon_6$                          & (016)   & $r$            & \boldmath$N_{a,r}$                                 & form. \\ \hline
$-a-b-2c-2d$   & $-\epsilon_7$                           & (017)   & $-a-b$         & $\gamma_1$                                         & (025)   \\ \hline
$-a-2b-2c-d$   & $-2\epsilon_8$                          & (018)   & $-a-b-c$       & $\delta_{13}\gamma_1$                              & (042)   \\ \hline
$-a-2b-2c-2d$  & $-\epsilon_9$                           & (019)   & $-a-b-c-d$     & $\epsilon_{25}\delta_{13}\gamma_1$                 & (052)   \\ \hline
$-a-2b-3c-2d$  & $-\epsilon_0 $                          & (020)   & $-a-b-2c$      & $\delta_{1234}\gamma_1$                            & (074)   \\ \hline
$r$            & \boldmath$N_{c,r}$                      & form. & $-a-b-2c-d$    & $\epsilon_{36}\delta_{1234}\gamma_1$               & (084)   \\ \hline
$-c-d$         & $\epsilon_1$                            & (001)   & $-a-b-2c-2d$   & $\epsilon_{3467}\delta_{1234}\gamma_1$             & (096)   \\ \hline
$-b-c$         & $-2\delta_1$                            & (002)   & $-2a-3b-4c-2d$ & $-\alpha_1$                                        & (028)   \\
\hline
\end{tabular}
\end{center}

\vskip5mm
\centerline{\textit{Table 4, continue}}
\vskip5mm

\begin{center}
\begin{tabular}{|c|c|c|c|c|c|}
\hline
$-b-2c-d$      & $\epsilon_{13}\delta_2$                 & (010)   & $-2a-3b-4c-2d$ & $\epsilon_{90}\delta_{3456}\gamma_{13}\alpha_1$    & (110)   \\ \hline
$-a-b-c$       & $-2\delta_{13}\gamma_1$                 & (042)   & $r$            & \boldmath$N_{a+b+2c+d,r}$                          & form. \\ \hline
$-a-2b-2c$     & $\delta_{14}\gamma_2$                   & (044)   & $-a-b-2c-2d$   & $\epsilon_7$                                       & (017)   \\ \hline
$-a-2b-2c-d$   & $\epsilon_{58}\delta_{14}\gamma_2$      & (058)   & $-a-2b-2c-d$   & $2\epsilon_{68}\gamma_2$                           & (046)   \\ \hline
$-a-2b-3c-d$   & $-\epsilon_{68}\delta_{15}\gamma_2$     & (066)   & $-a-2b-3c-d$   & $\epsilon_{68}\delta_{15}\gamma_2$                 & (066)   \\ \hline
$r$            & \boldmath$N_{b+c,r}$                    & form. & $r$            & \boldmath$N_{a+b+2c+d,r}$                          & form. \\ \hline
$-a-2b-3c-2d$  & $-2\epsilon_{6780}\delta_{15}\gamma_2$  & (068)   & $-a-2b-3c-2d$  & $-\epsilon_{2680}\delta_{15}\gamma_2$              & (070)   \\ \hline
$-a-3b-4c-2d$  & $\delta_{16}\gamma_3$                   & (048)   & $-a-2b-4c-2d$  & $-\epsilon_{3680}\delta_{1256}\gamma_2$            & (092)   \\
\hline
\end{tabular}
\end{center}

\vskip5mm
\centerline{\textit{Table 4, continue.}} 
\vskip5mm

\textbf{Proof.}
According to \cite[Theorem 4.1.2]{Car72} the structure constants of a simple Lie algebra of type $\Phi$ over
$\mathbb{C}$ satisfy the following relations:

(i) $\displaystyle{N_{s,r}=-N_{r,s},\ r,s\in\Phi;}$
\bigskip

(ii) $\displaystyle{\frac{N_{r_1,r_2}}{(r_3,r_3)}=
\frac{N_{r_2,r_3}}{(r_1,r_1)}=\frac{N_{r_3,r_1}}{(r_2,r_2)},}$
\medskip

\noindent
if $\displaystyle{r_1,r_2,r_3\in\Phi}$ satisfy
$\displaystyle{r_1+r_2+r_3=0;}$
\bigskip

(iii) $\displaystyle{N_{r,s}N_{-r,-s}=-(p+1)^2,}$
\medskip

\noindent
if $r,s,r+s\in\Phi;$
\bigskip

(iv)
$\displaystyle{\frac{N_{r_1,r_2}N_{r_3,r_4}}{(r_1+r_2,r_1+r_2)}+
\frac{N_{r_2,r_3}N_{r_1,r_4}}{(r_2+r_3,r_2+r_3)}+
\frac{N_{r_3,r_1}N_{r_2,r_4}}{(r_3+r_1,r_3+r_1)}=0,}$
\medskip

\noindent
if $\displaystyle{r_1,r_2,r_3,r_4\in \Phi}$ satisfy
$\displaystyle{r_1+r_2+r_3+r_4=0}$ and if no pair opposite.
\medskip

001. We have $d + (c) + (-c-d) = 0,$ therefore (see Table~2)
$$
N_{d,c} = N_{c,-c-d} = N_{-c-d,d} = \epsilon_1. \eqno{(001)}
$$

001. We have $d + (c) + (-c-d) = 0,$ therefore (see Table~2)
$$
N_{d,c} = N_{c,-c-d} = N_{-c-d,d} = \epsilon_1. \eqno{(001)}
$$

002. We have $c + (b) + (-b-c) = 0,$ therefore (see Table~2)
$$
\frac{N_{c,b}}{1} = \frac{N_{b,-b-c}}{1} = \frac{N_{-b-c,c}}{2} = \delta_1. \eqno{(002)}
$$

003. We have $d + (b+c) + (-b-c-d)  =  0,$ therefore (see Table~2)
$$
N_{d,b+c} = N_{b+c,-b-c-d} = N_{-b-c-d,d} = \epsilon_2. \eqno{(003)}
$$

004. We have $c + (b+c) + (-b-2c) = 0,$ therefore (see Table~2)
$$
\frac{N_{c,b+c}}{2} = \frac{N_{b+c,-b-2c}}{1} = \frac{N_{-b-2c,c}}{1} = \delta_2. \eqno{(004)}
$$

005. We have $d + (b+2c) + (-b-2c-d) = 0,$ therefore (see Table~2)
$$
\frac{N_{d,b+2c}}{1} = \frac{N_{b+2c,-b-2c-d}}{1} = \frac{N_{-b-2c-d,d}}{2} = \epsilon_3. \eqno{(005)}
$$

006. We have $d + (b+2c+d) + (-b-2c-2d) = 0,$ therefore (see Table~2)
$$
\frac{N_{d,b+2c+d}}{2} = \frac{N_{b+2c+d,-b-2c-2d}}{1} = \frac{N_{-b-2c-2d,d}}{1} = \epsilon_4.\eqno{(006)}
$$

007. We have $d + (c) + b + (-b-c-d) = 0$ and $b+d\notin\ F_4,$ means
$$
S = \frac{N_{d,c} N_{b,-b-c-d}}{1}+\frac{ N_{c,b} N_{d,-b-c-d}}{1} =  0.
$$
Since
$$
N_{d,c} N_{b,-b-c-d} = |\text{Table 2}| = \epsilon_1 N_{b,-b-c-d},
$$
$$
N_{c,b} N_{d,-b-c-d} = |\text{Table 2 and formula (003)}| = \delta_1 (-\epsilon_2) = -\epsilon_2\delta_1,
$$
we have
$$
S = \epsilon_1 N_{b,-b-c-d} -\epsilon_2\delta_1  =  0
$$
and therefore
$$
N_{b,-b-c-d} = \epsilon_{12}\delta_1. \eqno{(007)}
$$

008. We have $(c+d) + b + (-b-c-d) = 0,$ therefore (see formula (007))
$$
\frac{N_{c+d,b}}{1} = \frac{N_{b,-b-c-d}}{1} = \frac{N_{-b-c-d,c+d}}{2} = \epsilon_{12}\delta_1. \eqno{(008)}
$$ 

009. We have $(c+d) + (b+c) + (-b-2c) + (-d) = 0$ and $(-b-2c) + (c+d)\notin\ F_4,$ means
$$
S = \frac{N_{c+d,b+c} N_{-b-2c,-d}}{1} + \frac{N_{b+c,-b-2c} N_{c+d,-d}}{1} =  0.
$$
Since
$$
N_{c+d,b+c}N_{-b-2c,-d} = |\text{Table 2}| = \epsilon_3{N_{c+d,b+c}},
$$
$$
N_{b+c,-b-2c} N_{c+d,-d} = |\text{Formulas (004) and (001)}| = \delta_2 (-\epsilon_1) = -\epsilon_1\delta_2,
$$
we have
$$
S = \epsilon_3 N_{c+d,b+c} - \epsilon_1\delta_2  =  0
$$
and therefore
$$
N_{c+d,b+c} = \epsilon_{13}\delta_2. \eqno{(009)}
$$

010. We have $(c+d) + (b+c) + (-b-2c-d) = 0,$ therefore (see formula (009))
$$
N_{c+d,b+c} = N_{b+c,-b-2c-d} = N_{-b-2c-d,c+d} = \epsilon_{13}\delta_2. \eqno{(010)}
$$

011. We have $c + (b+c+d) + (-b-c) + (-c-d) = 0,$ means
$$
S = \frac{N_{c,b+c+d} N_{-b-c,-c-d}}{1} + \frac{ N_{b+c+d,-b-c} N_{c,-c-d}}{1} + \frac{N_{-b-c,c} N_{b+c+d,-c-d}}{2} =  0.
$$
Since
$$
N_{c,b+c+d} N_{-b-c,-c-d} = |\text{Formula (009)}| = \epsilon_{13}\delta_2 N_{c,b+c+d},
$$
$$
N_{b+c+d,-b-c} N_{c,-c-d} = |\text{Formulas (003) and (001)}| = \epsilon_2 (\epsilon_1) = \epsilon_{12},
$$
$$
N_{-b-c,c} N_{b+c+d,-c-d} = |\text{Formulas (002) and (008)}| = 2\delta_1 (-2\epsilon_{12}\delta_1) = -4\epsilon_{12},
$$
we have
$$
S = \epsilon_{13}\delta_2 N_{c,b+c+d} + \epsilon_{12} - \frac{4\epsilon_{12}}{2}  =  \epsilon_{13}\delta_2 N_{c,b+c+d} - \epsilon_{12}  =  0,
$$
and therefore
$$
N_{c,b+c+d} = \epsilon_{23}\delta_2.\eqno{(011)}
$$

012. We have $c + (b+c+d) + (-b-2c-d) = 0,$ therefore (see formula (011))
$$
N_{c,b+c+d} = N_{b+c+d,-b-2c-d} = N_{-b-2c-d,c} = \epsilon_{23}\delta_2. \eqno{(012)}.
$$

013. We have $(c+d) + (b+c+d) + (-b-2c-d)+(-d) = 0,$ means
$$
S = \frac{N_{c+d,b+c+d} N_{-b-2c-d,-d}}{2} + \frac{N_{b+c+d,-b-2c-d}N_{c+d,-d}}{1} + \frac{N_{-b-2c-d,c+d} N_{b+c+d,-d}}{1} = 0.
$$
Since
$$
N_{c+d,b+c+d} N_{-b-2c-d,-d} = |\text{Table 2}| = 2\epsilon_4 N_{c+d,b+c+d},
$$ 
$$
N_{b+c+d,-b-2c-d} N_{c+d,-d} = |\text{Formulas (012) and (001)}| = \epsilon_{23}\delta_2 (-\epsilon_1) = -\epsilon_{123}\delta_2,
$$
$$
N_{-b-2c-d,c+d} N_{b+c+d,-d} = |\text{Formulas (010) and (003)}| =  \epsilon_{13}\delta_2 (-\epsilon_2) = -\epsilon_{123}\delta_2,
$$
we have
$$
S = \frac{2\epsilon_4 N_{c+d,b+c+d}}{2} - \epsilon_{123}\delta_2-\epsilon_{123}\delta_2 = 0 
$$
and therefore
$$
N_{c+d,b+c+d} = 2\epsilon_{1234}\delta_2. \eqno{(013)}
$$

014. We have $(c+d) + (b+c+d) + (-b-2c-2d) = 0,$ therefore (see formula (013))
$$
\frac{N_{c+d,b+c+d}}{2} = \frac{N_{b+c+d,-b-2c-2d}}{1} = \frac{N_{-b-2c-2d,c+d}}{1}  = \epsilon_{1234}\delta_2. \eqno{(014)}
$$

015. We have $d + (a+b+c) + (-a-b-c-d) = 0,$ therefore (see Table 2)
$$
N_{d,a+b+c} = N_{a+b+c,-a-b-c-d} = N_{-a-b-c-d,d}  = \epsilon_5. \eqno{(015)}
$$

016. We have $d + (a+b+2c) + (-a-b-2c-d) = 0,$ therefore  (see Table 2)
$$
\frac{N_{d,a+b+2c}}{1} = \frac{N_{a+b+2c,-a-b-2c-d}}{1} = \frac{N_{-a-b-2c-d,d}}{2} = \epsilon_6. \eqno{(016)}
$$

017. We have $d + (a+b+2c+d) + (-a-b-2c-2d) = 0,$ therefore  (see Table 2)
$$
\frac{N_{d,a+b+2c+d}}{2} = \frac{N_{a+b+2c+d,-a-b-2c-2d}}{1} = \frac{N_{-a-b-2c-2d,d}}{1}  = \epsilon_7. \eqno{(017)}
$$

018. We have $d + (a+2b+2c) + (-a-2b-2c-d) = 0,$ therefore  (see Table 2)
$$
\frac{N_{d,a+2b+2c}}{1} = \frac{N_{a+2b+2c,-a-2b-2c-d}}{1} = \frac{N_{-a-2b-2c-d,d}}{2} = \epsilon_8. \eqno{(018)}
$$

019. We have $d + (a+2b+2c+d) + (-a-2b-2c-2d) = 0,$ therefore  (see Table 2)
$$
\frac{N_{d,a+2b+2c+d}}{2} = \frac{N_{a+2b+2c+d,-a-2b-2c-2d}}{1} = \frac{N_{-a-2b-2c-2d,d}}{1} = \epsilon_9. \eqno{(019)}
$$

020. We have $d + (a+2b+3c+d) + (-a-2b-3c-2d) = 0,$ therefore  (see Table 2)
$$
N_{d,a+2b+3c+d} = N_{a+2b+3c+d,-a-2b-3c-2d} = N_{-a-2b-3c-2d,d} = \epsilon_0. \eqno{(020)}
$$

021. We have $c + (a+b) + (-a-b-c) = 0,$ therefore  (see Table 2)
$$
\frac{N_{c,a+b}}{1} = \frac{N_{a+b,-a-b-c}}{1} = \frac{N_{-a-b-c,c}}{2} = \delta_3. \eqno{(021)}
$$

022. We have $c + (a+b+c) + (-a-b-2c) = 0,$ therefore  (see Table 2)
$$
\frac{N_{c,a+b+c}}{2} = \frac{N_{a+b+c,-a-b-2c}}{1} = \frac{N_{-a-b-2c,c}}{1} = \delta_4. \eqno{(022)}
$$

023. We have $c + (a+2b+2c+d) + (-a-2b-3c-d) = 0,$ therefore  (see Table 2)
$$
N_{c,a+2b+2c+d} = N_{a+2b+2c+d,-a-2b-3c-d} = N_{-a-2b-3c-d,c} = \delta_5. \eqno{(023)}
$$

024. We have $c + (a+2b+3c+2d) + (-a-2b-4c-2d) = 0,$ therefore  (see Table 2)
$$
\frac{N_{c,a+2b+3c+2d}}{2} = \frac{N_{a+2b+3c+2d,-a-2b-4c-2d}}{1} = \frac{N_{-a-2b-4c-2d,c}}{1} = \delta_6. \eqno{(024)}
$$

025. We have $b + (a) + (-a-b) = 0,$ therefore  (see Table 2)
$$
N_{b,a} = N_{a,-a-b} = N_{-a-b,b} = \gamma_1. \eqno{(025)}
$$

026. We have $b + (a+b+2c) + (-a-2b-2c) = 0,$ therefore  (see Table 2)
$$
N_{b,a+b+2c} = N_{a+b+2c,-a-2b-2c} = N_{-a-2b-2c,b} = \gamma_2. \eqno{(026)}
$$

027. We have $b + (a+2b+4c+2d) + (-a-3b-4c-2d) = 0,$ therefore  (see Table 2)
$$
N_{b,a+2b+4c+2d} = N_{a+2b+4c+2d,-a-3b-4c-2d} = N_{-a-3b-4c-2d,b} = \gamma_3. \eqno{(027)}
$$

028. We have $a + (a+3b+4c+2d) + (-2a-3b-4c-2d) = 0,$ therefore  (see Table 2)
$$
N_{a,a+3b+4c+2d} = N_{a+3b+4c+2d,-2a-3b-4c-2d} = N_{-2a-3b-4c-2d,a} = \alpha_1. \eqno{(028)}
$$

029. We have $c + (d) + (a+b) + (-a-b-c-d) = 0$ and $d + (a+b)\notin\ F_4,$ means
$$
S = \frac{N_{c,d} N_{a+b,-a-b-c-d}}{1} + \frac{N_{a+b,c} N_{d,-a-b-c-d}}{1} =  0.
$$
Since
$$
N_{c,d}N_{a+b,-a-b-c-d} = |\text{Table 2}| = -\epsilon_1 N_{a+b,-a-b-c-d},
$$
$$
N_{a+b,c}N_{d,-a-b-c-d} = |\text{Table 2 and formula (015)}| = -\delta_3 (-\epsilon_5) = \epsilon_5\delta_3,
$$
we have
$$
S = -\epsilon_1 N_{a+b,-a-b-c-d} + \epsilon_5\delta_3 = 0,
$$
and therefore
$$
N_{a+b,-a-b-c-d} = \epsilon_{15}\delta_3. \eqno{(029)}
$$

030. We have $(a+b) + (-a-b-c-d) + (c+d) = 0,$ therefore  (see formula (029))
$$ 
\frac{N_{a+b,-a-b-c-d}}{1} = \frac{N_{-a-b-c-d,c+d}}{2} = \frac{N_{c+d,a+b}}{1}  = \epsilon_{15}\delta_3. \eqno{(030)}
$$

031. We have $(c+d) + (a+b+c) + (-d) + (-a-b-2c) = 0$ and $(a+b+c) + (-d)\notin\ F_4,$ means
$$
S = \frac{N_{c+d,a+b+c} N_{-d,-a-b-2c}}{1} + \frac{N_{-d,c+d} N_{a+b+c,-a-b-2c}}{1} =  0.
$$
Since
$$
N_{c+d,a+b+c} N_{-d,-a-b-2c} = |\text{Table 2}| = -\epsilon_6 N_{c+d,a+b+c},
$$
$$
N_{-d,c+d} N_{a+b+c,-a-b-2c} = |\text{Formulas (001) and (022)}| = \epsilon_1\delta_4,
$$
we have
$$
S = -\epsilon_6 N_{c+d,a+b+c} + \epsilon_1\delta_4 = 0,
$$
and therefore
$$
N_{c+d,a+b+c} = \epsilon_{16}\delta_4. \eqno{(031)}
$$

032. We have $(c+d) + (a+b+c) + (-a-b-2c-d) = 0,$ therefore  (see formula (031))
$$
N_{c+d,a+b+c} = N_{a+b+c,-a-b-2c-d} = N_{-a-b-2c-d,c+d} = \epsilon_{16}\delta_4. \eqno{(032)}
$$

033. We have $c + (a+b+c+d) + (-d) + (-a-b-2c) = 0$ and $(-d) + c\notin\ F_4,$ means
$$
S = \frac{N_{c,a+b+c+d} N_{-d,-a-b-2c}}{1} + \frac{N_{a+b+c+d,-d} N_{c,-a-b-2c}}{1} =  0.
$$
Since
$$
N_{c,a+b+c+d} N_{-d,-a-b-2c} = |\text{Table 2}| = -\epsilon_6 N_{c,a+b+c+d},
$$
$$
N_{a+b+c+d,-d} N_{c,-a-b-2c} =  |\text{Formulas (015) and (022)}| = -\epsilon_5(-\delta_4) = \epsilon_5\delta_4,
$$
we have
$$
S = -\epsilon_6 N_{c,a+b+c+d} + \epsilon_5\delta_4  =  0,
$$
and therefore
$$
N_{c,a+b+c+d} = \epsilon_{56}\delta_4. \eqno{(033)}
$$

034. We have $c + (a+b+c+d) + (-a-b-2c-d) = 0,$ therefore  (see formula (033))
$$
N_{c,a+b+c+d} = N_{a+b+c+d,-a-b-2c-d} = N_{-a-b-2c-d,c} = \epsilon_{56}\delta_4. \eqno{(034)}
$$

035. We have $(c+d) + (a+2b+2c) + (-c) + (-a-2b-2c-d) = 0$ and $(a+2b+2c) + (-c)\notin\ F_4,$ means
$$
S = \frac{N_{c+d,a+2b+2c} N_{-c,-a-2b-2c-d}}{1 }+ \frac{N_{-c,c+d} N_{a+2b+2c,-a-2b-2c-d}}{1} =  0.
$$
Since
$$
N_{c+d,a+2b+2c} N_{-c,-a-2b-2c-d} = |\text{Table 2}| = -\delta_5 N_{c+d,a+2b+2c},
$$  
$$
N_{-c,c+d} N_{a+2b+2c,-a-2b-2c-d} = |\text{Formulas (001) and (018)}| = -\epsilon_1 (\epsilon_8) = -\epsilon_{18},
$$ 
we have
$$
S = -\delta_5 N_{c+d,a+2b+2c} - \epsilon_{18} = 0
$$
and therefore
$$
N_{c+d,a+2b+2c}  = -\epsilon_{18}\delta_5. \eqno{(035)}
$$

036. We have $(c+d) + (a+2b+2c) + (-a-2b-3c-d) = 0,$ therefore  (see formula (035))
$$
\frac{N_{c+d,a+2b+2c}}{1} = \frac{N_{a+2b+2c,-a-2b-3c-d}}{1} = \frac{N_{-a-2b-3c-d,c+d}}{2} = -\epsilon_{18}\delta_5. \eqno{(036)}
$$

037. We have $(c+d) + (a+2b+2c+d) + (-d) + (-a-2b-3c-d) = 0,$ means
$$
S = \frac{N_{c+d,a+2b+2c+d} N_{-d,-a-2b-3c-d}}{1} + \frac{N_{a+2b+2c+d,-d} N_{c+d,-a-2b-3c-d}}{2}+
$$
$$
+\frac{N_{-d,c+d} N_{a+2b+2c+d,-a-2b-3c-d}}{1} =  0.
$$
Since
$$
N_{c+d,a+2b+2c+d} N_{-d,-a-2b-3c-d} =|\text{Table 2}| = -\epsilon_0 N_{c+d,a+2b+2c+d},
$$ 
$$
N_{a+2b+2c+d,-d} N_{c+d,-a-2b-3c-d} = |\text{Formulas (018) and (036)}| = -2\epsilon_8 (2\epsilon_{18}\delta_5) = -4\epsilon_1\delta_5,
$$
$$
N_{-d,c+d} N_{a+2b+2c+d,-a-2b-3c-d} = |\text{Formulas (001) and (023)}| = \epsilon_1 (\delta_5) = \epsilon_1\delta_5,
$$
we have
$$
S = -\epsilon_0 N_{c+d,a+2b+2c+d} - \frac{4\epsilon_1\delta_5}{2} + \epsilon_1\delta_5  =  0
$$
and therefore
$$
N_{c+d,a+2b+2c+d}  = -\epsilon_{10}\delta_5. \eqno{(037)}
$$

038. We have $(c+d) + (a+2b+2c+d) + (-a-2b-3c-2d) = 0,$ therefore  (see formula (037))
$$
N_{c+d,a+2b+2c+d} = N_{a+2b+2c+d,-a-2b-3c-2d} = N_{-a-2b-3c-2d,c+d} = -\epsilon_{10}\delta_5. \eqno{(038)}
$$

039. We have $(c+d) + (a+2b+3c+d) + (-c) + (-a-2b-3c-2d) = 0,$ means 
$$
S = \frac{N_{c+d,a+2b+3c+d} N_{-c,-a-2b-3c-2d}}{2} + \frac{N_{a+2b+3c+d,-c} N_{c+d,-a-2b-3c-2d}}{1} +
$$
$$
+ \frac{N_{-c,c+d} N_{a+2b+3c+d,-a-2b-3c-2d}}{1} =  0,
$$
Since
$$
N_{c+d,a+2b+3c+d} N_{-c,-a-2b-3c-2d} = |\text{Table 2}| = -2\delta_6 N_{c+d,a+2b+3c+d},
$$
$$
N_{a+2b+3c+d,-c} N_{c+d,-a-2b-3c-2d} = |\text{Formulas (023) and (038)}| = -\delta_5 (\epsilon_{10}\delta_5)  =  -\epsilon_{10},
$$
$$
N_{-c,c+d} N_{a+2b+3c+d,-a-2b-3c-2d} = |\text{Formulas (001) and (020)}| = -\epsilon_1 (\epsilon_0) = -\epsilon_{10},
$$
we have
$$
S = -\frac{2\delta_6 N_{c+d,a+2b+3c+d}}{2} - \epsilon_{10}-\epsilon_{10} = 0
$$
and therefore
$$
N_{c+d,a+2b+3c+d}  = -2\epsilon_{10}\delta_6. \eqno{(039)}
$$

040. We have $(c+d) + (a+2b+3c+d) + (-a-2b-4c-2d) = 0,$ therefore  (see formula (039))
$$
\frac{N_{c+d,a+2b+3c+d}}{2} = \frac{N_{a+2b+3c+d,-a-2b-4c-2d}}{1} = \frac{N_{-a-2b-4c-2d,c+d}}{1} = -\epsilon_{10}\delta_6.\eqno{(040)}
$$

41. We have $(b+c) + a + (-c) + (-a-b) = 0$ and $a+(-c)\notin F_4,$ means
$$
S = \frac{ N_{b+c,a} N_{-c,-a-b}}{1} + \frac{ N_{-c,b+c} N_{a,-a-b}}{2} =  0.
$$
Since
$$
N_{b+c,a} N_{-c,-a-b} = |\text{Table 2}| = -\delta_3 N_{b+c,a},
$$
$$
N_{-c,b+c} N_{a,-a-b} = |\text{Formulas (002) and (025)}| = 2\delta_1 (\gamma_1)=2\delta_1\gamma_1,
$$
we have
$$
S = -\delta_3 N_{b+c,a} + \frac{2\delta_1\gamma_1}{2} = 0
$$
and therefore
$$
N_{b+c,a}  = \delta_{13}\gamma_1. \eqno{(041)}
$$

042. We have $(b+c) + a + (-a-b-c) = 0,$ therefore  (see formula (041))
$$
\frac{N_{b+c,a}}{1} = \frac{N_{a,-a-b-c}}{1} = \frac{N_{-a-b-c,b+c}}{2} = \delta_{13}\gamma_1. \eqno{(042)}
$$

043. We have $(b+c)+(a+b+c)+(-b)+(-a-b-2c) = 0$ and $(a+b+c) + (-b)\notin F_4,$ means
$$
S = \frac{N_{b+c,a+b+c} N_{-b,-a-b-2c}}{2} + \frac{N_{-b,b+c} N_{a+b+c,-a-b-2c}}{1} =  0.
$$
Since
$$
N_{b+c,a+b+c} N_{-b,-a-b-2c} = |\text{Table 2}| = -\gamma_2 N_{b+c,a+b+c},
$$
$$
N_{-b,b+c} N_{a+b+c,-a-b-2c} = |\text{Formulas (002) and (022)}| = -\delta_1 (\delta_4) = -\delta_{14}.
$$
we have
$$
S = -\frac{\gamma_2 N_{b+c,a+b+c}}{2} - \delta_{14}  =  0
$$
and therefore
$$
N_{b+c,a+b+c}  =  -2\delta_{14}\gamma_2. \eqno{(043)}
$$

044. We have $(b+c) + (a+b+c) + (-a-2b-2c) = 0,$ therefore  (see formula (043))
$$
\frac{N_{b+c,a+b+c}}{2} = \frac{N_{a+b+c,-a-2b-2c}}{1} = \frac{N_{-a-2b-2c,b+c}}{1} = -\delta_{14}\gamma_2. \eqno{(044)}
$$

045. We have $b + (a+b+2c+d) + (-d) + (-a-2b-2c) = 0$ and $(-d) + b \notin F_4,$ means
$$
S = \frac{ N_{b,a+b+2c+d} N_{-d,-a-2b-2c}}{1} + \frac{N_{a+b+2c+d,-d} N_{b,-a-2b-2c}}{2} =  0.
$$
Since
$$
N_{b,a+b+2c+d} N_{-d,-a-2b-2c} = |\text{Table 2}| = -\epsilon_8 N_{b,a+b+2c+d},
$$
$$
N_{a+b+2c+d,-d} N_{b,-a-2b-2c} = |\text{Formulas (016) and (026)}| = -2\epsilon_6 (-\gamma_2) = 2\epsilon_6\gamma_2,
$$
we have
$$
S = -\epsilon_8 N_{b,a+b+2c+d} + \frac{2\epsilon_6\gamma_2}{2}  =  0
$$
and therefore
$$
N_{b,a+b+2c+d} = \epsilon_{68}\gamma_2. \eqno{(045)}
$$

046. We have $b + (a+b+2c+d) + (-a-2b-2c-d) = 0,$ therefore  (see formula (045))
$$
\frac{N_{b,a+b+2c+d}}{1} = \frac{N_{a+b+2c+d,-a-2b-2c-d}}{2} = \frac{N_{-a-2b-2c-d,b}}{1} = \epsilon_{68}\gamma_2. \eqno{(046)}
$$

047. We have $(b+c)+(a+2b+3c+2d)+(-b)+(-a-2b-4c-2d) = 0$ and $(a+2b+3c+2d) + (-b)\notin F_4,$ means
$$
S = \frac{N_{b+c,a+2b+3c+2d} N_{-b,-a-2b-4c-2d}}{2} + \frac{N_{-b,b+c} N_{a+2b+3c+2d,-a-2b-4c-2d}}{1} =  0.
$$
Since
$$
N_{b+c,a+2b+3c+2d} N_{-b,-a-2b-4c-2d} = |\text{Table 2}| = -\gamma_3 N_{b+c,a+2b+3c+2d},
$$
$$
N_{-b,b+c} N_{a+2b+3c+2d,-a-2b-4c-2d} = |\text{Formulas (002) and (024)}| = -\delta_1 (\delta_6) = -\delta_{16}.
$$
we have
$$
S = -\frac{\gamma_3 N_{b+c,a+2b+3c+2d}}{2} - \delta_{16}  =  0
$$
and therefore
$$
N_{b+c,a+2b+3c+2d} =  -2\delta_{16}\gamma_3. \eqno{(047)}
$$

048. We have $(b+c) + (a+2b+3c+2d) + (-a-3b-4c-2d) = 0,$ therefore  (see formula (047))
$$
\frac{N_{b+c,a+2b+3c+2d}}{2} = \frac{N_{a+2b+3c+2d,-a-3b-4c-2d}}{1} = \frac{N_{-a-3b-4c-2d,b+c}}{1} = -\delta_{16}\gamma_3.\eqno{(048)}
$$

049. We have $(a+b) + (a+2b+4c+2d) + (-a) + (-a-3b-4c-2d) = 0$ and $(a+2b+4c+2d) + (-a) \notin F_4,$ means
$$
S = \frac{N_{a+b,a+2d+4c+2b} N_{-a,-a-3b-4c-2d}}{2} + \frac{N_{-a,a+b}N_{a+2b+4c+2d,-a-3b-4c-2d}}{2} =  0.
$$
Since
$$
N_{a+b,a+2b+4c+2d}N_{-a,-a-3b-4c-2d} = |\text{Table 2}| = -\alpha_1 N_{a+b,a+2b+4c+2d},
$$
$$
N_{-a,a+b} N_{a+2b+4c+2d,-a-3b-4c-2d} = |\text{Formulas (025) and (027)}| = -\gamma_1 (\gamma_3) = -\gamma_{13},
$$
we have
$$
S = -\frac{\alpha_1 N_{a+b,a+2b+4c+2d}}{2} - \frac{\gamma_{13}}{2} = 0
$$
and therefore
$$
N_{a+b,a+2b+4c+2d} = -\gamma_{13}\alpha_1. \eqno{(049)}
$$

050. We have $(a+b) + (a+2b+4c+2d) + (-2a-3b-4c-2d) = 0,$ therefore  (see formula (049))
$$
N_{a+b,a+2b+4c+2b} = N_{a+2b+4c+2b,-2a-3b-4c-2b} = N_{-2a-3b-4c-2b,a+b} = -\gamma_{13}\alpha_1. \eqno{(050)}
$$

051. We have $(b+c+d) + a + (-d) + (-a-b-c) = 0$ and $a + (-d) \notin F_4,$ means
$$
S = \frac{N_{b+c+d,a} N_{-d,-a-b-c}}{1} + \frac{N_{-d,b+c+d} N_{a,-a-b-c}}{1} =  0.
$$
Since
$$
N_{b+c+d,a} N_{-d,-a-b-c} = |\text{Table 2}| = -\epsilon_5 N_{b+c+d,a},
$$
$$
N_{-d,b+c+d} N_{a,-a-b-c} = |\text{Formulas (003) and (042)}| = \epsilon_2 (\delta_{13}\gamma_1),
$$
we have
$$
S = -\epsilon_5 N_{b+c+d,a} + \epsilon_2\delta_{13}\gamma_1 = 0
$$
and therefore
$$
N_{b+c+d,a} = \epsilon_{25}\delta_{13}\gamma_1. \eqno{(051)}
$$ 

052. We have $(b+c+d) + a + (-a-b-c-d) = 0,$ therefore  (see formula 051))
$$
\frac{N_{b+c+d,a}}{1} = \frac{N_{a,-a-b-c-d}}{1} = \frac{N_{-a-b-c-d,b+c+d}}{2} = \epsilon_{25}\delta_{13}\gamma_1. \eqno{(052)}
$$

053. We have $(c+d) + (a+b+c+d) + (-d)+(-a-b-2c-d) = 0,$ means
$$
S = \frac{N_{c+d,a+b+c+d} N_{-d,-a-b-2c-d}}{2} + \frac{N_{a+b+c+d,-d} N_{c+d,-a-b-2c-d}}{1} + 
$$
$$
+\frac{N_{-d,c+d} N_{a+b+c+d,-a-b-2c-d}}{1} =  0.
$$
Since
$$
N_{c+d,a+b+c+d} N_{-d,-a-b-2c-d} = |\text{Table 2}| = -2\epsilon_7 N_{c+d,a+b+c+d},
$$
$$
N_{a+b+c+d,-d} N_{c+d,-a-b-2c-d} = |\text{Formulas (015) and (032)}| = -\epsilon_5 (-\epsilon_{16}\delta_4) = \epsilon_{156}\delta_4,
$$
$$
N_{-d,c+d} N_{a+b+c+d,-a-b-2c-d} = |\text{Formulas (001) and (034)}| = \epsilon_1 (\epsilon_{56}\delta_4) = \epsilon_{156}\delta_4,
$$
we have
$$
S = -\frac{2\epsilon_7 N_{c+d,a+b+c+d}}{2} + \epsilon_{156}\delta_4+\epsilon_{156}\delta_4 = 0
$$
and therefore
$$
N_{c+d,a+b+c+d}  = 2\epsilon_{1567}\delta_4. \eqno{(053)}
$$

054. We have $(c+d) + (a+b+c+d) + (-a-b-2c-2d) = 0,$ therefore  (see formula (053))
$$
\frac{N_{c+d,a+b+c+d}}{2} = \frac{N_{a+b+c+d,-a-b-2c-2d}}{1} = \frac{N_{-a-b-2c-2d,c+d}}{1} = \epsilon_{1567}\delta_4. \eqno{(054)}
$$

055. We have $b + (a+b+2c+2d) + (-d) + (-a-2b-2c-d) = 0$ and $(-d) + b \notin F_4,$ means
$$
S = \frac{N_{b,a+b+2c+2d} N_{-d,-a-2b-2c-d}}{2} + \frac{ N_{a+b+2c+2d,-d} N_{b,-a-2b-2c-d}}{1} =  0.
$$
Since
$$
N_{b,a+b+2c+2d} N_{-d,-a-2b-2c-d} = |\text{Table 2}| = -2\epsilon_9 N_{b,a+b+2c+2d},
$$
$$
N_{a+b+2c+2d,-d} N_{b,-a-2b-2c-d} = |\text{Formulas (017) and (046)}| = -\epsilon_7 (-\epsilon_{68}\gamma_2)  =  \epsilon_{678}\gamma_2,
$$
we have
$$ 
S = -\frac{2\epsilon_9 N_{b,a+b+2c+2d}}{2} + \epsilon_{678}\gamma_2 = 0
$$
and therefore 
$$
N_{b,a+b+2c+2d}  = \epsilon_{6789}\gamma_2. \eqno{(055)}
$$

056. We have $b + (a+b+2c+2d) + (-a-2b-2c-2d) = 0,$ therefore  (see formula (055))
$$
N_{b,a+b+2c+2d} = N_{a+b+2c+2d,-a-2b-2c-2d} = N_{-a-2b-2c-2d,b} = \epsilon_{6789}\gamma_2. \eqno{(056)}
$$

057. We have $(b+c) + (a+b+c+d) + (-d) + (-a-2b-2c) = 0$ and $(-d) + (b+c) \notin F_4,$ means
$$
S = \frac{N_{b+c,a+b+c+d} N_{-d,-a-2b-2c}}{1} + \frac{N_{a+b+c+d,-d}N_{b+c,-a-2b-2c}}{1} =  0.
$$
Since
$$
N_{b+c,a+b+c+d} N_{-d,-a-2b-2c} = |\text{Table 2}| = -\epsilon_8 N_{b+c,a+b+c+d},
$$
$$
N_{a+b+c+d,-d} N_{b+c,-a-2b-2c} = |\text{Formulas (015) and (044)}| = -\epsilon_5 (\delta_{14}\gamma_2) = -\epsilon_5\delta_{14}\gamma_2,
$$
we have
$$
S = -\epsilon_8 N_{b+c,a+b+c+d} - \epsilon_5\delta_{14}\gamma_2 = 0
$$
and therefore 
$$
N_{b+c,a+b+c+d}  = -\epsilon_{58}\delta_{14}\gamma_2. \eqno{(057)}
$$

058. We have $(b+c) + (a+b+c+d) + (-a-2b-2c-d) = 0,$ therefore  (see formula (057))
$$
N_{b+c,a+b+c+d} = N_{a+b+c+d,-a-2b-2c-d} = N_{-a-2b-2c-d,b+c} = -\epsilon_{58}\delta_{14}\gamma_2. \eqno{(058)}
$$

059. We have $(b+c+d) + (a+b+c) + (-d) + (-a-2b-2c) = 0$ and $(a+b+c) + (-d) \notin F_4,$  means
$$
S = \frac{N_{b+c+d,a+b+c} N_{-d,-a-2b-2c}}{1} + \frac{N_{-d,b+c+d} N_{a+b+c,-a-2b-2c}}{1} =  0.
$$
Since
$$
N_{b+c+d,a+b+c} N_{-d,-a-2b-2c} = |\text{Table 2}| = -\epsilon_8 N_{b+c+d,a+b+c},
$$
$$
N_{-d,b+c+d} N_{a+b+c,-a-2b-2c} = |\text{Formulas (003) and (044)}| = \epsilon_2 (-\delta_{14}\gamma_2) = -\epsilon_2\delta_{14}\gamma_2,
$$
we have
$$
S = -\epsilon_8 N_{b+c+d,a+b+c} - \epsilon_2\delta_{14}\gamma_2 = 0
$$
and therefore
$$
N_{b+c+d,a+b+c} = -\epsilon_{28}\delta_{14}\gamma_2. \eqno{(059)}
$$

060. We have $(b+c+d) + (a+b+c) + (-a-2b-2c-d) = 0,$ therefore  (see formula (059))
$$
N_{b+c+d,a+b+c} = N_{a+b+c,-a-2b-2c-d} = N_{-a-2b-2c-d,b+c+d} = -\epsilon_{28}\delta_{14}\gamma_2. \eqno{(060)}
$$

061. We have $(b+c+d) + (a+b+c+d) + (-d) + (-a-2b-2c-d) = 0,$ means
$$
S = \frac{N_{b+c+d,a+b+c+d} N_{-d,-a-2b-2c-d}}{2} + \frac{N_{a+b+c+d,-d} N_{b+c+d,-a-2b-2c-d}}{1} +
$$
$$
+ \frac{N_{-d,b+c+d} N_{a+b+c+d,-a-2b-2c-d}}{1} =  0.
$$
Since
$$
N_{b+c+d,a+b+c+d} N_{-d,-a-2b-2c-d} = |\text{Table 2}| = -2\epsilon_9 N_{b+c+d,a+b+c+d},
$$
$$
N_{a+b+c+d,-d} N_{b+c+d,-a-2b-2c-d} = |\text{Formulas (015) and (060)}| = -\epsilon_5 (\epsilon_{28}\delta_{14}\gamma_2) = -\epsilon_{258}\delta_{14}\gamma_2,
$$
$$
N_{-d,b+c+d} N_{a+b+c+d,-a-2b-2c-d} = |\text{Formulas (003) and (058)}|  = \epsilon_2 (-\epsilon_{58}\delta_{14}\gamma_2) = -\epsilon_{258}\delta_{14}\gamma_2,
$$
we have
$$
S = -\frac{2\epsilon_9 N_{b+c+d,a+b+c+d}}{2} - \epsilon_{258}\delta_{14}\gamma_2 - \epsilon_{258}\delta_{14}\gamma_2 = 0
$$
and therefore
$$
N_{b+c+d,a+b+c+d}  =  -2\epsilon_{2589}\delta_{14}\gamma_2. \eqno{(061)}
$$

062. We have $(b+c+d) + (a+b+c+d) + (-a-2b-2c-2d) = 0,$ therefore  (see formula (061))
$$
\frac{N_{b+c+d,a+b+c+d}}{2} = \frac{N_{a+b+c+d,-a-2b-2c-2d}}{1} = \frac{N_{-a-2b-2c-2d,b+c+d}}{1} = -\epsilon_{2589}\delta_{14}\gamma_2. \eqno{(062)}
$$

063. We have $(b+c+d) + (a+b+2c) + (-c) + (-a-2b-2c-d) = 0$ and $(-c) + (b+c+d) \notin F_4,$ means
$$
S = \frac{N_{b+c+d,a+b+2c} N_{-c,-a-2b-2c-d}}{1} + \frac{N_{a+b+2c,-c} N_{b+c+d,-a-2b-2c-d}}{1} =  0.
$$
Since
$$
N_{b+c+d,a+b+2c} N_{-c,-a-2b-2c-d} = |\text{Table 2}| = -\delta_5 N_{b+c+d,a+b+2c},
$$
$$
N_{a+b+2c,-c} N_{b+c+d,-a-2b-2c-d} = |\text{Formulas (022) and (060)}| = -\delta_4 (\epsilon_{28}\delta_{14}\gamma_2) = -\epsilon_{28}\delta_1\gamma_2,
$$
we have
$$
S = -\delta_5N_{b+c+d,a+b+2c} - \epsilon_{28}\delta_1\gamma_2 = 0
$$
and therefore
$$
N_{b+c+d,a+b+2c} = -\epsilon_{28}\delta_{15}\gamma_2. \eqno{(063)}
$$

064. We have $(b+c+d) + (a+b+2c) + (-a-2b-3c-d) = 0,$ therefore  (see formula (063))
$$
\frac{N_{b+c+d,a+b+2c}}{1} = \frac{N_{a+b+2c,-a-2b-3c-d}}{1} = \frac{N_{-a-2b-3c-d,b+c+d}}{2} = -\epsilon_{28}\delta_{15}\gamma_2. \eqno{(064)}
$$

065. We have $(b+c) + (a+b+2c+d) + (-c) + (-a-2b-2c-d) = 0,$ means
$$
S = \frac{N_{b+c,a+b+2c+d} N_{-c,-a-2b-2c-d}}{1}+\frac{N_{a+b+2c+d,-c}N_{b+c,-a-2b-2c-d}}{1}+
$$
$$
+\frac{N_{-c,b+c} N_{a+b+2c+d,-a-2b-2c-d}}{2} =  0.
$$
Since
$$
N_{b+c,a+b+2c+d} N_{-c,-a-2b-2c-d} = |\text{Table 2}| = -\delta_5 N_{b+c,a+b+2c+2d},
$$
$$
N_{a+b+2c+d,-c} N_{b+c,-a-2b-2c-d} = |\text{Formulas (034) and (058)}| = -\epsilon_{56}\delta_4 (\epsilon_{58}\delta_{14}\gamma_2) = -\epsilon_{68}\delta_1\gamma_2,
$$
$$
N_{-c,b+c} N_{a+b+2c+d,-a-2b-2c-d} = |\text{Formulas (002) and (046)}| = 2\delta_1 (2\epsilon_{68}\gamma_2) = 4\epsilon_{68}\delta_1\gamma_2,
$$
we have
$$
S = -\delta_5 N_{b+c,a+b+2c+d} - \epsilon_{68}\delta_1\gamma_2 + \frac{4\epsilon_{68}\delta_1\gamma_2}{2} = 0
$$
and therefore
$$
N_{b+c,a+b+2c+d} = \epsilon_{68}\delta_{15}\gamma_2. \eqno{(065)}
$$

066. We have $(b+c) + (a+b+2c+d) + (-a-2b-3c-d) = 0,$ therefore  (see formula (065))
$$
N_{b+c,a+b+2c+d} = N_{a+b+2c+d,-a-2b-3c-d} = N_{-a-2b-3c-d,b+c} = \epsilon_{68}\delta_{15}\gamma_2. \eqno{(066)}
$$

067. We have $(b+c) + (a+b+2c+2d) + (-d) + (-a-2b-3c-d) = 0$ and $(-d) + (b+c) \notin F_4,$ means
$$
S = \frac{N_{b+c,a+b+2c+2d} N_{-d,-a-2b-3c-d}}{1} + \frac{N_{a+b+2c+2d,-d} N_{b+c,-a-2b-3c-d}}{1} = 0.
$$
Since
$$
N_{b+c,a+b+2c+2d} N_{-d,-a-2b-3c-d} = |\text{Table 2}| = -\epsilon_0 N_{b+c,a+b+2c+2d},
$$
$$
N_{a+b+2c+2d,-d} N_{b+c,-a-2b-3c-d} = |\text{Formulas (017) and (066)}| = -\epsilon_7 (-\epsilon_{68}\delta_{15}\gamma_2) = \epsilon_{678}\delta_{15}\gamma_2,
$$
we have
$$
S = -\epsilon_0 N_{b+c,a+b+2c+2d} + \epsilon_{678}\delta_{15}\gamma_2 = 0,
$$
and therefore
$$
N_{b+c,a+b+2c+2d} = \epsilon_{6780}\delta_{15}\gamma_2. \eqno{(067)}
$$

068. We have $(b+c) + (a+b+2c+2d) + (-a-2b-3c-2d) = 0,$ therefore  (see formula (067))
$$
\frac{N_{b+c,a+b+2c+2d}}{1} = \frac{N_{a+b+2c+2d,-a-2b-3c-2d}}{1} = \frac{N_{-a-2b-3c-2d,b+c}}{2} =  \epsilon_{6780}\delta_{15}\gamma_2. \eqno{(068)}
$$

069. We have $(b+c+d) + (a+b+2c+d) + (-d) + (-a-2b-3c-d) = 0,$ means
$$
S = \frac{N_{b+c+d,a+b+2c+d} N_{-d,-a-2b-3c-d}}{1} + \frac{N_{a+b+2c+d,-d} N_{b+c+d,-a-2b-3c-d}}{2} +
$$
$$
+ \frac{N_{-d,b+c+d}N_{a+b+2c+d,-a-2b-3c-d}}{1} =  0.
$$
Since
$$
N_{b+c+d,a+b+2c+d} N_{-d,-a-2b-3c-d} = |\text{Table 2}| = -\epsilon_0 N_{b+c+d,a+b+2c+d},
$$
$$
N_{a+b+2c+d,-d} N_{b+c+d,-a-2b-3c-d} = |\text{Formulas (016) and (064)}| = -2\epsilon_6 (2\epsilon_{28}\delta_{15}\gamma_2) = 
$$
$$
=-4\epsilon_{268}\delta_{15}\gamma_2,
$$
$$
N_{-d,b+c+d} N_{a+b+2c+d,-a-2b-3c-d} = |\text{Formulas (003) and (066)}| = \epsilon_2 (\epsilon_{68}\delta_{15}\gamma_2) = \epsilon_{268}\delta_{15}\gamma_2,
$$
we have
$$
S = -\epsilon_0 N_{b+c+d,a+b+2c+d} - \frac{4\epsilon_{268}\delta_{15}\gamma_2}{2} + \epsilon_{268}\delta_{15}\gamma_2 = 0
$$
and therefore
$$
N_{b+c+d,a+b+2c+d} = -\epsilon_{2680}\delta_{15}\gamma_2. \eqno{(069)}
$$

070. We have $(b+c+d) + (a+b+2c+d) + (-a-2b-3c-2d) = 0,$ therefore  (see formula (069))
$$
N_{b+c+d,a+b+2c+d} = N_{a+b+2c+d,-a-2b-3c-2d} = N_{-a-2b-3c-2d,b+c+d} = -\epsilon_{2680}\delta_{15}\gamma_2. \eqno{(070)}
$$

071. We have $(b+c+d) + (a+2b+3c+d) + (-b) + (-a-2b-4c-2d) = 0$ and $(a+2b+3c+d) + (-b) \notin F_4,$ means
$$
S = \frac{N_{b+c+d,a+2b+3c+d} N_{-b,-a-2b-4c-2d}}{2} + \frac{N_{-b,b+c+d} N_{a+2b+3c+d,-a-2b-4c-2d}}{1} =  0.
$$
Since
$$
N_{b+c+d,a+2b+3c+d} N_{-b,-a-2b-4c-2d} = |\text{Table 2}| = -\gamma_3 N_{b+c+d,a+2b+3c+d},
$$
$$
N_{-b,b+c+d} N_{a+2b+3c+d,-a-2b-4c-2d} = |\text{Formulas (007) and (040)}| = -\epsilon_{12}\delta_1 (-\epsilon_{10}\delta_6) = \epsilon_{20}\delta_{16},
$$
we have
$$
S = -\frac{\gamma_3 N_{b+c+d,a+2b+3c+d}}{2} + \epsilon_{20}\delta_{16} = 0
$$
and therefore
$$
N_{b+c+d,a+2b+3c+d} = 2\epsilon_{20}\delta_{16}\gamma_3. \eqno{(071)}
$$

072. We have $(b+c+d) + (a+2b+3c+d) + (-a-3b-4c-2d) = 0,$ therefore  (see formula (071))
$$
\frac{N_{b+c+d,a+2b+3c+d}}{2} = \frac{N_{a+2b+3c+d,-a-3b-4c-2d}}{1} = \frac{N_{-a-3b-4c-2d,b+c+d}}{1} = \epsilon_{20}\delta_{16}\gamma_3.\eqno{(072)}
$$

073. We have $(b+2c) + a+ (-c) + (-a-b-c) = 0$ and $a + (-c) \notin F_4,$ means
$$
S = \frac{N_{b+2c,a} N_{-c,-a-b-c}}{2} + \frac{N_{-c,b+2c} N_{a,-a-b-c}}{1} =  0.
$$
Since
$$
N_{b+2c,a} N_{-c,-a-b-c} = |\text{Table 2}| = -2\delta_4 N_{b+2c,a},
$$
$$
N_{-c,b+2c} N_{a,-a-b-c} = |\text{Formulas (004) and (042)}| = \delta_2 (\delta_{13}\gamma_1) = \delta_{123}\gamma_1,
$$
we have
$$
S = -\frac{2\delta_4 N_{b+2c,a}}{2} + \delta_{123}\gamma_1 = 0
$$
and therefore
$$
N_{b+2c,a} = \delta_{1234}\gamma_1. \eqno{(073)}
$$

074. We have $(b+2c) + a + (-a-b-2c) = 0,$ therefore  (see formula (073))
$$
N_{b+2c,a} = N_{a,-a-b-2c} = N_{-a-b-2c,b+2c} = \delta_{1234}\gamma_1. \eqno{(074)}
$$

075. We have $(b+2c) + (a+b) + (-b) + (-a-b-2c) = 0$ and  $(-b) + (b+2c) \notin F_4,$ means
$$
S = \frac{N_{b+2c,a+b} N_{-b,-a-b-2c}}{2} + \frac{N_{a+b,-b} N_{b+2c,-a-b-2c}}{2} =  0.
$$
Since
$$
N_{b+2c,a+b} N_{-b,-a-b-2c} = |\text{Table 2}| = -\gamma_2 N_{b+2c,a+b},
$$
$$
N_{a+b,-b} N_{b+2c,-a-b-2c} = |\text{Formulas (025) and (074)}| = -\gamma_1 (-\delta_{1234}\gamma_1) = \delta_{1234},
$$
we have
$$
S = -\frac{\gamma_2 N_{b+2c,a+b}}{2} + \frac{\delta_{1234}}{2} = 0
$$
and therefore
$$
N_{b+2c,a+b} = \delta_{1234}\gamma_2. \eqno{(075)}
$$

076. We have $(b+2c) + (a+b) + (-a-2b-2c) = 0,$ therefore  (see formula (075))
$$
N_{b+2c,a+b} = N_{a+b,-a-2b-2c} = N_{-a-2b-2c,b+2c} =  \delta_{1234}\gamma_2. \eqno{(076)}
$$

077. We have $(b+2c) + (a+b+c+d) + (-c) + (-a-2b-2c-d) = 0$ and $(a+b+c+d) + (-c) \notin F_4,$ means
$$
S = \frac{N_{b+2c,a+b+c+d} N_{-c,-a-2b-2c-d}}{1} + \frac{N_{-c,b+2c} N_{a+b+c+d,-a-2b-2c-d}}{1} =  0.
$$
Since
$$
N_{b+2c,a+b+c+d} N_{-c,-a-2b-2c-d} = |\text{Table 2}| = -\delta_5 N_{b+2c,a+b+c+d},
$$
$$
N_{-c,b+2c} N_{a+b+c+d,-a-2b-2c-d} = |\text{Formulas (004) and (058)}| = \delta_2 (-\epsilon_{58}\delta_{14}\gamma_2) = -\epsilon_{58}\delta_{124}\gamma_2,
$$
we have
$$
S = -\delta_5 N_{b+2c,a+b+c+d} - \epsilon_{58}\delta_{124}\gamma_2 = 0
$$
and therefore
$$
N_{b+2c,a+b+c+d} = -\epsilon_{58}\delta_{1245}\gamma_2. \eqno{(077)}
$$

078. We have $(b+2c) + (a+b+c+d) + (-a-2b-3c-d) = 0,$ therefore  (see formula (077))
$$
\frac{N_{b+2c,a+b+c+d}}{1} = \frac{N_{a+b+c+d,-a-2b-3c-d}}{2} = \frac{N_{-a-2b-3c-d,b+2c}}{1} = -\epsilon_{58}\delta_{1245}\gamma_2. \eqno{(078)}
$$

079. We have $(b+2c) + (a+b+2c+2d) + (-c) + (-a-2b-3c-2d) = 0$ and $(a+b+2c+2d) + (-c) \notin F_4,$ means
$$
S = \frac{N_{b+2c,a+b+2c+2d} N_{-c,-a-2b-3c-2d}}{2} + \frac{N_{-c,b+2c} N_{a+b+2c+2d,-a-2b-3c-2d}}{1} = 0.
$$
Since
$$
N_{b+2c,a+b+2c+2d} N_{-c,-a-2b-3c-2d} = |\text{Table 2}| = -2\delta_6 N_{b+2c,a+b+2c+2d},
$$
$$
N_{-c,b+2c} N_{a+b+2c+2d,-a-2b-3c-2d} = |\text{Formulas (004) and (068)}| = \delta_2 (\epsilon_{6780}\delta_{15}\gamma_2) = \epsilon_{6780}\delta_{125}\gamma_2,
$$
we have
$$
S = -\frac{2\delta_6 N_{b+2c,a+b+2c+2d}}{2} + \epsilon_{6780}\delta_{125}\gamma_2 = 0
$$
and therefore
$$
N_{b+2c,a+b+2c+2d} = \epsilon_{6780}\delta_{1256}\gamma_2. \eqno{(079)}
$$

080. We have $(b+2c) + (a+b+2c+2d) + (-a-2b-4c-2d) = 0,$ therefore (see formula (079))
$$
N_{b+2c,a+b+2c+2d} = N_{a+b+2c+2d,-a-2b-4c-2d} = N_{-a-2b-4c-2d,b+2c} =  \epsilon_{6780}\delta_{1256}\gamma_2. \eqno{(080)}
$$

081. We have $(b+2c) + (a+2b+2c+2d) + (-b) + (-a-2b-4c-2d) = 0$ and $(-b) + (b+2c) \notin F_4,$ means
$$
S = \frac{N_{b+2c,a+2b+2c+2d} N_{-b,-a-2b-4c-2d}}{2} + \frac{N_{a+2b+2c+2d,-b} N_{b+2c,-a-2b-4c-2d}}{2} = 0.
$$
Since
$$
N_{b+2c,a+2b+2c+2d} N_{-b,-a-2b-4c-2d} = |\text{Table 2}| = -\gamma_3 N_{b+2c,a+2b+2c+2d},
$$
$$
N_{a+2b+2c+2d,-b} N_{b+2c,-a-2b-4c-2d} = |\text{Formulas (056) and (080)}| = -\epsilon_{6789}\gamma_2 (-\epsilon_{6780}\delta_{1256}\gamma_2) = 
$$
$$
 = \epsilon_{90}\delta_{1256},
$$
we have
$$
S = -\frac{\gamma_3 N_{b+2c,a+2b+2c+2d}}{2} + \frac{\epsilon_{90}\delta_{1256}}{2} = 0
$$
and therefore
$$
N_{b+2c,a+2b+2c+2d}  = \epsilon_{90}\delta_{1256}\gamma_3. \eqno{(081)}
$$

082. We have $(b+2c) + (a+2b+2c+2d) + (-a-3b-4c-2d) = 0,$ therefore  (see formula (081))
$$
N_{b+2c,a+2b+2c+2d} = N_{a+2b+2c+2d,-a-3b-4c-2d} = N_{-a-3b-4c-2d,b+2c} = \epsilon_{90}\delta_{1256}\gamma_3. \eqno{(082)}
$$

083. We have $(b+2c+d) + a + (-d) + (-a-b-2c) = 0$ and $a + (-d) \notin F_4,$ means
$$
S = \frac{N_{b+2c+d,a} N_{-d,-a-b-2c}}{1} + \frac{N_{-d,b+2c+d} N_{a,-a-b-2c}}{2} = 0.
$$
Since
$$
N_{b+2c+d,a} N_{-d,-a-b-2c} = |\text{Table 2}| = -\epsilon_6 N_{b+2c+d,a},
$$
$$
N_{-d,b+2c+d} N_{a,-a-b-2c} = |\text{Formulas (005) and (074)}| = 2\epsilon_3 (\delta_{1234}\gamma_1) = 2\epsilon_3\delta_{1234}\gamma_1,
$$
we have
$$
S = -\epsilon_6 N_{b+2c+d,a} + \frac{2\epsilon_3\delta_{1234}\gamma_1}{2} = 0
$$
and therefore
$$
N_{b+2c+d,a} = \epsilon_{36}\delta_{1234}\gamma_1. \eqno{(083)}
$$

084. We have $(b+2c+d) + a + (-a-b-2c-d) = 0,$ therefore  (see formula (083))
$$
\frac{N_{b+2c+d,a}}{1} = \frac{N_{a,-a-b-2c-d}}{1} = \frac{N_{-a-b-2c-d,b+2c+d}}{2} = \epsilon_{36}\delta_{1234}\gamma_1. \eqno{(084)}
$$

085. We have $(b+2c+d) + (a+b) + (-d) + (-a-2b-2c) = 0$ and $(a+b) + (-d) \notin F_4,$ means
$$
S = \frac{N_{b+2c+d,a+b} N_{-d,-a-2b-2c}}{1} + \frac{N_{-d,b+2c+d} N_{a+b,-a-2b-2c}}{2} =  0.
$$
Since
$$
N_{b+2c+d,a+b} N_{-d,-a-2b-2c} = |\text{Table 2}| = -\epsilon_8 N_{b+2c+d,a+b},
$$
$$
N_{-d,b+2c+d} N_{a+b,-a-2b-2c} = |\text{Formulas (005) and (076)}| = 2\epsilon_3 (\delta_{1234}\gamma_2) = 2\epsilon_3\delta_{1234}\gamma_2,
$$
we have
$$
S = -\epsilon_8 N_{b+2c+d,a+b} + \frac{2\epsilon_3\delta_{1234}\gamma_2}{2} = 0
$$
and therefore
$$
N_{b+2c+d,a+b} =  \epsilon_{38}\delta_{1234}\gamma_2. \eqno{(085)}
$$

086. We have $(b+2c+d) + (a+b) + (-a-2b-2c-d) = 0,$ therefore  (see formula (085))
$$
\frac{N_{b+2c+d,a+b}}{1} = \frac{N_{a+b,-a-2b-2c-d}}{1} = \frac{N_{-a-2b-2c-d,b+2c+d}}{2} = \epsilon_{38}\delta_{1234}\gamma_2. \eqno{(086)}
$$

087. We have $(b+2c+d) + (a+b+c) + (-c) + (-a-2b-2c-d) = 0,$ means
$$
S = \frac{N_{b+2c+d,a+b+c} N_{-c,-a-2b-2c-d}}{1}+\frac{N_{a+b+c,-c} N_{b+2c+d,-a-2b-2c-d}}{2}+
$$
$$
+\frac{N_{-c,b+2c+d} N_{a+b+c,-a-2b-2c-d}}{1} =  0.
$$
Since
$$
N_{b+2c+d,a+b+c} N_{-c,-a-2b-2c-d} = |\text{Table 2}|  = -\delta_5 N_{b+2c+d,a+b+c},
$$
$$
N_{a+b+c,-c} N_{b+2c+d,-a-2b-2c-d} = |\text{Formulas (021) and (086)}| = -2\delta_3 (-2\epsilon_{38}\delta_{1234}\gamma_2) = 
$$
$$
= 4\epsilon_{38}\delta_{124}\gamma_2,
$$
$$
N_{-c,b+2c+d} N_{a+b+c,-a-2b-2c-d} = |\text{Formulas (012) and (060)}| = \epsilon_{23}\delta_2 (-\epsilon_{28}\delta_{14}\gamma_2) = -\epsilon_{38}\delta_{124}\gamma_2,
$$
we have
$$
S = -\delta_5 N_{b+2c+d,a+b+c} + \frac{4\epsilon_{38}\delta_{124}\gamma_2}{2} - \epsilon_{38}\delta_{124}\gamma_2 = 0
$$
and therefore
$$
N_{b+2c+d,a+b+c} = \epsilon_{38}\delta_{1245}\gamma_2. \eqno{(087)}
$$

088. We have $(b+2c+d) + (a+b+c) + (-a-2b-3c-d) = 0,$ therefore  (see formula (087))
$$
N_{b+2c+d,a+b+c} = N_{a+b+c,-a-2b-3c-d} = N_{-a-2b-3c-d,b+2c+d} = \epsilon_{38}\delta_{1245}\gamma_2. \eqno{(088)}
$$

089. We have $(b+2c+d) + (a+b+c+d) + (-d) + (-a-2b-3c-d) = 0,$ means
$$
S = \frac{N_{b+2c+d,a+b+c+d} N_{-d,-a-2b-3c-d}}{1} + \frac{N_{a+b+c+d,-d} N_{b+2c+d,-a-2b-3c-d}}{1} +
$$
$$
+ \frac{N_{-d,b+2c+d} N_{a+b+c+d,-a-2b-3c-d}}{2} =  0.
$$
Since
$$
N_{b+2c+d,a+b+c+d} N_{-d,-a-2b-3c-d} = |\text{Table 2}| = -\epsilon_0 N_{b+2c+d,a+b+c+d},
$$
$$
N_{a+b+c+d,-d} N_{b+2c+d,-a-2b-3c-d} = |\text{Formulas (015) and (088)}| = -\epsilon_5 (-\epsilon_{38}\delta_{1245}\gamma_2) = 
$$
$$
= \epsilon_{358}\delta_{1245}\gamma_2,
$$
$$
N_{-d,b+2c+d} N_{a+b+c+d,-a-2b-3c-d} = |\text{Formulas (005) and (078)}| = 2\epsilon_3 (-2\epsilon_{58}\delta_{1245}\gamma_2) = 
$$
$$
= - 4\epsilon_{358}\delta_{1245}\gamma_2,
$$
we have
$$
S = -\epsilon_0 N_{b+2c+d,a+b+c+d} + \epsilon_{358}\delta_{1245}\gamma_2 - \frac{4\epsilon_{358}\delta_{1245}\gamma_2}{2}  = 0
$$
and therefore
$$
N_{b+2c+d,a+b+c+d} = -\epsilon_{3580}\delta_{1245}\gamma_2. \eqno{(089)}
$$

090. We have $(b+2c+d) + (a+b+c+d) + (-a-2b-3c-2d) = 0,$ therefore  (see formula (089))
$$
N_{b+2c+d,a+b+c+d} = N_{a+b+c+d,-a-2b-3c-2d} = N_{-a-2b-3c-2d,b+2c+d} = -\epsilon_{3580}\delta_{1245}\gamma_2. \eqno{(090)}
$$

091. We have $(b+2c+d) + (a+b+2c+d) + (-c) + (-a-2b-3c-2d) = 0,$ means
$$
S = \frac{N_{b+2c+d,a+b+2c+d} N_{-c,-a-2b-3c-2d}}{2} + \frac{N_{a+b+2c+d,-c} N_{b+2c+d,-a-2b-3c-2d}}{1} +
$$
$$
+ \frac{N_{-c,b+2c+d} N_{a+b+2c+d,-a-2b-3c-2d}}{1} =  0.
$$
Since
$$
N_{b+2c+d,a+b+2c+d} N_{-c,-a-2b-3c-2d} = |\text{Table 2}| = -2\delta_6 N_{b+2c+d,a+b+2c+d},
$$
$$
N_{a+b+2c+d,-c} N_{b+2c+d,-a-2b-3c-2d} = |\text{Formulas (034) and (090)}| = -\epsilon_{56}\delta_4 (\epsilon_{3580}\delta_{1245}\gamma_2) = 
$$
$$
= -\epsilon_{3680}\delta_{125}\gamma_2,
$$
$$
N_{-c,b+2c+d} N_{a+b+2c+d,-a-2b-3c-2d} = |\text{Formulas (012) and (070)}| = \epsilon_{23}\delta_2 (-\epsilon_{2680}\delta_{15}\gamma_2) = 
$$
$$
= -\epsilon_{3680}\delta_{125}\gamma_2,
$$
we have
$$
S = -\frac{2\delta_6 N_{b+2c+d,a+b+2c+d}}{2} - \epsilon_{3680}\delta_{125}\gamma_2 - \epsilon_{3680}\delta_{125}\gamma_2 = 0,
$$
and therefore
$$
N_{b+2c+d,a+b+2c+d} = -2\epsilon_{3680}\delta_{1256}\gamma_2. \eqno{(091)}
$$

092. We have $(b+2c+d) + (a+b+2c+d) + (-a-2b-4c-2d) = 0,$ therefore  (see formula (091))
$$
\frac{N_{b+2c+d,a+b+2c+d}}{2} = \frac{N_{a+b+2c+d,-a-2b-4c-2d}}{1} = \frac{N_{-a-2b-4c-2d,b+2c+d}}{1} = -\epsilon_{3680}\delta_{1256}\gamma_2. \eqno{(092)}
$$

093. We have $(b+2c+d) + (a+2b+2c+d) + (-b) + (-a-2b-4c-2d) = 0$ and $(-b) + (b+2c+d) \notin F_4,$ means
$$
S = \frac{N_{b+2c+d,a+2b+2c+d} N_{-b,-a-2b-4c-2d}}{2} + \frac{N_{a+2b+2c+d,-b} N_{b+2c+d,-a-2b-4c-2d}}{1} = 0.
$$
Since
$$
N_{b+2c+d,a+2b+2c+d} N_{-b,-a-2b-4c-2d} = |\text{Table 2}| = -\gamma_3 N_{b+2c+d,a+2b+2c+d},
$$
$$
N_{a+2b+2c+d,-b} N_{b+2c+d,-a-2b-4c-2d} = |\text{Formulas (046) and (092)}| = -\epsilon_{68}\gamma_2 (\epsilon_{3680}\delta_{1256}\gamma_2) = 
$$
$$
= -\epsilon_{30}\delta_{1256},
$$
we have
$$
S = -\frac{\gamma_3 N_{b+2c+d,a+2b+2c+d}}{2} - \epsilon_{30}\delta_{1256} = 0,
$$
and therefore
$$
N_{b+2c+d,a+2b+2c+d} = -2\epsilon_{30}\delta_{1256}\gamma_3. \eqno{(093)}
$$

094. We have $(b+2c+d) + (a+2b+2c+d) + (-a-3b-4c-2d) = 0,$ therefore  (see formula (093))
$$
\frac{N_{b+2c+d,a+2b+2c+d}}{2} = \frac{N_{a+2b+2c+d,-a-3b-4c-2d}}{1} = \frac{N_{-a-3b-4c-2d,b+2c+d}}{1} = -\epsilon_{30}\delta_{1256}\gamma_3. \eqno{(094)}
$$

095. We have $(b+2c+2d) + a + (-d) + (-a-b-2c-d) = 0$ and $a + (-d) \notin F_4,$ means
$$
S = \frac{N_{b+2c+2d,a} N_{-d,-a-b-2c-d}}{2} + \frac{N_{-d,b+2c+2d} N_{a,-a-b-2c-d}}{1} = 0.
$$
Since
$$
N_{b+2c+2d,a}N_{-d,-a-b-2c-d} = |\text{Table 2}| = -2\epsilon_7 N_{b+2c+2d,a},
$$
$$
N_{-d,b+2c+2d} N_{a,-a-b-2c-d} = |\text{Formulas (006) and (084)}| = \epsilon_4 (\epsilon_{36}\delta_{1234}\gamma_1) = \epsilon_{346}\delta_{1234}\gamma_1,
$$
we have
$$
S = -\frac{2\epsilon_7 N_{b+2c+2d,a}}{2} + \epsilon_{346}\delta_{1234}\gamma_1 = 0
$$
and therefore
$$
N_{b+2c+2d,a} = \epsilon_{3467}\delta_{1234}\gamma_1. \eqno{(095)}
$$

096. We have $(b+2c+2d) + a + (-a-b-2c-2d) = 0,$ therefore  (see formula (095))
$$
N_{b+2c+2d,a} = N_{a,-a-b-2c-2d} = N_{-a-b-2c-2d,b+2c+2d} = \epsilon_{3467}\delta_{1234}\gamma_1. \eqno{(096)}
$$

097. We have $(b+2c+2d) + (a+b) + (-d) + (-a-2b-2c-d) = 0$ and $(a+b)+(-d) \notin F_4,$ means
$$
S = \frac{N_{b+2c+2d,a+b} N_{-d,-a-2b-2c-d}}{2} + \frac{N_{-d,b+2c+2d} N_{a+b,-a-2b-2c-d}}{1} = 0.
$$
Since
$$
N_{b+2c+2d,a+b} N_{-d,-a-2b-2c-d} = |\text{Table 2}| = -2\epsilon_9 N_{b+2c+2d,a+b},
$$
$$
N_{-d,b+2c+2d} N_{a+b,-a-2b-2c-d} = |\text{Formulas (006) and (086)}| = \epsilon_4 (\epsilon_{38}\delta_{1234}\gamma_2) = \epsilon_{348}\delta_{1234}\gamma_2,
$$
we have
$$
S = -\frac{2\epsilon_9 N_{b+2c+2d,a+b}}{2} + \epsilon_{348}\delta_{1234}\gamma_2 = 0
$$
and therefore
$$
N_{b+2c+2d,a+b} = \epsilon_{3489}\delta_{1234}\gamma_2. \eqno{(097)}
$$

098. We have $(b+2c+2d) + (a+b) + (-a-2b-2c-2d) = 0,$ therefore  (see formula (097))
$$
N_{b+2c+2d,a+b} = N_{a+b,-a-2b-2c-2d} = N_{-a-2b-2c-2d,b+2c+2d} = \epsilon_{3489}\delta_{1234}\gamma_2. \eqno{(098)}
$$ 

099. We have $(b+2c+2d) + (a+b+c) + (-d) + (-a-2b-3c-d) = 0$ and $(a+b+c) + (-d)\notin F_4,$ means
$$
S = \frac{N_{b+2c+2d,a+b+c} N_{-d,-a-2b-3c-d}}{1} + \frac{N_{-d,b+2c+2d} N_{a+b+c,-a-2b-3c-d}}{1} =  0.
$$
Since
$$
N_{b+2c+2d,a+b+c} N_{-d,-a-2b-3c-d} = |\text{Table 2}| = -\epsilon_0 N_{b+2c+2d,a+b+c},
$$
$$
N_{-d,b+2c+2d} N_{a+b+c,-a-2b-3c-d} = |\text{Formulas (006) and (088)}| = \epsilon_4 (\epsilon_{38}\delta_{1245}\gamma_2) = \epsilon_{348}\delta_{1245}\gamma_2,
$$
we have
$$
-\epsilon_0 N_{b+2c+2d,a+b+c} + \epsilon_{348}\delta_{1245}\gamma_2 = 0
$$
and therefore
$$
N_{b+2c+2d,a+b+c} = \epsilon_{3480}\delta_{1245}\gamma_2. \eqno{(099)}
$$

100. We have $(b+2c+2d) + (a+b+c) + (-a-2b-3c-2d) = 0,$ therefore  (see formula (099))
$$
\frac{N_{b+2c+2d,a+b+c}}{1} = \frac{N_{a+b+c,-a-2b-3c-2d}}{2} = \frac{N_{-a-2b-3c-2d,b+2c+2d}}{1} =  \epsilon_{3480}\delta_{1245}\gamma_2. \eqno{(100)}
$$

101. We have $(b+2c+2d) + (a+b+2c) + (-c) + (-a-2b-3c-2d) = 0$ and $(-c) + (b+2c+2d) \notin F_4,$ means
$$
S = \frac{N_{b+2c+2d,a+b+2c} N_{-c,-a-2b-3c-2d}}{2} + \frac{N_{a+b+2c,-c} N_{b+2c+2d,-a-2b-3c-2d}}{1} = 0.
$$
Since
$$
N_{b+2c+2d,a+b+2c} N_{-c,-a-2b-3c-2d} = |\text{Table 2}| = -2\delta_6 N_{b+2c+2d,a+b+2c},
$$
$$
N_{a+b+2c,-c} N_{b+2c+2d,-a-2b-3c-2d} = |\text{Formulas (022) and (100)}| = -\delta_4 (-\epsilon_{3480}\delta_{1245}\gamma_2) = 
$$
$$
= \epsilon_{3480}\delta_{125}\gamma_2,
$$
we have
$$
S = -\frac{2\delta_6 N_{b+2c+2d,a+b+2c}}{2} + \epsilon_{3480}\delta_{125}\gamma_2 = 0
$$
and therefore
$$
N_{b+2c+2d,a+b+2c} = \epsilon_{3480}\delta_{1256}\gamma_2. \eqno{(101)}
$$

102. We have $(b+2c+2d) + (a+b+2c) + (-a-2b-4c-2d) = 0,$ therefore  (see formula (101))
$$
N_{b+2c+2d,a+b+2c} = N_{a+b+2c,-a-2b-4c-2d} = N_{-a-2b-4c-2d,b+2c+2d} = \epsilon_{3480}\delta_{1256}\gamma_2. \eqno{(102)}
$$

103. We have $(b+2c+2d)+(a+2b+2c)+(-b)+(-a-2b-4c-2d) = 0$ and $(b+2c+2d) + (-b) \notin F_4,$ means
$$
S = \frac{N_{b+2c+2d,a+2b+2c} N_{-b,-a-2b-4c-2d}}{2} + \frac{N_{a+2b+2c,-b} N_{b+2c+2d,-a-2b-4c-2d}}{2} = 0.
$$
Since
$$
N_{b+2c+2d,a+2b+2c} N_{-b,-a-2b-4c-2d} = |\text{Table 2}| = -\gamma_3 N_{b+2c+2d,a+2b+2c},
$$
$$
N_{a+2b+2c,-b} N_{b+2c+2d,-a-2b-4c-2d} = |\text{Formulas (026) and (102)}| = -\gamma_2 (-\epsilon_{3480}\delta_{1256}\gamma_2) = 
$$
$$
= \epsilon_{3480}\delta_{1256},
$$
we have
$$
S = -\frac{\gamma_3 N_{b+2c+2d,a+2b+2c}}{2} + \frac{\epsilon_{3480}\delta_{1256}}{2} = 0
$$
and therefore
$$
N_{b+2c+2d,a+2b+2c} = \epsilon_{3480}\delta_{1256}\gamma_3. \eqno{(103)}
$$

104. We have $(b+2c+2d) + (a+2b+2c) + (-a-3b-4c-2d) = 0,$ therefore  (see formula (103))
$$
N_{b+2c+2d,a+2b+2c} = N_{a+2b+2c,-a-3b-4c-2d} = N_{-a-3b-4c-2d,b+2c+2d} = \epsilon_{3480}\delta_{1256}\gamma_3. \eqno{(104)}
$$

105. We have $(a+b+c) + (a+2b+3c+2d) + (-a) + (-a-3b-4c-2d) = 0$ and $(a+2b+3c+2d) + (-a) \notin F_4,$ means
$$
S = \frac{N_{a+b+c,a+2b+3c+2d} N_{-a,-a-3b-4c-2d}}{2} + \frac{N_{-a,a+b+c} N_{a+2b+3c+2d,-a-3b-4c-2d}}{1} = 0.
$$
Since
$$
N_{a+b+c,a+2b+3c+2d} N_{-a,-a-3b-4c-2d} = |\text{Table 2}| = -\alpha_1 N_{a+b+c,a+2b+3c+2d},
$$
$$
N_{-a,a+b+c} N_{a+2b+3c+2d,-a-3b-4c-2d} = |\text{Formulas (042) and (048)}| = -\delta_{13}\gamma_1 (-\delta_{16}\gamma_3)  = 
$$
$$
= \delta_{36}\gamma_{13},
$$
we have
$$
S = -\frac{\alpha_1 N_{a+b+c,a+2b+3c+2d}}{2} + \delta_{36}\gamma_{13} = 0
$$
and therefore
$$
N_{a+b+c,a+2b+3c+2d} = 2\delta_{36}\gamma_{13}\alpha_1. \eqno{(105)}
$$

106. We have $(a+b+c) + (a+2b+3c+2d) + (-2a-3b-4c-2d) = 0,$ therefore  (see formula (105))
$$
\frac{N_{a+b+c,a+2b+3c+2d}}{2} = \frac{N_{a+2b+3c+2d,-2a-3b-4c-2d}}{1} = \frac{N_{-2a-3b-4c-2d,a+b+c}}{1} = \delta_{36}\gamma_{13}\alpha_1. \eqno{(106)}
$$

107. We have $(a+b+c+d) + (a+2b+3c+d) + (-a) + (-a-3b-4c-2d) = 0$ and $(a+2b+3c+d) + (-a) \notin F_4,$ means
$$
S = \frac{N_{a+b+c+d,a+2b+3c+d} N_{-a,-a-3b-4c-2d}}{2} + \frac{N_{-a,a+b+c+d} N_{a+2b+3c+d,-a-3b-4c-2d}}{1} = 0.
$$
Since
$$
N_{a+b+c+d,a+2b+3c+d} N_{-a,-a-3b-4c-2d} = |\text{Table 2}| = -\alpha_1 N_{a+b+c+d,a+2b+3c+d},
$$
$$
N_{-a,a+b+c+d} N_{a+2b+3c+d,-a-3b-4c-2d} = |\text{Formulas (052) and (072)}| = -\epsilon_{25}\delta_{13}\gamma_1 (\epsilon_{20}\delta_{16}\gamma_3) = 
$$
$$
 = -\epsilon_{50}\delta_{36}\gamma_{13},
$$
we have
$$
S = -\frac{\alpha_1 N_{a+b+c+d,a+2b+3c+d}}{2} - \epsilon_{50}\delta_{36}\gamma_{13} = 0
$$
and therefore
$$
N_{a+b+c+d,a+2b+3c+d} = -2\epsilon_{50}\delta_{36}\gamma_{13}\alpha_1. \eqno{(107)}
$$

108. We have $(a+b+c+d) + (a+2b+3c+d) + (-2a-3b-4c-2d) = 0,$ therefore  (see formula (107))
$$
\frac{N_{a+b+c+d,a+2b+3c+d}}{2} = \frac{N_{a+2b+3c+d,-2a-3b-4c-2d}}{1} = \frac{N_{-2a-3b-4c-2d,a+b+c+d}}{1} = -\epsilon_{50}\delta_{36}\gamma_{13}\alpha_1. \eqno{(108)}
$$

109. We have $(a+b+2c) + (a+2b+2c+2d) + (-a) + (-a-3b-4c-2d) = 0$ and $(a+2b+2c+2d) + (-a) \notin F_4,$ means
$$
S = \frac{N_{a+b+2c,a+2b+2c+2d} N_{-a,-a-3b-4c-2d}}{2} + \frac{N_{-a,a+b+2c} N_{a+2b+2c+2d,-a-3b-4c-2d}}{2} = 0.
$$
Since
$$
N_{a+b+2c,a+2b+2c+2d} N_{-a,-a-3b-4c-2d} = |\text{Table 2}| = -\alpha_1 N_{a+b+2c,a+2b+2c+2d},
$$
$$
N_{-a,a+b+2c} N_{a+2b+2c+2d,-a-3b-4c-2d} = |\text{Formulas (074) and (082)}| = -\delta_{1234}\gamma_1 (\epsilon_{90}\delta_{1256}\gamma_3) = 
$$
$$
 = -\epsilon_{90}\delta_{3456}\gamma_{13},
$$
we have
$$
S = -\frac{\alpha_1 N_{a+b+2c,a+2b+2c+2d}}{2} - \frac{\epsilon_{90}\delta_{3456}\gamma_{13}}{2} = 0
$$
and therefore
$$
N_{a+b+2c,a+2b+2c+2d} = -\epsilon_{90}\delta_{3456}\gamma_{13}\alpha_1. \eqno{(109)}
$$

110. We have $(a+b+2c) + (a+2b+2c+2d) +  (-2a-3b-4c-2d) = 0,$ therefore  (see formula (109))
$$
N_{a+b+2c,a+2b+2c+2d} = N_{a+2b+2c+2d,-2a-3b-4c-2d} = N_{-2a-3b-4c-2d,a+b+2c} = -\epsilon_{90}\delta_{3456}\gamma_{13}\alpha_1.\eqno{(110)}
$$

111. We have $(a+b+2c+d) + (a+2b+2c+d) + (-a) + (-a-3b-4c-2d) = 0$ and $(a+2b+2c+d) + (-a) \notin F_4,$ means
$$
S = \frac{N_{a+b+2c+d,a+2b+2c+d} N_{-a,-a-3b-4c-2d}}{2} + \frac{N_{-a,a+b+2c+d} N_{a+2b+2c+d,-a-3b-4c-2d}}{1} = 0.
$$
Since
$$
N_{a+b+2c+d,a+2b+2c+d} N_{-a,-a-3b-4c-2d} = |\text{Table 2}| = -\alpha_1 N_{a+b+2c+d,a+2b+2c+d},
$$
$$
N_{-a,a+b+2c+d} N_{a+2b+2c+d,-a-3b-4c-2d} = |\text{Formulas (084) and (094)}| = 
$$
$$
= -\epsilon_{36}\delta_{1234}\gamma_1 (-\epsilon_{30}\delta_{1256}\gamma_3)  = \epsilon_{60}\delta_{3456}\gamma_{13},
$$
we have
$$
S = -\frac{\alpha_1 N_{a+b+2c+d,a+2b+2c+d}}{2} + \epsilon_{60}\delta_{3456}\gamma_{13} = 0
$$
and therefore
$$
N_{a+b+2c+d,a+2b+2c+d} = 2\epsilon_{60}\delta_{3456}\gamma_{13}\alpha_1. \eqno{(111)}
$$

112. We have $(a+b+2c+d)+(a+2b+2c+d)+(-2a-3b-4c-2d) = 0,$ therefore  (see formula (111))
$$
\frac{N_{a+b+2c+d,a+2b+2c+d}}{2} = \frac{N_{a+2b+2c+d,-2a-3b-4c-2d}}{1} = \frac{N_{-2a-3b-4c-2d,a+b+2c+d}}{1} = \epsilon_{60}\delta_{3456}\gamma_{13}\alpha_1.\eqno{(112)}
$$

113. We have $(a+b+2c+2d) + (a+2b+2c) + (-a) + (-a-3b-4c-2d) = 0$ and $(a+2b+2c) + (-a) \notin F_4,$ means
$$
S = \frac{N_{a+b+2c+2d,a+2b+2c} N_{-a,-a-3b-4c-2d}}{2} + \frac{N_{-a,a+b+2c+2d} N_{a+2b+2c,-a-3b-4c-2d}}{2} = 0.
$$
Since
$$
N_{a+b+2c+2d,a+2b+2c} N_{-a,-a-3b-4c-2d} = |\text{Table 2}| = -\alpha_1 N_{a+b+2c+2d,a+2b+2c},
$$
$$
N_{-a,a+b+2c+2d} N_{a+2b+2c,-a-3b-4c-2d} = |\text{Formulas (096) and (104)}| = 
$$
$$
= (-\epsilon_{3467}\delta_{1234}\gamma_1)(\epsilon_{3480}\delta_{1256}\gamma_3) = -\epsilon_{6780}\delta_{3456}\gamma_{13},
$$
we have
$$
S = -\frac{\alpha_1 N_{a+b+2c+2d,a+2b+2c}}{2} - \frac{\epsilon_{6780}\delta_{3456}\gamma_{13}}{2} = 0
$$
and therefore
$$
N_{a+b+2c+2d,a+2b+2c}  = -\epsilon_{6780}\delta_{3456}\gamma_{13}\alpha_1. \eqno{(113)}
$$

114. We have $(a+b+2c+2d) + (a+2b+2c) + (-2a-3b-4c-2d) = 0,$ therefore  (see formula (113))
$$
N_{a+b+2c+2d,a+2b+2c} = N_{a+2b+2c,-2a-3b-4c-2d} = N_{-2a-3b-4c-2d,a+b+2c+2d} = -\epsilon_{6780}\delta_{3456}\gamma_{13}\alpha_1. \eqno{(114)}
$$

115. We have $c+(a+2b+2c+2d) + (-b-c-d) + (-a-b-2c-d) = 0$ and $(-b-c-d) + c \notin F_4,$ means
$$
S = \frac{N_{c,a+2b+2c+2d} N_{-b-c-d,-a-b-2c-d}}{1} + \frac{N_{a+2b+2c+2d,-b-c-d} N_{c,-a-b-2c-d}}{1} = 0.
$$
Since
$$
N_{c,a+2b+2c+2d} N_{-b-c-d,-a-b-2c-d} = |\text{Formula (069)}| = \epsilon_{2680}\delta_{15}\gamma_2 N_{c,a+2b+2c+2d},
$$
$$
N_{a+2b+2c+2d,-b-c-d} N_{c,-a-b-2c-d} = |\text{Formulas (062) and (034)}| = \epsilon_{2589}\delta_{14}\gamma_2 (-\epsilon_{56}\delta_4) = 
$$
$$
= -\epsilon_{2689}\delta_1\gamma_2,
$$
we have
$$
S = \epsilon_{2680}\delta_{15}\gamma_2 N_{c,a+2b+2c+2d} - \epsilon_{2689}\delta_1\gamma_2 = 0,
$$
and therefore
$$
N_{c,a+2b+2c+2d} = \epsilon_{90}\delta_5. \eqno{(115)}
$$

116. We have $c + (a+2b+2c+2d) + (-a-2b-3c-2d) = 0,$ therefore  (see formula (115))
$$
\frac{N_{c,a+2b+2c+2d}}{1} = \frac{N_{a+2b+2c+2d,-a-2b-3c-2d}}{1} = \frac{N_{-a-2b-3c-2d,c}}{2} = \epsilon_{90}\delta_5.\eqno{(116)}
$$

To conclude the proof, we present the diagrams that show the dependencies of the formulas.

\newpage


\thicklines
{

\end{center}

\centerline{Picture 6.}
\medskip

The theorem is completely proven. \hfill$\square$

\subsection{Chevalley's Commutator Formulas}

Let $\Phi$ be a reduced indecomposable root system, $K$ is a
field. According to \cite[Theorem 5.2.2]{Car72} the commutator
$$
[x_s(u),x_r(t)]=x_s(u)^{-1}x_r(t)^{-1}x_s(u)x_r(t),
$$
where $r,s\in\Phi$ and $u,t\in K,$ is equal to one if $r+s\notin\Phi$ and $r\ne-s,$ and expands
into the product of root elements according to the formula
$$
\left[x_s(u),x_r(t)\right]=
 \prod_{{\scriptsize\begin{array}{c} ir+js\in\Phi,\\ i,j>0\end{array}}}x_{ir+js}(C_{ij,rs}( -t)^iu^j),
$$
if $r+s\in\Phi.$ The product is taken over all pairs of positive integers $i,j$ for which $ir+js$ is a root, in order of increasing $i+j.$ The constant $C_{ij,rs}$
are integers and are determined by the formulas \cite[theorem 5.2.2]{Car72}:
\begin{align*}
    C_{i1,rs} &= M_{r,s,i}, \\[3mm]
    C_{1j,rs} &= (-1)^j M_{s,r,j}, \\[3mm]
    C_{32,rs} &= \frac{1}{3}M_{r+s,r,2}, \\[3mm]
    C_{23,rs} &= -\frac{2}{3}M_{s+r,s,2}.
\end{align*}
In turn, the numbers $M_{r,s,i}$ are expressed in terms of structure
constants $N_{r,s}$ of the corresponding Lie algebra according to the formula \cite[p. 61]{Car72}
$$
M_{r,s,i}=\frac{1}{i!}N_{r,s}N_{r,r+s}\ldots N_{r,(i-1)r+s},
$$

Let $r,s,r+s\in\Phi.$ We have 
$$
C_{11,rs} = N_{rs} = -N_{sr}
$$
and
$$
C_{11,-r,-s} = N_{-r,-s} = N_{sr}.
$$
Further, if $2r+s\in\Phi,$ then
$$
C_{21,rs}=\frac12 N_{rs}N_{r,r+s}=-\frac12 N_{sr}N_{r,r+s}
$$
and
$$
C_{21,-r,-s}=\frac12 N_{-r,-s}N_{-r,-r-s}=-\frac12 N_{sr}N_{r,r+s},
$$
and if $r+2s\in\Phi,$ then
$$
C_{12,rs}=\frac12 N_{sr}N_{s,r+s}
$$
and
$$
C_{12,-r,-s}=\frac12 N_{-s,-r}N_{-s,-r-s}=\frac12 N_{sr}N_{s,r+s}.
$$
It follows that:

1) if $2r+s,r+2s\notin\Phi,$ then
$$
[x_s(u),x_r(t)] = x_{r+s}(C_{11,rs}(-t)u) = x_{r+s}(N_{sr}tu)
$$ 
and
$$
[x_{-s}(u),x_{-r}(t)] = x_{-r-s}(C_{11,-r,-s}(-t)u) = x_{r+s}(-N_{sr}tu).
$$

2) if $2r+s\in\Phi,$ $r+2s\notin\Phi,$ then
$$
\begin{array}{rcl}
[x_s(u),x_r(t)] & = & x_{r+s}(C_{11,rs}(-t)u)\ x_{2r+s}(C_{21,rs}(-t)^2u) = \\[4mm]
                & = & x_{r+s}(N_{sr}tu)\ x_{2r+s}(-(1/2) N_{sr} N_{r,r+s} t^2u).\\
\end{array}
$$ 
and
$$
\begin{array}{rcl}
[x_{-s}(u),x_{-r}(t)] & = & x_{-r-s}(C_{11,-r,-s}(-t)u)\ x_{-2r-s}(C_{21,-r,-s}(-t)^2u) = \\[4mm]
                      & = & x_{-r-s}(-N_{sr}tu)\ x_{-2r-s}(-(1/2) N_{sr} N_{r,r+s} t^2u).\\
\end{array}
$$

3) if $2r+s\notin\Phi,$ $r+2s\in\Phi,$ then
$$
\begin{array}{rcl}
[x_s(u),x_r(t)] & = & x_{r+s}(C_{11,rs}(-t)u)\ x_{r+2s}(C_{12,rs}(-t)u^2) = \\[4mm]
                & = & x_{r+s}(N_{sr}tu)\ x_{r+2s}(-(1/2) N_{sr} N_{s,r+s} tu^2).\\
\end{array}
$$ 
and
$$
\begin{array}{rcl}
[x_{-s}(u),x_{-r}(t)] & = & x_{-r-s}(C_{11,-r,-s}(-t)u)\ x_{-r-2s}(C_{12,-r,-s}(-t)u^2) = \\[4mm]
                      & = & x_{-r-s}(-N_{sr}tu)\ x_{-r-2s}(-(1/2) N_{sr} N_{s,r+s} tu^2).\\
\end{array}
$$
\subsection{List of Formulas}

\vskip5mm
\begin{center}
\textit{Positive roots}
\end{center}

{\boldmath$d:$}

$$
[x_d(u),x_c(t)] = x_{c+d}(N_{d,c}tu) = x_{c+d}(\epsilon_1 tu).\eqno{(01)}
$$

$$
[x_d(u),x_{b+c}(t)] = x_{b+c+d}(N_{d,b+c}tu) = x_{b+c+d}(\epsilon_2 tu).\eqno{(02)}
$$

$$
\begin{array}{rcl}
[x_d(u),x_{b+2c}(t)] & = & x_{b+2c+d}(N_{d,b+2c}tu)\ x_{b+2c+2d}(-(1/2) N_{d,b+2c} N_{d,b+2c+d}tu^2)  \\[4mm]
                     & = &x_{b+2c+d}(\epsilon_3 tu)x_{b+2c+2d}(-\epsilon_{34} tu^2).
\end{array}\eqno{(03)}
$$
Explanation: 
$$
-(1/2) N_{d,b+2c} N_{d,b+2c+d} = -(1/2) (\epsilon_3) (2\epsilon_4) = -\epsilon_{34}.
$$

$$
[x_d(u),x_{b+2c+d}(t)] = x_{b+2c+2d}(N_{d,b+2c+d}tu) = x_{b+2c+2d}(2\epsilon_4 tu).\eqno{(04)}
$$

$$
[x_d(u),x_{a+b+c}(t)] = x_{a+b+c+d}(N_{d,a+b+c}tu) = x_{a+b+c+d}(\epsilon_5 tu).\eqno{(05)}
$$

$$
\begin{array}{rcl}
[x_d(u),x_{a+b+2c}(t)] & = & x_{a+b+2c+d}(N_{d,a+b+2c}tu)\times \\[4mm]
                       &   & x_{a+b+2c+2d}(-(1/2) N_{d,a+b+2c} N_{d,a+b+2c+d}tu^2)  \\[4mm]
                       & = & x_{a+b+2c+d}(\epsilon_6 tu)\ x_{a+b+2c+2d}(-\epsilon_{67} tu^2).
\end{array} \eqno{(06)}
$$
Explanation: 
$$
-(1/2) N_{d,a+b+2c} N_{d,a+b+2c+d} = -(1/2) (\epsilon_6) (2\epsilon_7) = -\epsilon_{67}.
$$

$$
[x_d(u),x_{a+b+2c+d}(t)] = x_{a+b+2c+2d}(N_{d,a+b+2c+d}tu) = x_{a+b+2c+2d}(2\epsilon_7 tu).\eqno{(07)}
$$

$$
\begin{array}{rcl}
[x_d(u),x_{a+2b+2c}(t)] & = & x_{a+2b+2c+d}(N_{d,a+2b+2c}tu)\times\\[4mm]
                        &   & x_{a+2b+2c+2d}(-(1/2) N_{d,a+2b+2c} N_{d,a+2b+2c+d}tu^2)  \\[4mm]
                        & = & x_{a+2b+2c+d}(\epsilon_8 tu)\ x_{a+2b+2c+2d}(-\epsilon_{89} tu^2).
\end{array} \eqno{(08)}
$$
Explanation: 
$$
-(1/2) N_{d,a+2b+2c} N_{d,a+2b+2c+d} = -(1/2) (\epsilon_8) (2\epsilon_9) = -\epsilon_{89}.
$$

$$
[x_d(u),x_{a+2b+2c+d}(t)] = x_{a+2b+2c+2d}(N_{d,a+2b+2c+d}tu) = x_{a+2b+2c+2d}(2\epsilon_9 tu).\eqno{(09)}
$$

$$
[x_d(u),x_{a+2b+3c+d}(t)] = x_{a+2b+3c+2d}(N_{d,a+2b+3c+d}tu) = x_{a+2b+3c+2d}(\epsilon_0 tu).\eqno{(10)}
$$

{\boldmath$c:$}

$$
\begin{array}{rcl}
[x_c(u),x_{b}(t)] & = & x_{b+c}(N_{c,b}tu)\ x_{b+2c}(-(1/2) N_{c,b} N_{c,b+c}tu^2) \\[4mm]
                  & = & x_{b+c}(\delta_1 tu)\ x_{b+2c}(-\delta_{12} tu^2).
\end{array} \eqno{(11)}
$$
Explanation: 
$$
-(1/2) N_{c,b} N_{c,b+c} = -(1/2) (\delta_1) (2\delta_2) = -\delta_{12}.
$$

$$
[x_c(u),x_{b+c}(t)] = x_{b+2c}(N_{c,b+c}tu) = x_{b+2c}(2\delta_2 tu).\eqno{(12)}
$$

$$
[x_c(u),x_{b+c+d}(t)] = x_{b+2c+d}(N_{c,b+c+d}tu) = x_{b+2c+d}(\epsilon_{23}\delta_2 tu).\eqno{(13)}
$$

$$
\begin{array}{rcl}
[x_c(u),x_{a+b}(t)] & = & x_{a+b+c}(N_{c,a+b}tu)\ x_{a+b+2c}(-(1/2) N_{c,a+b} N_{c,a+b+c}tu^2)  \\[4mm]
                    & = & x_{a+b+c}(\delta_3 tu)\ x_{a+b+2c}(-\delta_{34} tu^2).
\end{array} \eqno{(14)}
$$
Explanation: 
$$
-(1/2) N_{c,a+b} N_{c,a+b+c} = -(1/2) (\delta_3) (2\delta_4) = -\delta_{34}.
$$

$$
[x_c(u),x_{a+b+c}(t)] = x_{a+b+2c}(N_{c,a+b+c}tu) = x_{a+b+2c}(2\delta_4 tu).\eqno{(15)}
$$

$$
[x_c(u),x_{a+b+c+d}(t)] = x_{a+b+2c+d}(N_{c,a+b+c+d}tu) = x_{a+b+2c+d}(\epsilon_{56}\delta_4 tu).\eqno{(16)}
$$

$$
[x_c(u),x_{a+2b+2c+d}(t)] = x_{a+2b+3c+d}(N_{c,a+2b+2c+d}tu) = x_{a+2b+3c+d}(\delta_5 tu).\eqno{(17)}
$$

$$
\begin{array}{rcl}
[x_c(u),x_{a+2b+2c+2d}(t)] & = & x_{a+2b+3c+2d}(N_{c,a+2b+2c+2d}tu)\times\\[4mm]
                           & & x_{a+2b+4c+2d}(-(1/2) N_{c,a+2b+2c+2d} N_{c,a+2b+3c+2d} tu^2)  \\[4mm]
                           & = & x_{a+2b+3c+2d}(\epsilon_{90}\delta_5 tu)\ x_{a+2b+4c+2d}(-\epsilon_{90}\delta_{56} tu^2).
\end{array} \eqno{(18)}
$$
Explanation: 
$$
-(1/2) N_{c,a+2b+2c+2d} N_{c,a+2b+3c+2d} = -(1/2) (\epsilon_{90}\delta_5) (2\delta_6) = -\epsilon_{90}\delta_{56}.
$$

$$
[x_c(u),x_{a+2b+3c+2d}(t)] = x_{a+2b+4c+2d}(N_{c,a+2b+3c+2d}tu) = x_{a+2b+4c+2d}(2\delta_6 tu).\eqno{(19)}
$$

{\boldmath$c+d:$}

$$
\begin{array}{rcl}
[x_{c+d}(u),x_{b}(t)] & = & x_{b+c+d}(N_{c+d,b} tu)\ x_{b+2c+2d}(-(1/2) N_{c+d,b} N_{c+d,b+c+d} tu^2)  \\[4mm]
                      & = & x_{b+c+d}(\epsilon_{12}\delta_1 tu)\ x_{b+2c+2d}(-\epsilon_{34}\delta_{12} tu^2).\\
\end{array} \eqno{(20)}
$$
Explanation: 
$$
-(1/2) N_{c+d,b} N_{c+d,b+c+d} = -(1/2) (\epsilon_{12}\delta_1) (2\epsilon_{1234}\delta_2) = -\epsilon_{34}\delta_{12}.
$$

$$
[x_{c+d}(u),x_{b+c+d}(t)] = x_{b+2c+2d}(N_{c+d,b+c+d}tu) = x_{b+2c+2d}(2\epsilon_{1234}\delta_2 tu).\eqno{(21)}
$$

$$
[x_{c+d}(u),x_{b+c}(t)] = x_{b+2c+d}(N_{c+d,b+c}tu) = x_{b+2c+d}(\epsilon_{13}\delta_2 tu).\eqno{(22)}
$$

$$
\begin{array}{rcl}
[x_{c+d}(u),x_{a+b}(t)] & = & x_{a+b+c+d}(N_{c+d,a+b} tu)\ x_{a+b+2c+2d}(-(1/2) N_{c+d,a+b} N_{c+d,a+b+c+d} tu^2)  \\[4mm]
                        & = & x_{a+b+c+d}(\epsilon_{15}\delta_3 tu)\ x_{a+b+2c+2d}(-\epsilon_{67}\delta_{34} tu^2).\\
\end{array} \eqno{(23)}
$$
Explanation: 
$$
-(1/2) N_{c+d,a+b} N_{c+d,a+b+c+d} = -(1/2) (\epsilon_{15}\delta_3) (2\epsilon_{1567}\delta_4) = -\epsilon_{67}\delta_{34}.
$$

$$
[x_{c+d}(u),x_{a+b+c+d}(t)] = x_{a+b+2c+2d}(N_{c+d,a+b+c+d} tu) = x_{a+b+2c+2d}(2\epsilon_{1567}\delta_4 tu).\eqno{(24)}
$$

$$
[x_{c+d}(u),x_{a+b+c}(t)] = x_{a+b+2c+d}(N_{c+d,a+b+c} tu) = x_{a+b+2c+d}(\epsilon_{16}\delta_4 tu).\eqno{(25)}
$$

$$
\begin{array}{rcl}
[x_{c+d}(u),x_{a+2b+2c}(t)] & = & x_{a+2b+3c+d}(N_{c+d,a+2b+2c} tu)\times\\[4mm]
                            & & x_{a+2b+4c+2d}(-(1/2) N_{c+d,a+2b+2c} N_{c+d,a+2b+3c+d} tu^2) \\[4mm]
                            & = & x_{a+2b+3c+d}(-\epsilon_{18}\delta_5 tu)\ x_{a+2b+4c+2d}(-\epsilon_{80}\delta_{56} tu^2).\\
\end{array} \eqno{(26)}
$$
Explanation: 
$$
-(1/2) N_{c+d,a+2b+2c} N_{c+d,a+2b+3c+d}  = -(1/2) (-\epsilon_{18}\delta_5) (-2\epsilon_{10}\delta_6) = -\epsilon_{80}\delta_{56}.
$$

$$
[x_{c+d}(u),x_{a+2b+3c+d}(t)] = x_{a+2b+4c+2d}(N_{c+d,a+2b+3c+d} tu) = x_{a+2b+4c+2d}(-2\epsilon_{10}\delta_6 tu).\eqno{(27)}
$$

$$
[x_{c+d}(u),x_{a+2b+2c+d}(t)] = x_{a+2b+3c+2d}(N_{c+d,a+2b+2c+d} tu) = x_{a+2b+3c+2d}(-\epsilon_{10}\delta_5 tu).\eqno{(28)}
$$

{\boldmath$b:$}

$$
[x_b(u),x_a(t)] = x_{a+b}(N_{b,a} tu) = x_{a+b}(\gamma_1 tu).\eqno{(29)}
$$

$$
[x_b(u),x_{a+b+2c}(t)] = x_{a+2b+2c}(N_{b,a+b+2c} tu) = x_{a+2b+2c}(\gamma_2 tu).\eqno{(30)}
$$

$$
[x_b(u),x_{a+b+2c+d}(t)] = x_{a+2b+2c+d}(N_{b,a+b+2c+d} tu) = x_{a+2b+2c+d}(\epsilon_{68}\gamma_2 tu).\eqno{(31)}
$$

$$
[x_b(u),x_{a+b+2c+2d}(t)] = x_{a+2b+2c+2d}(N_{b,a+b+2c+2d} tu) = x_{a+2b+2c+2d}(\epsilon_{6789}\gamma_2 tu).\eqno{(32)}
$$

$$
[x_b(u),x_{a+2b+4c+2d}(t)] = x_{a+3b+4c+2d}(N_{b,a+2b+4c+2d} tu) = x_{a+3b+4c+2d}(\gamma_3 tu).\eqno{(33)}
$$

{\boldmath$b+c:$}

$$
\begin{array}{rcl}
[x_{b+c}(u),x_a(t)] & = & x_{a+b+c}(N_{b+c,a} tu)\ x_{a+2b+2c}(-(1/2) N_{b+c,a} N_{b+c,a+b+c} tu^2)  \\[4mm]
                    & = & x_{a+b+c}(\delta_{13}\gamma_1 tu)\ x_{a+2b+2c}(\delta_{34}\gamma_{12} tu^2).\\
\end{array} \eqno{(34)}
$$
Explanation: 
$$
-(1/2) N_{b+c,a} N_{b+c,a+b+c}  = -(1/2) (\delta_{13}\gamma_1) (-2\delta_{14}\gamma_2) = \delta_{34}\gamma_{12}.
$$

$$
[x_{b+c}(u),x_{a+b+c}(t)] = x_{a+2b+2c}(N_{b+c,a+b+c} tu) = x_{a+2b+2c}(-2\delta_{14}\gamma_2 tu).\eqno{(35)}
$$

$$
[x_{b+c}(u),x_{a+b+c+d}(t)] = x_{a+2b+2c+d}(N_{b+c,a+b+c+d} tu) = x_{a+2b+2c+d}(-\epsilon_{58}\delta_{14}\gamma_2 tu).\eqno{(36)}
$$

$$
[x_{b+c}(u),x_{a+b+2c+d}(t)] = x_{a+2b+3c+d}(N_{b+c,a+b+2c+d} tu) = x_{a+2b+3c+d}(\epsilon_{68}\delta_{15}\gamma_2 tu).\eqno{(37)}
$$

$$
\begin{array}{rcl}
[x_{b+c}(u),x_{a+b+2c+2d}(t)] & = & x_{a+2b+3c+2d}(N_{b+c,a+b+2c+2d} tu)\times\\[4mm]
                              & & x_{a+3b+4c+2d}(-(1/2) N_{b+c,a+b+2c+2d} N_{b+c,a+2b+3c+2d} tu^2)  \\[4mm]
                              & = & x_{a+2b+3c+2d}(\epsilon_{6780}\delta_{15}\gamma_2 tu)\ x_{a+3b+4c+2d}(\epsilon_{6780}\delta_{56}\gamma_{23} tu^2).\\
\end{array} \eqno{(38)}
$$
Explanation: 
$$
-(1/2) N_{b+c,a+b+2c+2d} N_{b+c,a+2b+3c+2d}  = -(1/2) (\epsilon_{6780}\delta_{15}\gamma_2) (-2\delta_{16}\gamma_3) = \epsilon_{6780}\delta_{56}\gamma_{23}.
$$

$$
[x_{b+c}(u),x_{a+2b+3c+2d}(t)] = x_{a+3b+4c+2d}(N_{b+c,a+2b+3c+2d} tu) = x_{a+b}(-2\delta_{16}\gamma_3 tu).\eqno{(39)}
$$

{\boldmath$b+c+d:$}

$$
\begin{array}{rcl}
[x_{b+c+d}(u),x_{a}(t)] & = & x_{a+b+c+d}(N_{b+c+d,a} tu)\times \\[4mm]
                        &   & x_{a+2b+2c+2d}(-(1/2) N_{b+c+d,a} N_{b+c+d,a+b+c+d} tu^2)  \\[4mm]
                        & = & x_{a+b+c+d}(\epsilon_{25}\delta_{13}\gamma_1 tu)\ x_{a+2b+2c+2d}(\epsilon_{89}\delta_{34}\gamma_{12} tu^2).\\
\end{array} \eqno{(40)}
$$
Explanation: 
$$
-(1/2) N_{b+c+d,a} N_{b+c+d,a+b+c+d}  = -(1/2) (\epsilon_{25}\delta_{13}\gamma_1) (-2\epsilon_{2589}\delta_{14}\gamma_2) = \epsilon_{89}\delta_{34}\gamma_{12}.
$$

$$
[x_{b+c+d}(u),x_{a+b+c+d}(t)] = x_{a+2b+2c+2d}(N_{b+c+d,a+b+c+d} tu) = x_{a+2b+2c+2d}(-2\epsilon_{2589}\delta_{14}\gamma_2 tu).\eqno{(41)}
$$

$$
[x_{b+c+d}(u),x_{a+b+c}(t)] = x_{a+2b+2c+d}(N_{b+c+d,a+b+c} tu) = x_{a+2b+2c+d}(-\epsilon_{28}\delta_{14}\gamma_2 tu).\eqno{(42)}
$$

$$
\begin{array}{rcl}
[x_{b+c+d}(u),x_{a+b+2c}(t)] & = & x_{a+2b+3c+d}(N_{b+c+d,a+b+2c} tu)\times\\[4mm] 
                             & & x_{a+3b+4c+2d}(-(1/2) N_{b+c+d,a+b+2c} N_{b+c+d,a+2b+3c+d} tu^2)  \\[4mm]
                             & = & x_{a+2b+3c+d}(-\epsilon_{28}\delta_{15}\gamma_2 tu)\ x_{a+3b+4c+2d}(\epsilon_{80}\delta_{56}\gamma_{23} tu^2).\\
\end{array} \eqno{(43)}
$$
Explanation: 
$$
-(1/2) N_{b+c+d,a+b+2c} N_{b+c+d,a+2b+3c+d}  = -(1/2) (-\epsilon_{28}\delta_{15}\gamma_2) (2\epsilon_{20}\delta_{16}\gamma_3) = \epsilon_{80}\delta_{56}\gamma_{23}.
$$

$$
\begin{array}{rcl}
[x_{b+c+d}(u),x_{a+2b+3c+d}(t)] & = & x_{a+3b+4c+2d}(N_{b+c+d,a+2b+3c+d} tu) \\[4mm]
                                & = & x_{a+3b+4c+2d}(2\epsilon_{20}\delta_{16}\gamma_3 tu).\\
\end{array}\eqno{(44)}
$$

$$
\begin{array}{rcl}
[x_{b+c+d}(u),x_{a+b+2c+d}(t)] & = & x_{a+2b+3c+2d}(N_{b+c+d,a+b+2c+d} tu) \\[4mm]
                               & = & x_{a+2b+3c+2d}(-\epsilon_{2680}\delta_{15}\gamma_2 tu).
\end{array} \eqno{(45)}
$$

{\boldmath$b+2c:$}

$$
[x_{b+2c}(u),x_{a}(t)] = x_{a+b+2c}(N_{b+2c,a} tu) = x_{a+b+2c}(\delta_{1234}\gamma_1 tu).\eqno{(46)}
$$

$$
[x_{b+2c}(u),x_{a+b}(t)] = x_{a+2b+2c}(N_{b+2c,a+b} tu) = x_{a+2b+2c}(\delta_{1234}\gamma_2 tu).\eqno{(47)}
$$

$$
[x_{b+2c}(u),x_{a+b+c+d}(t)] = x_{a+2b+3c+d}(N_{b+2c,a+b+c+d} tu) = x_{a+2b+3c+d}(-\epsilon_{58}\delta_{1245}\gamma_2 tu).\eqno{(48)}
$$

$$
\begin{array}{rcl}
[x_{b+2c}(u),x_{a+b+2c+2d}(t)] & = & x_{a+2b+4c+2d}(N_{b+2c,a+b+2c+2d} tu) \\[4mm]
                               & = & x_{a+2b+4c+2d}(\epsilon_{6780}\delta_{1256}\gamma_2 tu).
\end{array} \eqno{(49)}
$$

$$
\begin{array}{rcl}
[x_{b+2c}(u),x_{a+2b+2c+2d}(t)] & = & x_{a+3b+4c+2d}(N_{b+2c,a+2b+2c+2d} tu) \\[4mm]
                                & = & x_{a+3b+4c+2d}(\epsilon_{90}\delta_{1256}\gamma_3 tu).
\end{array} \eqno{(50)}
$$

{\boldmath$b+2c+d:$}

$$
\begin{array}{rcl}
[x_{b+2c+d}(u),x_{a}(t)] & = & x_{a+b+2c+d}(N_{b+2c+d,a} tu)\times\\[4mm]
                         &   & x_{a+2b+4c+2d}(-(1/2) N_{b+2c+d,a} N_{b+2c+d,a+b+2c+d} tu^2) \\[4mm]
												 & = & x_{a+b+2c+d}(\epsilon_{36}\delta_{1234}\gamma_1 tu)\ x_{a+2b+4c+2d}(\epsilon_{80}\delta_{3456}\gamma_{12} tu^2).
\end{array} \eqno{(43)}
$$
Explanation: 
$$
-(1/2)  N_{b+2c+d,a} N_{b+2c+d,a+b+2c+d}  = -(1/2) (\epsilon_{36}\delta_{1234}\gamma_1) (-2\epsilon_{3680}\delta_{1256}\gamma_2) = \epsilon_{80}\delta_{3456}\gamma_{12}.
$$

$$
\begin{array}{rcl}
[x_{b+2c+d}(u),x_{a+b+2c+d}(t)] & = & x_{a+2b+4c+2d}(N_{b+2c+d,a+b+2c+d} tu) \\[4mm]
                                & = & x_{a+2b+4c+2d}(-2\epsilon_{3680}\delta_{1256}\gamma_2 tu).
\end{array} \eqno{(52)}
$$

$$
\begin{array}{rcl}
[x_{b+2c+d}(u),x_{a+b}(t)] & = & x_{a+2b+2c+d}(N_{b+2c+d,a+b} tu)\times\\[4mm] 
                             & & x_{a+3b+4c+2d}(-(1/2) N_{b+2c+d,a+b} N_{b+2c+d,a+2b+2c+d} tu^2)  \\[4mm]
                             & = & x_{a+2b+2c+d}(\epsilon_{38}\delta_{1234}\gamma_2 tu)\ x_{a+3b+4c+2d}(\epsilon_{80}\delta_{3456}\gamma_{23} tu^2).\\
\end{array} \eqno{(53)}
$$
Explanation: 
$$
-(1/2) N_{b+2c+d,a+b} N_{b+2c+d,a+2b+2c+d}  = -(1/2) (\epsilon_{38}\delta_{1234}\gamma_2) (-2\epsilon_{30}\delta_{1256}\gamma_3) = \epsilon_{80}\delta_{3456}\gamma_{23}.
$$

$$
\begin{array}{rcl}
[x_{b+2c+d}(u),x_{a+2b+2c+d}(t)] & = & x_{a+3b+4c+2d}(N_{b+2c+d,a+2b+2c+d} tu) \\[4mm]
                                 & = & x_{a+3b+4c+2d}(-2\epsilon_{30}\delta_{1256}\gamma_3 tu).
\end{array} \eqno{(54)}
$$

$$
[x_{b+2c+d}(u),x_{a+b+c}(t)] = x_{a+2b+3c+d}(N_{b+2c+d,a+b+c} tu) = x_{a+2b+3c+d}(\epsilon_{38}\delta_{1245}\gamma_2 tu).\eqno{(55)}
$$

$$
\begin{array}{rcl}
[x_{b+2c+d}(u),x_{a+b+c+d}(t)] & = & x_{a+2b+3c+2d}(N_{b+2c+d,a+b+c+d} tu) \\[4mm]
                               & = & x_{a+2b+3c+2d}(-\epsilon_{3580}\delta_{1245}\gamma_2 tu).
\end{array} \eqno{(56)}
$$

{\boldmath$b+2c+2d:$}

$$
[x_{b+2c+2d}(u),x_a(t)] = x_{a+b+2c+2d}(N_{b+2c+2d,a} tu) = x_{a+b+2c+2d}(\epsilon_{3467}\delta_{1234}\gamma_1 tu).\eqno{(57)}
$$

$$
[x_{b+2c+2d}(u),x_{a+b}(t)] = x_{a+2b+2c+2d}(N_{b+2c+2d,a+b} tu) = x_{a+2b+2c+2d}(\epsilon_{3489}\delta_{1234}\gamma_2 tu).\eqno{(58)}
$$

$$
\begin{array}{rcl}
[x_{b+2c+2d}(u),x_{a+b+c}(t)] & = & x_{a+2b+3c+2d}(N_{b+2c+2d,a+b+c} tu) \\[4mm]
                              & = & x_{a+2b+3c+2d}(\epsilon_{3480}\delta_{1245}\gamma_2 tu).
\end{array}\eqno{(59)}
$$

$$
\begin{array}{rcl}
[x_{b+2c+2d}(u),x_{a+b+2c}(t)] & = & x_{a+2b+4c+2d}(N_{b+2c+2d,a+b+2c} tu) \\[4mm]
                               & = & x_{a+2b+4c+2d}(\epsilon_{3480}\delta_{1256}\gamma_2 tu).
\end{array} \eqno{(60)}
$$

$$
\begin{array}{rcl}
[x_{b+2c+2d}(u),x_{a+2b+2c}(t)] & = & x_{a+3b+4c+2d}(N_{b+2c+2d,a+2b+2c} tu) \\[4mm]
                                & = & x_{a+3b+4c+2d}(\epsilon_{3480}\delta_{1256}\gamma_3 tu).
\end{array} \eqno{(61)}
$$

{\boldmath$a:$}

$$
\begin{array}{rcl}
[x_a(u),x_{a+3b+4c+2d}(t)] & = & x_{2a+3b+4c+2d}(N_{a,a+3b+4c+2d} tu) \\[4mm]
                           & = & x_{2a+3b+4c+2d}(\alpha_1 tu).
\end{array} \eqno{(62)}
$$

{\boldmath$a+b:$}

$$
\begin{array}{rcl}
[x_{a+b}(u),x_{a+2b+4c+2d}(t)] & = & x_{2a+3b+4c+2d}(N_{a+b,a+2b+4c+2d} tu)\\[4mm]
                               & = & x_{2a+3b+4c+2d}(-\gamma_{13}\alpha_1 tu).
\end{array}\eqno{(63)}
$$

{\boldmath$a+b+c:$}

$$
\begin{array}{rcl}
[x_{a+b+c}(u),x_{a+2b+3c+2d}(t)] & = & x_{2a+3b+4c+2d}(N_{a+b+c,a+2b+3c+2d} tu) \\[4mm]
                                 & = & x_{2a+3b+4c+2d}(2\delta_{36}\gamma_{13}\alpha_1 tu). 
\end{array}\eqno{(64)}
$$

{\boldmath$a+b+c+d:$}

$$
\begin{array}{rcl}
[x_{a+b+c+d}(u),x_{a+2b+3c+d}(t)] & = & x_{2a+3b+4c+2d}(N_{a+b+c+d,a+2b+3c+d} tu)\\[4mm] 
                                  & = & x_{2a+3b+4c+2d}(-2\epsilon_{50}\delta_{36}\gamma_{13}\alpha_1 tu).
\end{array} \eqno{(65)}
$$

{\boldmath$a+b+2c:$}

$$
\begin{array}{rcl}
[x_{a+b+2c}(u),x_{a+2b+2c+2d}(t)] & = & x_{2a+3b+4c+2d}(N_{a+b+2c,a+2b+2c+2d} tu)\\[4mm]
                                  & = & x_{2a+3b+4c+2d}(-\epsilon_{90}\delta_{3456}\gamma_{13}\alpha_1 tu).
\end{array} \eqno{(66)}
$$

{\boldmath$a+b+2c+d:$}

$$
\begin{array}{rcl}
[x_{a+b+2c+d}(u),x_{a+2b+2c+d}(t)] & = & x_{2a+3b+4c+2d}(N_{a+b+2c+d,a+2b+2c+d} tu)\\[4mm]
                                   & = & x_{2a+3b+4c+2d}(2\epsilon_{60}\delta_{3456}\gamma_{13}\alpha_1 tu).
\end{array} \eqno{(67)}
$$

{\boldmath$a+b+2c+2d:$}

$$
\begin{array}{rcl}
[x_{a+b+2c+2d}(u),x_{a+2b+2c}(t)] & = & x_{2a+3b+4c+2d}(N_{a+b+2c+2d,a+2b+2c} tu) \\[4mm]
                                  & = & x_{2a+3b+4c+2d}(-\epsilon_{6780}\delta_{3456}\gamma_{13}\alpha_1 tu).
\end{array} \eqno{(68)}
$$

\vskip5mm
\begin{center}
\textit{Negative roots}
\end{center}

{\boldmath$-d:$}

$$
[x_{-d}(u),x_{-c}(t)] = x_{-c-d}(-N_{d,c}tu) = x_{-c-d}(-\epsilon_1 tu).\eqno{(01)}
$$

$$
[x_{-d}(u),x_{-b-c}(t)] = x_{-b-c-d}(-N_{d,b+c}tu) = x_{-b-c-d}(-\epsilon_2 tu).\eqno{(02)}
$$

$$
\begin{array}{rcl}
[x_{-d}(u),x_{-b-2c}(t)] & = & x_{-b-2c-d}(-N_{d,b+2c}tu)\ x_{-b-2c-2d}(-(1/2) N_{d,b+2c} N_{d,b+2c+d}tu^2) \\[4mm]
                     & = &x_{-b-2c-d}(-\epsilon_3 tu)x_{-b-2c-2d}(-\epsilon_{34} tu).
\end{array}\eqno{(03)}
$$
Explanation: 
$$
-(1/2) N_{d,b+2c} N_{d,b+2c+d} = -(1/2) (\epsilon_3) (2\epsilon_4) = -\epsilon_{34}.
$$

$$
[x_{-d}(u),x_{-b-2c-d}(t)] = x_{-b-2c-2d}(-N_{d,b+2c+d}tu) = x_{-b-2c-2d}(-2\epsilon_4 tu).\eqno{(04)}
$$

$$
[x_{-d}(u),x_{-a-b-c}(t)] = x_{-a-b-c-d}(-N_{d,a+b+c}tu) = x_{-a-b-c-d}(-\epsilon_5 tu).\eqno{(05)}
$$

$$
\begin{array}{rcl}
[x_{-d}(u),x_{-a-b-2c}(t)] & = & x_{-a-b-2c-d}(-N_{d,a+b+2c}tu)\times \\[4mm]
                           &   & x_{-a-b-2c-2d}(-(1/2) N_{d,a+b+2c} N_{d,a+b+2c+d}tu^2)  \\[4mm]
                           & = & x_{-a-b-2c-d}(-\epsilon_6 tu)\ x_{-a-b-2c-2d}(-\epsilon_{67} tu^2).
\end{array} \eqno{(06)}
$$
Explanation: 
$$
-(1/2) N_{d,a+b+2c} N_{d,a+b+2c+d} = -(1/2) (\epsilon_6) (2\epsilon_7) = -\epsilon_{67}.
$$

$$
[x_{-d}(u),x_{-a-b-2c-d}(t)] = x_{-a-b-2c-2d}(-N_{d,a+b+2c+d}tu) = x_{-a-b-2c-2d}(-2\epsilon_7 tu).\eqno{(07)}
$$

$$
\begin{array}{rcl}
[x_{-d}(u),x_{-a-2b-2c}(t)] & = & x_{-a-2b-2c-d}(-N_{d,a+2b+2c}tu)\times\\[4mm]
                            &   & x_{-a-2b-2c-2d}(-(1/2) N_{d,a+2b+2c} N_{d,a+2b+2c+d}tu^2) \\[4mm]
                            & = & x_{-a-2b-2c-d}(-\epsilon_8 tu)\ x_{-a-2b-2c-2d}(-\epsilon_{89} tu^2).
\end{array} \eqno{(08)}
$$
Explanation: 
$$
-(1/2) N_{d,a+2b+2c} N_{d,a+2b+2c+d} = -(1/2) (\epsilon_8) (2\epsilon_9) = -\epsilon_{89}.
$$

$$
[x_{-d}(u),x_{-a-2b-2c-d}(t)] = x_{-a-2b-2c-2d}(-N_{d,a+2b+2c+d}tu) = x_{-a-2b-2c-2d}(-2\epsilon_9 tu).\eqno{(09)}
$$

$$
[x_{-d}(u),x_{-a-2b-3c-d}(t)] = x_{-a-2b-3c-2d}(-N_{d,a+2b+3c+d}tu) = x_{-a-2b-3c-2d}(-\epsilon_0 tu).\eqno{(10)}
$$

{\boldmath$-c:$}

$$
\begin{array}{rcl}
[x_{-c}(u),x_{-b}(t)] & = & x_{-b-c}(-N_{c,b}tu)\ x_{-b-2c}(-(1/2) N_{c,b} N_{c,b+c}tu^2) \\[4mm]
                      & = & x_{-b-c}(-\delta_1 tu)\ x_{-b-2c}(-\delta_{12} tu^2).
\end{array} \eqno{(11)}
$$
Explanation: 
$$
-(1/2) N_{c,b} N_{c,b+c} = -(1/2) (\delta_1) (2\delta_2) = -\delta_{12}.
$$

$$
[x_{-c}(u),x_{-b-c}(t)] = x_{-b-2c}(-N_{c,b+c}tu) = x_{-b-2c}(-2\delta_2 tu).\eqno{(12)}
$$

$$
[x_{-c}(u),x_{-b-c-d}(t)] = x_{-b-2c-d}(-N_{c,b+c+d}tu) = x_{-b-2c-d}(-\epsilon_{23}\delta_2 tu).\eqno{(13)}
$$

$$
\begin{array}{rcl}
[x_{-c}(u),x_{-a-b}(t)] & = & x_{-a-b-c}(-N_{c,a+b}tu)\ x_{-a-b-2c}(-(1/2) N_{c,a+b} N_{c,a+b+c}tu^2)  \\[4mm]
                        & = & x_{-a-b-c}(-\delta_3 tu)\ x_{-a-b-2c}(-\delta_{34} tu^2).
\end{array} \eqno{(14)}
$$
Explanation: 
$$
-(1/2) N_{c,a+b} N_{c,a+b+c} = -(1/2) (\delta_3) (2\delta_4) = -\delta_{34}.
$$

$$
[x_{-c}(u),x_{-a-b-c}(t)] = x_{-a-b-2c}(-N_{c,a+b+c}tu) = x_{-a-b-2c}(-2\delta_4 tu).\eqno{(15)}
$$

$$
[x_{-c}(u),x_{-a-b-c-d}(t)] = x_{-a-b-2c-d}(-N_{c,a+b+c+d}tu) = x_{-a-b-2c-d}(-\epsilon_{56}\delta_4 tu).\eqno{(16)}
$$

$$
[x_{-c}(u),x_{-a-2b-2c-d}(t)] = x_{-a-2b-3c-d}(-N_{c,a+2b+2c+d}tu) = x_{-a-2b-3c-d}(-\delta_5 tu).\eqno{(17)}
$$

$$
\begin{array}{rcl}
[x_{-c}(u), x_{-a-2b-2c-2d}(t)]  & = & x_{-a-2b-3c-2d}(-N_{c,a+2b+2c+2d}tu)\times\\[4mm]
      &   & x_{-a-2b-4c-2d}(-(1/2) N_{c,a+2b+2c+2d} N_{c,a+2b+3c+2d} tu^2) \\[4mm]
      & = & x_{-a-2b-3c-2d}(-\epsilon_{90}\delta_5 tu)\ x_{-a-2b-4c-2d}(-\epsilon_{90}\delta_{56} tu^2).
\end{array} \eqno{(18)}
$$
Explanation: 
$$
-(1/2) N_{c,a+2b+2c+2d} N_{c,a+2b+3c+2d} = -(1/2) (\epsilon_{90}\delta_5) (2\delta_6) = -\epsilon_{90}\delta_{56}.
$$

$$
[x_{-c}(u),x_{-a-2b-3c-2d}(t)] = x_{-a-2b-4c-2d}(-N_{c,a+2b+3c+2d}tu) = x_{-a-2b-4c-2d}(-2\delta_6 tu).\eqno{(19)}
$$

{\boldmath$-c-d:$}

$$
\begin{array}{rcl}
[x_{-c-d}(u),x_{-b}(t)] & = & x_{-b-c-d}(-N_{c+d,b} tu)\ x_{-b-2c-2d}(-(1/2) N_{c+d,b} N_{c+d,b+c+d} tu^2) \\[4mm]
                      & = & x_{-b-c-d}(-\epsilon_{12}\delta_1 tu)\ x_{b-2c-2d}(-\epsilon_{34}\delta_{12} tu^2).\\
\end{array} \eqno{(20)}
$$
Explanation: 
$$
-(1/2) N_{c+d,b} N_{c+d,b+c+d} = -(1/2) (\epsilon_{12}\delta_1) (2\epsilon_{1234}\delta_2) = -\epsilon_{34}\delta_{12}.
$$

$$
[x_{-c-d}(u),x_{-b-c-d}(t)] = x_{-b-2c-2d}(-N_{c+d,b+c+d}tu) = x_{-b-2c-2d}(-2\epsilon_{1234}\delta_2 tu).\eqno{(21)}
$$

$$
[x_{-c-d}(u),x_{-b-c}(t)] = x_{-b-2c-d}(-N_{c+d,b+c}tu) = x_{-b-2c-d}(-\epsilon_{13}\delta_2 tu).\eqno{(22)}
$$

$$
\begin{array}{rcl}
[x_{-c-d}(u),x_{-a-b}(t)] & = & x_{-a-b-c-d}(-N_{c+d,a+b} tu)\times\\[4mm]
                          & & x_{-a-b-2c-2d}(-(1/2) N_{c+d,a+b} N_{c+d,a+b+c+d} tu^2) \\[4mm]
                        & = & x_{-a-b-c-d}(-\epsilon_{15}\delta_3 tu)\ x_{-a-b-2c-2d}(-\epsilon_{67}\delta_{34} tu^2).\\
\end{array} \eqno{(23)}
$$
Explanation: 
$$
-(1/2) N_{c+d,a+b} N_{c+d,a+b+c+d} = -(1/2) (\epsilon_{15}\delta_3) (2\epsilon_{1567}\delta_4) = -\epsilon_{67}\delta_{34}.
$$

$$
\begin{array}{rcl}
[x_{-c-d}(u),x_{-a-b-c-d}(t)] & = & x_{-a-b-2c-2d}(-N_{c+d,a+b+c+d} tu) \\[4mm] 
                              & = & x_{-a-b-2c-2d}(-2\epsilon_{1567}\delta_4 tu).
\end{array} \eqno{(24)}
$$

$$
[x_{-c-d}(u),x_{-a-b-c}(t)] = x_{-a-b-2c-d}(-N_{c+d,a+b+c} tu) = x_{-a-b-2c-d}(-\epsilon_{16}\delta_4 tu).\eqno{(25)}
$$

$$
\begin{array}{rcl}
[x_{-c-d}(u),x_{-a-2b-2c}(t)] & = & x_{-a-2b-3c-d}(-N_{c+d,a+2b+2c} tu)\times\\[4mm]
                            & & x_{-a-2b-4c-2d}(-(1/2) N_{c+d,a+2b+2c} N_{c+d,a+2b+3c+d} tu^2) \\[4mm]
                            & = & x_{-a-2b-3c-d}(\epsilon_{18}\delta_5 tu)\ x_{-a-2b-4c-2d}(-\epsilon_{80}\delta_{56} tu^2).\\
\end{array} \eqno{(26)}
$$
Explanation: 
$$
-(1/2) N_{c+d,a+2b+2c} N_{c+d,a+2b+3c+d}  = -(1/2) (-\epsilon_{18}\delta_5) (-2\epsilon_{10}\delta_6) = -\epsilon_{80}\delta_{56}.
$$

$$
\begin{array}{rcl}
[x_{-c-d}(u),x_{-a-2b-3c-d}(t)] & = & x_{-a-2b-4c-2d}(-N_{c+d,a+2b+3c+d} tu) \\[4mm]
                                & = & x_{-a-2b-4c-2d}(2\epsilon_{10}\delta_6 tu).
\end{array} \eqno{(27)}
$$

$$
\begin{array}{rcl}
[x_{-c-d}(u),x_{-a-2b-2c-d}(t)] & = & x_{-a-2b-3c-2d}(-N_{c+d,a+2b+2c+d} tu) \\[4mm] 
                                & = & x_{-a-2b-3c-2d}(\epsilon_{10}\delta_5 tu).
\end{array} \eqno{(28)}
$$

{\boldmath$-b:$}

$$
[x_{-b}(u),x_{-a}(t)] = x_{-a-b}(-N_{b,a} tu) = x_{-a-b}(-\gamma_1 tu).\eqno{(29)}
$$

$$
[x_{-b}(u),x_{-a-b-2c}(t)] = x_{-a-2b-2c}(-N_{b,a+b+2c} tu) = x_{-a-2b-2c}(-\gamma_2 tu).\eqno{(30)}
$$

$$
\begin{array}{rcl}
[x_{-b}(u),x_{-a-b-2c-d}(t)] & = & x_{-a-2b-2c-d}(-N_{b,a+b+2c+d} tu) \\[4mm] 
                             & = & x_{-a-2b-2c-d}(-\epsilon_{68}\gamma_2 tu).
\end{array}\eqno{(31)}
$$

$$
\begin{array}{rcl}
[x_{-b}(u),x_{-a-b-2c-2d}(t)] & = & x_{-a-2b-2c-2d}(-N_{b,a+b+2c+2d} tu) \\[4mm]
                              & = & x_{-a-2b-2c-2d}(-\epsilon_{6789}\gamma_2 tu).
\end{array}\eqno{(32)}
$$

$$
[x_{-b}(u),x_{-a-2b-4c-2d}(t)] = x_{-a-3b-4c-2d}(-N_{b,a+2b+4c+2d} tu) = x_{-a-3b-4c-2d}(-\gamma_3 tu).\eqno{(33)}
$$

{\boldmath$-b-c:$}

$$
\begin{array}{rcl}
[x_{-b-c}(u),x_{-a}(t)] & = & x_{-a-b-c}(-N_{b+c,a} tu)\ x_{-a-2b-2c}(-(1/2) N_{b+c,a} N_{b+c,a+b+c} tu^2) \\[4mm]
                    & = & x_{-a-b-c}(-\delta_{13}\gamma_1 tu)\ x_{-a-2b-2c}(\delta_{34}\gamma_{12} tu^2).\\
\end{array} \eqno{(34)}
$$
Explanation: 
$$
-(1/2) N_{b+c,a} N_{b+c,a+b+c}  = -(1/2) (\delta_{13}\gamma_1) (-2\delta_{14}\gamma_2) = \delta_{34}\gamma_{12}.
$$

$$
[x_{-b-c}(u),x_{-a-b-c}(t)] = x_{-a-2b-2c}(-N_{b+c,a+b+c} tu) = x_{-a-2b-2c}(2\delta_{14}\gamma_2 tu).\eqno{(35)}
$$

$$
[x_{-b-c}(u),x_{-a-b-c-d}(t)] = x_{-a-2b-2c-d}(-N_{b+c,a+b+c+d} tu) = x_{-a-2b-2c-d}(\epsilon_{58}\delta_{14}\gamma_2 tu).\eqno{(36)}
$$

$$
\begin{array}{rcl}
[x_{-b-c}(u),x_{-a-b-2c-d}(t)] & = & x_{-a-2b-3c-d}(-N_{b+c,a+b+2c+d} tu) \\[4mm]
                               & = & x_{-a-2b-3c-d}(-\epsilon_{68}\delta_{15}\gamma_2 tu).
\end{array} \eqno{(37)}
$$

$$
\begin{array}{rcl}
[x_{-b-c}(u),x_{-a-b-2c-2d}(t)] & = & x_{-a-2b-3c-2d}(-N_{b+c,a+b+2c+2d} tu)\times\\[4mm]
                              & & x_{-a-3b-4c-2d}(-(1/2) N_{b+c,a+b+2c+2d} N_{b+c,a+2b+3c+2d} tu^2) \\[4mm]
                              & = & x_{-a-2b-3c-2d}(-\epsilon_{6780}\delta_{15}\gamma_2 tu)\ x_{-a-3b-4c-2d}(\epsilon_{6780}\delta_{56}\gamma_{23} tu^2).\\
\end{array} \eqno{(38)}
$$
Explanation: 
$$
-(1/2) N_{b+c,a+b+2c+2d} N_{b+c,a+2b+3c+2d}  = -(1/2) (\epsilon_{6780}\delta_{15}\gamma_2) (-2\delta_{16}\gamma_3) = \epsilon_{6780}\delta_{56}\gamma_{23}.
$$

$$
\begin{array}{rcl}
[x_{-b-c}(u),x_{-a-2b-3c-2d}(t)] & = & x_{-a-3b-4c-2d}(-N_{b+c,a+2b+3c+2d} tu) \\[4mm]
                                 & = & x_{-a-3b-4c-2d}(2\delta_{16}\gamma_3 tu).
\end{array} \eqno{(39)}
$$

{\boldmath$-b-c-d:$}

$$
\begin{array}{rcl}
[x_{-b-c-d}(u),x_{-a}(t)] & = & x_{-a-b-c-d}(-N_{b+c+d,a} tu)\times\\[4mm]
                          &   & x_{-a-2b-2c-2d}(-(1/2) N_{b+c+d,a} N_{b+c+d,a+b+c+d} tu^2) \\[4mm]
                          & = & x_{-a-b-c-d}(-\epsilon_{25}\delta_{13}\gamma_1 tu)\ x_{-a-2b-2c-2d}(\epsilon_{89}\delta_{34}\gamma_{12} tu^2).\\
\end{array} \eqno{(40)}
$$
Explanation: 
$$
-(1/2) N_{b+c+d,a} N_{b+c+d,a+b+c+d}  = -(1/2) (\epsilon_{25}\delta_{13}\gamma_1) (-2\epsilon_{2589}\delta_{14}\gamma_2) = \epsilon_{89}\delta_{34}\gamma_{12}.
$$

$$
\begin{array}{rcl}
[x_{-b-c-d}(u),x_{-a-b-c-d}(t)] & = & x_{-a-2b-2c-2d}(-N_{b+c+d,a+b+c+d} tu) \\[4mm]
                                & = & x_{-a-2b-2c-2d}(2\epsilon_{2589}\delta_{14}\gamma_2 tu).
\end{array} \eqno{(41)}
$$

$$
[x_{-b-c-d}(u),x_{-a-b-c}(t)] = x_{-a-2b-2c-d}(-N_{b+c+d,a+b+c} tu) = x_{-a-2b-2c-d}(\epsilon_{28}\delta_{14}\gamma_2 tu).\eqno{(42)}
$$

$$
\begin{array}{rcl}
[x_{-b-c-d}(u),x_{-a-b-2c}(t)] & = & x_{-a-2b-3c-d}(-N_{b+c+d,a+b+2c} tu)\times\\[4mm] 
                             & & x_{-a-3b-4c-2d}(-(1/2) N_{b+c+d,a+b+2c} N_{b+c+d,a+2b+3c+d} tu^2) \\[4mm]
                             & = & x_{-a-2b-3c-d}(\epsilon_{28}\delta_{15}\gamma_2 tu)\ x_{-a-3b-4c-2d}(\epsilon_{80}\delta_{56}\gamma_{23} tu^2).\\
\end{array} \eqno{(43)}
$$
Explanation: 
$$
-(1/2) N_{b+c+d,a+b+2c} N_{b+c+d,a+2b+3c+d}  = -(1/2) (-\epsilon_{28}\delta_{15}\gamma_2) (2\epsilon_{20}\delta_{16}\gamma_3) = \epsilon_{80}\delta_{56}\gamma_{23}.
$$

$$
\begin{array}{rcl}
[x_{-b-c-d}(u),x_{-a-2b-3c-d}(t)] & = & x_{-a-3b-4c-2d}(-N_{b+c+d,a+2b+3c+d} tu) \\[4mm]
                                  & = & x_{-a-3b-4c-2d}(-2\epsilon_{20}\delta_{16}\gamma_3 tu).
\end{array} \eqno{(44)}
$$

$$
\begin{array}{rcl}
[x_{-b-c-d}(u),x_{-a-b-2c-d}(t)] & = & x_{-a-2b-3c-2d}(-N_{b+c+d,a+b+2c+d} tu) \\[4mm] 
                                 & = & x_{-a-2b-3c-2d}(\epsilon_{2680}\delta_{15}\gamma_2 tu).
\end{array} \eqno{(45)}
$$

{\boldmath$-b-2c:$}

$$
[x_{-b-2c}(u),x_{-a}(t)] = x_{-a-b-2c}(-N_{b+2c,a} tu) = x_{-a-b-2c}(-\delta_{1234}\gamma_1 tu).\eqno{(46)}
$$

$$
[x_{-b-2c}(u),x_{-a-b}(t)] = x_{-a-2b-2c}(-N_{b+2c,a+b} tu) = x_{-a-2b-2c}(-\delta_{1234}\gamma_2 tu).\eqno{(47)}
$$

$$
\begin{array}{rcl}
[x_{-b-2c}(u),x_{-a-b-c-d}(t)] & = & x_{-a-2b-3c-d}(-N_{b+2c,a+b+c+d} tu) \\[4mm]
                               & = & x_{-a-2b-3c-d}(\epsilon_{58}\delta_{1245}\gamma_2 tu).
\end{array} \eqno{(48)}
$$

$$
\begin{array}{rcl}
[x_{-b-2c}(u),x_{-a-b-2c-2d}(t)] & = & x_{-a-2b-4c-2d}(-N_{b+2c,a+b+2c+2d} tu) \\[4mm] 
                                 & = & x_{-a-2b-4c-2d}(-\epsilon_{6780}\delta_{1256}\gamma_2 tu).
\end{array} \eqno{(49)}
$$

$$
\begin{array}{rcl}
[x_{-b-2c}(u),x_{-a-2b-2c-2d}(t)] & = & x_{-a-3b-4c-2d}(-N_{b+2c,a+2b+2c+2d} tu) \\[4mm]
                                  & = & x_{-a-3b-4c-2d}(-\epsilon_{90}\delta_{1256}\gamma_3 tu).
\end{array} \eqno{(50)}
$$

{\boldmath$-b-2c-d:$}

$$
\begin{array}{rcl}
[x_{-b-2c-d}(u),x_{-a}(t)] & = & x_{-a-b-2c-d}(-N_{b+2c+d,a} tu)\times\\[4mm]
                         & & x_{-a-2b-4c-2d}(-(1/2) N_{b+2c+d,a} N_{b+2c+d,a+b+2c+d} tu^2) \\[4mm]
												 & = & x_{-a-b-2c-d}(-\epsilon_{36}\delta_{1234}\gamma_1 tu)\ x_{-a-2b-4c-2d}(\epsilon_{80}\delta_{3456}\gamma_{12} tu^2).
\end{array} \eqno{(51)}
$$
Explanation: 
$$
-(1/2)  N_{b+2c+d,a} N_{b+2c+d,a+b+2c+d}  = -(1/2) (\epsilon_{36}\delta_{1234}\gamma_1) (-2\epsilon_{3680}\delta_{1256}\gamma_2) = \epsilon_{80}\delta_{3456}\gamma_{12}.
$$

$$
\begin{array}{rcl}
[x_{-b-2c-d}(u),x_{-a-b-2c-d}(t)] & = & x_{-a-2b-4c-2d}(-N_{b+2c+d,a+b+2c+d} tu) \\[4mm]
                                  & = & x_{-a-2b-4c-2d}(2\epsilon_{3680}\delta_{1256}\gamma_2 tu).
\end{array} \eqno{(52)}
$$

$$
\begin{array}{rcl}
[x_{-b-2c-d}(u),x_{-a-b}(t)] & = & x_{-a-2b-2c-d}(-N_{b+2c+d,a+b} tu)\times\\[4mm] 
                             & & x_{-a-3b-4c-2d}(-(1/2) N_{b+2c+d,a+b} N_{b+2c+d,a+2b+2c+d} tu^2) \\[4mm]
                             & = & x_{-a-2b-2c-d}(-\epsilon_{38}\delta_{1234}\gamma_2 tu)\ x_{-a-3b-4c-2d}(\epsilon_{80}\delta_{3456}\gamma_{23} tu^2).\\
\end{array} \eqno{(53)}
$$
Explanation: 
$$
-(1/2) N_{b+2c+d,a+b} N_{b+2c+d,a+2b+2c+d}  = -(1/2) (\epsilon_{38}\delta_{1234}\gamma_2) (-2\epsilon_{30}\delta_{1256}\gamma_3) = \epsilon_{80}\delta_{3456}\gamma_{23}.
$$

$$
\begin{array}{rcl}
[x_{-b-2c-d}(u),x_{-a-2b-2c-d}(t)] & = & x_{-a-3b-4c-2d}(-N_{b+2c+d,a+2b+2c+d} tu) \\[4mm]
                                   & = & x_{-a-3b-4c-2d}(2\epsilon_{30}\delta_{1256}\gamma_3 tu).
\end{array} \eqno{(54)}
$$

$$
\begin{array}{rcl}
[x_{-b-2c-d}(u),x_{-a-b-c}(t)] & = & x_{-a-2b-3c-d}(-N_{b+2c+d,a+b+c} tu) \\[4mm]
                               & = & x_{-a-2b-3c-d}(-\epsilon_{38}\delta_{1245}\gamma_2 tu).
\end{array} \eqno{(55)}
$$

$$
\begin{array}{rcl}
[x_{-b-2c-d}(u),x_{-a-b-c-d}(t)] & = & x_{-a-2b-3c-2d}(-N_{b+2c+d,a+b+c+d} tu) \\[4mm]
                                 & = & x_{-a-2b-3c-2d}(\epsilon_{3580}\delta_{1245}\gamma_2 tu).
\end{array} \eqno{(56)}
$$

{\boldmath$-b-2c-2d:$}

$$
\begin{array}{rcl}
[x_{-b-2c-2d}(u),x_{-a}(t)] & = & x_{-a-b-2c-2d}(-N_{b+2c+2d,a} tu) \\[4mm] 
                            & = & x_{-a-b-2c-2d}(-\epsilon_{3467}\delta_{1234}\gamma_1 tu).
\end{array} \eqno{(57)}
$$

$$
\begin{array}{rcl}
[x_{-b-2c-2d}(u),x_{-a-b}(t)] & = & x_{-a-2b-2c-2d}(-N_{b+2c+2d,a+b} tu) \\[4mm]
                              & = & x_{-a-2b-2c-2d}(-\epsilon_{3489}\delta_{1234}\gamma_2 tu).
\end{array} \eqno{(58)}
$$

$$
\begin{array}{rcl}
[x_{-b-2c-2d}(u),x_{-a-b-c}(t)] & = & x_{-a-2b-3c-2d}(-N_{b+2c+2d,a+b+c} tu) \\[4mm]
                                & = & x_{-a-2b-3c-2d}(-\epsilon_{3480}\delta_{1245}\gamma_2 tu).
\end{array} \eqno{(59)}
$$

$$
\begin{array}{rcl}
[x_{-b-2c-2d}(u),x_{-a-b-2c}(t)] & = & x_{-a-2b-4c-2d}(-N_{b+2c+2d,a+b+2c} tu) \\[4mm]
                                 & = & x_{-a-2b-4c-2d}(-\epsilon_{3480}\delta_{1256}\gamma_2 tu).
\end{array} \eqno{(60)}
$$

$$
\begin{array}{rcl}
[x_{-b-2c-2d}(u),x_{-a-2b-2c}(t)] & = & x_{-a-3b-4c-2d}(-N_{b+2c+2d,a+2b+2c} tu) \\[4mm]
                                  & = & x_{-a-3b-4c-2d}(-\epsilon_{3480}\delta_{1256}\gamma_3 tu).
\end{array} \eqno{(61)}
$$

{\boldmath$-a:$}

$$
\begin{array}{rcl}
[x_{-a}(u),x_{-a-3b-4c-2d}(t)] & = & x_{-2a-3b-4c-2d}(-N_{a,a+3b+4c+2d} tu) \\[4mm]
                               & = & x_{-2a-3b-4c-2d}(-\alpha_1 tu).
\end{array} \eqno{(62)}
$$

$\bold{-a-b:}$

$$
\begin{array}{rcl}
[x_{-a-b}(u),x_{-a-2b-4c-2d}(t)] & = & x_{-2a-3b-4c-2d}(-N_{a+b,a+2b+4c+2d} tu) \\[4mm]
                                 & = & x_{-2a-3b-4c-2d}(\gamma_{13}\alpha_1 tu).
\end{array} \eqno{(63)}
$$

{\boldmath$-a-b-c:$}

$$
\begin{array}{rcl}
[x_{-a-b-c}(u),x_{-a-2b-3c-2d}(t)] & = & x_{-2a-3b-4c-2d}(-N_{a+b+c,a+2b+3c+2d} tu) \\[4mm]
                                   & = & x_{-2a-3b-4c-2d}(-2\delta_{36}\gamma_{13}\alpha_1 tu).
\end{array} \eqno{(64)}
$$

{\boldmath$-a-b-c-d:$}

$$
\begin{array}{rcl}
[x_{-a-b-c-d}(u),x_{-a-2b-3c-d}(t)] & = & x_{-2a-3b-4c-2d}(-N_{a+b+c+d,a+2b+3c+d} tu) \\[4mm]
                                    & = & x_{-2a-3b-4c-2d}(2\epsilon_{50}\delta_{36}\gamma_{13}\alpha_1 tu).
\end{array} \eqno{(65)}
$$

{\boldmath$-a-b-2c:$}

$$
\begin{array}{rcl}
[x_{-a-b-2c}(u),x_{-a-2b-2c-2d}(t)] & = & x_{-2a-3b-4c-2d}(-N_{a+b+2c,a+2b+2c+2d} tu) \\[4mm]
                                    & = & x_{-2a-3b-4c-2d}(\epsilon_{90}\delta_{3456}\gamma_{13}\alpha_1 tu).
\end{array} \eqno{(66)}
$$

{\boldmath$-a-b-2c-d:$}

$$
\begin{array}{rcl}
[x_{-a-b-2c-d}(u),x_{-a-2b-2c-d}(t)] & = & x_{-2a-3b-4c-2d}(-N_{a+b+2c+d,a+2b+2c+d} tu) \\[4mm]
                                     & = & x_{-2a-3b-4c-2d}(-2\epsilon_{60}\delta_{3456}\gamma_{13}\alpha_1 tu).
\end{array} \eqno{(67)}
$$

{\boldmath$-a-b-2c-2d:$}

$$
\begin{array}{rcl}
[x_{-a-b-2c-2d}(u),x_{-a-2b-2c}(t)] & = & x_{-2a-3b-4c-2d}(-N_{a+b+2c+2d,a+2b+2c} tu) \\[4mm]
                                    & = & x_{-2a-3b-4c-2d}(\epsilon_{6780}\delta_{3456}\gamma_{13}\alpha_1 tu).
\end{array} \eqno{(68)}
$$

\vskip5mm
\begin{center}
\textit{Positive and negative roots}
\end{center}

{\boldmath$d:$}

$$
[x_d(u),x_{-c-d}(t)] = x_{-c}(N_{d,-c-d} tu) = x_{-c}(-\epsilon_1 tu).\eqno{(001)}
$$

$$
[x_d(u),x_{-b-c-d}(t)] = x_{-b-c}(N_{d,-b-c-d} tu) = x_{-b-c}(-\epsilon_2 tu).\eqno{(002)}
$$

$$
[x_d(u),x_{-b-2c-d}(t)] = x_{-b-2c}(N_{d,-b-2c-d} tu) = x_{-b-2c}(-2\epsilon_3 tu).\eqno{(003)}
$$

$$
\begin{array}{rcl}
[x_d(u),x_{-b-2c-2d}(t)] & = & x_{-b-2c-d}(N_{d,-b-2c-2d} tu)\times \\[4mm]
                         &   & x_{-b-2c}(-(1/2) N_{d,-b-2c-2d} N_{d,-b-2c-d}  tu^2)\\[4mm]
                         & = & x_{-b-2c-d}(-\epsilon_4 tu)\ x_{-b-2c}(-\epsilon_{34} tu^2).
\end{array} \eqno{(004)}
$$
Explanation: 
$$
-(1/2) N_{d,-b-2c-2d} N_{d,-b-2c-d}  = -(1/2) (-\epsilon_4) (-2\epsilon_3) = -\epsilon_{34}.
$$

$$
[x_d(u),x_{-a-b-c-d}(t)] = x_{-a-b-c}(N_{d,-a-b-c-d} tu) = x_{-a-b-c}(-\epsilon_5 tu).\eqno{(005)}
$$

$$
[x_d(u),x_{-a-b-2c-d}(t)] = x_{-a-b-2c}(N_{d,-a-b-2c-d} tu) = x_{-a-b-2c}(-2\epsilon_6 tu).\eqno{(006)}
$$

$$
\begin{array}{rcl}
[x_d(u),x_{-a-b-2c-2d}(t)] & = & x_{-a-b-2c-d}(N_{d,-a-b-2c-2d} tu)\times\\[4mm]
                           & & x_{-a-b-2c}(-(1/2) N_{d,-a-b-2c-2d} N_{d,-a-b-2c-d}  tu^2)\\[4mm]
                           & = & x_{-a-b-2c-d}(-\epsilon_7 tu)\ x_{-a-b-2c}(-\epsilon_{67} tu^2).
\end{array} \eqno{(007)}
$$
Explanation: 
$$
-(1/2) N_{d,-a-b-2c-2d} N_{d,-a-b-2c-d}  = -(1/2) (-\epsilon_7) (-2\epsilon_6) = -\epsilon_{67}.
$$

$$
[x_d(u),x_{-a-2b-2c-d}(t)] = x_{-a-2b-2c}(N_{d,-a-2b-2c-d} tu) = x_{-a-2b-2c}(-2\epsilon_8 tu).\eqno{(008)}
$$

$$
\begin{array}{rcl}
[x_d(u),x_{-a-2b-2c-2d}(t)] & = & x_{-a-2b-2c-d}(N_{d,-a-2b-2c-2d} tu)\times\\[4mm]
                            & & x_{-a-2b-2c}(-(1/2) N_{d,-a-2b-2c-2d} N_{d,-a-2b-2c-d}  tu^2)\\[4mm]
                            & = & x_{-a-2b-2c-d}(-\epsilon_9 tu)\ x_{-a-2b-2c}(-\epsilon_{89} tu^2).
\end{array} \eqno{(009)}
$$
Explanation: 
$$
-(1/2) N_{d,-a-2b-2c-2d} N_{d,-a-2b-2c-d}  = -(1/2) (-\epsilon_9) (-2\epsilon_8) = -\epsilon_{89}.
$$

$$
[x_d(u),x_{-a-2b-3c-2d}(t)] = x_{-a-2b-3c-d}(N_{d,-a-2b-3c-2d} tu) = x_{-a-2b-3c-d}(-\epsilon_0 tu).\eqno{(010)}
$$

{\boldmath$c:$}

$$
[x_c(u),x_{-c-d}(t)] = x_{-d}(N_{c,-c-d} tu) = x_{-d}(\epsilon_1 tu).\eqno{(011)}
$$

$$
[x_c(u),x_{-b-c}(t)] = x_{-b}(N_{c,-b-c} tu) = x_{-b}(-2\delta_1 tu).\eqno{(012)}
$$

$$
\begin{array}{rcl}
[x_c(u),x_{-b-2c}(t)] & = & x_{-b-c}(N_{c,-b-2c} tu)\ x_{-b}(-(1/2) N_{c,-b-2c} N_{c,-b-c} tu^2) \\[4mm]
                      & = & x_{-b-c}(-\delta_2 tu)\ x_{-b}(-\delta_{12} tu^2).
\end{array} \eqno{(013)}
$$
Explanation: 
$$
-(1/2) N_{c,-b-2c} N_{c,-b-c}  = -(1/2) (-\delta_2) (-2\delta_1) = -\delta_{12}.
$$

$$
[x_c(u),x_{-b-2c-d}(t)] = x_{-b-c-d}(N_{c,-b-2c-d} tu) = x_{-b-c-d}(-\epsilon_{23}\delta_2 tu).\eqno{(014)}
$$

$$
[x_c(u),x_{-a-b-c}(t)] = x_{-a-b}(N_{c,-a-b-c} tu) = x_{-a-b}(-2\delta_3 tu).\eqno{(015)}
$$

$$
\begin{array}{rcl}
[x_c(u),x_{-a-b-2c}(t)] & = & x_{-a-b-c}(N_{c,-a-b-2c} tu)\ x_{-a-b}(-(1/2) N_{c,-a-b-2c} N_{c,-a-b-c} tu^2) \\[4mm]
                        & = & x_{-a-b-c}(-\delta_4 tu)\ x_{-a-b}(-\delta_{34} tu^2).
\end{array} \eqno{(016)}
$$
Explanation: 
$$
-(1/2) N_{c,-a-b-2c} N_{c,-a-b-c}  = -(1/2) (-\delta_4) (-2\delta_3) = -\delta_{34}.
$$

$$
[x_c(u),x_{-a-b-2c-d}(t)] = x_{-a-b-c-d}(N_{c,-a-b-2c-d} tu) = x_{-a-b-c-d}(-\epsilon_{56}\delta_4 tu).\eqno{(017)}
$$

$$
[x_c(u),x_{-a-2b-3c-d}(t)] = x_{-a-2b-2c-d}(N_{c,-a-2b-3c-d} tu) = x_{-a-2b-2c-d}(-\delta_5 tu).\eqno{(018)}
$$

$$
\begin{array}{rcl}
[x_c(u),x_{-a-2b-3c-2d}(t)] & = & x_{-a-2b-2c-2d}(N_{c,-a-2b-3c-2d} tu) \\[4mm]
                            & = & x_{-a-2b-2c-2d}(-2\epsilon_{90}\delta_5 tu).
\end{array} \eqno{(019)}
$$

$$
\begin{array}{rcl}
[x_c(u),x_{-a-2b-4c-2d}(t)] & = & x_{-a-2b-3c-2d}(N_{c,-a-2b-4c-2d} tu)\times\\[4mm]
                            &   & x_{-a-2b-2c-2d}(-(1/2) N_{c,-a-2b-4c-2d} N_{c,-a-2b-3c-2d} tu^2) \\[4mm]
                            & = & x_{-a-2b-3c-2d}(-\delta_6 tu)\ x_{-a-2b-2c-2d}(-\epsilon_{90}\delta_{56} tu^2).
\end{array} \eqno{(020)}
$$
Explanation: 
$$
-(1/2) N_{c,-a-2b-4c-2d} N_{c,-a-2b-3c-2d}  = -(1/2) (-\delta_6) (-2\epsilon_{90}\delta_5) = -\epsilon_{90}\delta_{56}.
$$

{\boldmath$c+d:$}

$$
[x_{c+d}(u),x_{-b-c-d}(t)] = x_{-b}(N_{c+d,-b-c-d} tu) = x_{-b}(-2\epsilon_{12}\delta_1 tu).\eqno{(021)}
$$

$$
\begin{array}{rcl}
[x_{c+d}(u),x_{-b-2c-2d}(t)] & = & x_{-b-c-d}(N_{c+d,-b-2c-2d} tu)\times \\[4mm]
                             &   & x_{-b}(-(1/2) N_{c+d,-b-2c-2d} N_{c+d,-b-c-d} tu^2) \\[4mm]
                             & = & x_{-b-c-d}(-\epsilon_{1234}\delta_2 tu)\ x_{-b}(-\epsilon_{34}\delta_{12} tu^2).
\end{array} \eqno{(022)}
$$
Explanation: 
$$
-(1/2) N_{c+d,-b-2c-2d} N_{c+d,-b-c-d}  = -(1/2) (-\epsilon_{1234}\delta_2) (-2\epsilon_{12}\delta_1) = -\epsilon_{34}\delta_{12}.
$$

$$
[x_{c+d}(u),x_{-b-2c-d}(t)] = x_{-b-c}(N_{c+d,-b-2c-d} tu) = x_{-b-c}(-\epsilon_{13}\delta_2 tu).\eqno{(023)}
$$

$$
[x_{c+d}(u),x_{-a-b-c-d}(t)] = x_{-a-b}(N_{c+d,-a-b-c-d} tu) = x_{-a-b}(-2\epsilon_{15}\delta_3 tu).\eqno{(024)}
$$

$$
\begin{array}{rcl}
[x_{c+d}(u),x_{-a-b-2c-2d}(t)] & = & x_{-a-b-c-d}(N_{c+d,-a-b-2c-2d} tu)\times\\[4mm]
                               &   & x_{-a-b}(-(1/2) N_{c+d,-a-b-2c-2d} N_{c+d,-a-b-c-d} tu^2) \\[4mm]
                               & = & x_{-a-b-c-d}(-\epsilon_{1567}\delta_4 tu)\ x_{-a-b}(-\epsilon_{67}\delta_{34} tu^2).
\end{array} \eqno{(025)}
$$
Explanation: 
$$
-(1/2) N_{c+d,-a-b-2c-2d} N_{c+d,-a-b-c-d}  = -(1/2) (-\epsilon_{1567}\delta_4) (-2\epsilon_{15}\delta_3) = -\epsilon_{67}\delta_{34}.
$$

$$
[x_{c+d}(u),x_{-a-b-2c-d}(t)] = x_{-a-b-c}(N_{c+d,-a-b-2c-d} tu) = x_{-a-b-c}(-\epsilon_{16}\delta_4 tu).\eqno{(026)}
$$

$$
[x_{c+d}(u),x_{-a-2b-3c-d}(t)] = x_{-a-2b-2c}(N_{c+d,-a-2b-3c-d} tu) = x_{-a-2b-2c}(2\epsilon_{18}\delta_5 tu).\eqno{(027)}
$$

$$
\begin{array}{rcl}
[x_{c+d}(u),x_{-a-2b-4c-2d}(t)] & = & x_{-a-2b-3c-d}(N_{c+d,-a-2b-4c-2d} tu)\times\\[4mm]
                                & & x_{-a-2b-2c}(-(1/2) N_{c+d,-a-2b-4c-2d} N_{c+d,-a-2b-3c-d} tu^2) \\[4mm]
                                & = & x_{-a-2b-3c-d}(\epsilon_{10}\delta_6 tu)\ x_{-a-2b-2c}(-\epsilon_{80}\delta_{56} tu^2).
\end{array} \eqno{(028)}
$$
Explanation: 
$$
-(1/2) N_{c+d,-a-2b-4c-2d} N_{c+d,-a-2b-3c-d}  = -(1/2) (\epsilon_{10}\delta_6) (2\epsilon_{18}\delta_5) = -\epsilon_{80}\delta_{56}.
$$

$$
\begin{array}{rcl}
[x_{c+d}(u),x_{-a-2b-3c-2d}(t)] & = & x_{-a-2b-2c-d}(N_{c+d,-a-2b-3c-2d} tu) \\[4mm]
                                & = & x_{-a-2b-2c-d}(\epsilon_{10}\delta_5 tu).
\end{array}\eqno{(029)}
$$

{\boldmath$b:$}

$$
[x_b(u),x_{-b-c}(t)] = x_{-c}(N_{b,-b-c} tu) = x_{-c}(\delta_1 tu).\eqno{(030)}
$$

$$
[x_b(u),x_{-a-b}(t)] = x_{-a}(N_{b,-a-b} tu) = x_{-a}(-\gamma_1 tu).\eqno{(031)}
$$

$$
[x_b(u),x_{-b-c-d}(t)] = x_{-c-d}(N_{b,-b-c-d} tu) = x_{-c-d}(\epsilon_{12}\delta_1 tu).\eqno{(032)}
$$

$$
\begin{array}{rcl}
[x_b(u),x_{-a-2b-2c}(t)] & = & x_{-a-b-2c}(N_{b,-a-2b-2c} tu) \\[4mm] 
                         & = & x_{-a-b-2c}(-\gamma_2 tu).
\end{array} \eqno{(033)}
$$

$$
\begin{array}{rcl}
[x_b(u),x_{-a-2b-2c-d}(t)] & = & x_{-a-b-2c-d}(N_{b,-a-2b-2c-d} tu) \\[4mm]
                           & = & x_{-a-b-2c-d}(-\epsilon_{68}\gamma_2 tu).
\end{array} \eqno{(034)}
$$

$$
\begin{array}{rcl}
[x_b(u),x_{-a-2b-2c-2d}(t)] & = & x_{-a-b-2c-2d}(N_{b,-a-2b-2c-2d} tu) \\[4mm] 
                            & = & x_{-a-b-2c-2d}(-\epsilon_{6789}\gamma_2 tu).
\end{array} \eqno{(035)}
$$

$$
\begin{array}{rcl}
[x_b(u),x_{-a-3b-4c-2d}(t)] & = & x_{-a-2b-4c-2d}(N_{b,-a-3b-4c-2d} tu) \\[4mm] 
                            & = & x_{-a-2b-4c-2d}(-\gamma_3 tu).
\end{array} \eqno{(036)}
$$

{\boldmath$b+c:$}

$$
[x_{b+c}(u),x_{-b-c-d}(t)] = x_{-d}(N_{b+c,-b-c-d} tu) = x_{-d}(\epsilon_2 tu).\eqno{(037)}
$$

$$
[x_{b+c}(u),x_{-b-2c}(t)] = x_{-c}(N_{b+c,-b-2c} tu) = x_{-c}(\delta_2 tu).\eqno{(038)}
$$

$$
[x_{b+c}(u),x_{-a-b-c}(t)] = x_{-a}(N_{b+c,-a-b-c} tu) = x_{-a}(-2\delta_{13}\gamma_1 tu).\eqno{(039)}
$$

$$
\begin{array}{rcl}
[x_{b+c}(u),x_{-a-2b-2c}(t)] & = & x_{-a-b-c}(N_{b+c,-a-2b-2c} tu)\times\\[4mm]
                                & & x_{-a}(-(1/2) N_{b+c,-a-2b-2c} N_{b+c,-a-b-c} tu^2) \\[4mm]
                                & = & x_{-a-b-c}(\delta_{14}\gamma_2 tu)\ x_{-a}(\delta_{34}\gamma_{12} tu^2).
\end{array} \eqno{(040)}
$$
Explanation: 
$$
-(1/2) N_{b+c,-a-2b-2c} N_{b+c,-a-b-c}  = -(1/2) (\delta_{14}\gamma_2) (-2\delta_{13}\gamma_1) = \delta_{34}\gamma_{12}.
$$

$$
[x_{b+c}(u),x_{-b-2c-d}(t)] = x_{-c-d}(N_{b+c,-b-2c-d} tu) = x_{-c-d}(\epsilon_{13}\delta_2 tu).\eqno{(041)}
$$

$$
\begin{array}{rcl}
[x_{b+c}(u),x_{-a-2b-2c-d}(t)] & = & x_{-a-b-c-d}(N_{b+c,-a-2b-2c-d} tu) \\[4mm]
                               & = & x_{-a-b-c-d}(\epsilon_{58}\delta_{14}\gamma_2 tu).
\end{array} \eqno{(042)}
$$

$$
\begin{array}{rcl}
[x_{b+c}(u),x_{-a-2b-3c-d}(t)] & = & x_{-a-b-2c-d}(N_{b+c,-a-2b-3c-d} tu) \\[4mm]
                               & = & x_{-a-b-2c-d}(-\epsilon_{68}\delta_{15}\gamma_2 tu).
\end{array} \eqno{(043)}
$$

$$
\begin{array}{rcl}
[x_{b+c}(u),x_{-a-2b-3c-2d}(t)] & = & x_{-a-b-2c-2d}(N_{b+c,-a-2b-3c-2d} tu) \\[4mm]
                                & = & x_{-a-b-2c-2d}(-2\epsilon_{6780}\delta_{15}\gamma_2 tu).
\end{array} \eqno{(044)}
$$

$$
\begin{array}{rcl}
[x_{b+c}(u),x_{-a-3b-4c-2d}(t)] & = & x_{-a-2b-3c-2d}(N_{b+c,-a-3b-4c-2d} tu)\times\\[4mm]
                                & & x_{-a-b-2c-2d}(-(1/2) N_{b+c,-a-3b-4c-2d} N_{b+c,-a-2b-3c-2d} tu^2)  \\[4mm]
                                & = & x_{-a-2b-3c-2d}(\delta_{16}\gamma_3 tu)\ x_{-a-b-2c-2d}(\epsilon_{6780}\delta_{56}\gamma_{23} tu^2).
\end{array} \eqno{(045)}
$$
Explanation: 
$$
-(1/2) N_{b+c,-a-3b-4c-2d} N_{b+c,-a-2b-3c-2d}  = -(1/2) (\delta_{16}\gamma_3) (-2\epsilon_{6780}\delta_{15}\gamma_2) = \epsilon_{6780}\delta_{56}\gamma_{23}.
$$

{\boldmath$b+c+d:$}

$$
[x_{b+c+d}(u),x_{-a-b-c-d}(t)] = x_{-a}(N_{b+c+d,-a-b-c-d} tu) = x_{-a}(-2\epsilon_{25}\delta_{13}\gamma_1 tu).\eqno{(046)}
$$

$$
\begin{array}{rcl}
[x_{b+c+d}(u),x_{-a-2b-2c-2d}(t)] & = & x_{-a-b-c-d}(N_{b+c+d,-a-2b-2c-2d} tu)\times\\[4mm]
                                & & x_{-a}(-(1/2) N_{b+c+d,-a-2b-2c-2d} N_{b+c+d,-a-b-c-d} tu^2) \\[4mm]
                                & = & x_{-a-b-c-d}(\epsilon_{2589}\delta_{14}\gamma_2 tu)\ x_{-a}(\epsilon_{89}\delta_{34}\gamma_{12} tu^2).
\end{array} \eqno{(047)}
$$
Explanation: 
$$
-(1/2) N_{b+c+d,-a-2b-2c-2d} N_{b+c+d,-a-b-c-d}  = -(1/2) (\epsilon_{2589}\delta_{14}\gamma_2) (-2\epsilon_{25}\delta_{13}\gamma_1) = \epsilon_{89}\delta_{34}\gamma_{12}.
$$

$$
[x_{b+c+d}(u),x_{-b-2c-d}(t)] = x_{-c}(N_{b+c+d,-b-2c-d} tu) = x_{-c}(\epsilon_{23}\delta_2 tu).\eqno{(048)}
$$

$$
[x_{b+c+d}(u),x_{-b-2c-2d}(t)] = x_{-c-d}(N_{b+c+d,-b-2c-2d} tu) = x_{-c-d}(\epsilon_{1234}\delta_2 tu).\eqno{(049)}
$$

$$
\begin{array}{rcl}
[x_{b+c+d}(u),x_{-a-2b-2c-d}(t)] & = & x_{-a-b-c}(N_{b+c+d,-a-2b-2c-d} tu) \\[4mm]
                                 & = & x_{-a-b-c}(\epsilon_{28}\delta_{14}\gamma_2 tu).
\end{array} \eqno{(050)}
$$

$$
\begin{array}{rcl}
[x_{b+c+d}(u),x_{-a-2b-3c-d}(t)] & = & x_{-a-b-2c}(N_{b+c+d,-a-2b-3c-d} tu) \\[4mm]
                                 & = & x_{-a-b-2c}(2\epsilon_{28}\delta_{15}\gamma_2 tu).
\end{array} \eqno{(051)}
$$

$$
\begin{array}{rcl}
[x_{b+c+d}(u),x_{-a-3b-4c-2d}(t)] & = & x_{-a-2b-3c-d}(N_{b+c+d,-a-3b-4c-2d} tu)\times\\[4mm]
                                  &   & x_{-a-b-2c}(-(1/2) N_{b+c+d,-a-3b-4c-2d} N_{b+c+d,-a-2b-3c-d} tu^2) \\[4mm]
                                & = & x_{-a-2b-3c-d}(-\epsilon_{20}\delta_{16}\gamma_3 tu)\ x_{-a-b-2c}(\epsilon_{80}\delta_{56}\gamma_{23} tu^2).
\end{array} \eqno{(052)}
$$
Explanation: 
$$
-(1/2)  N_{b+c+d,-a-3b-4c-2d} N_{b+c+d,-a-2b-3c-d}  = -(1/2) (-\epsilon_{20}\delta_{16}\gamma_3) (2\epsilon_{28}\delta_{15}\gamma_2) = \epsilon_{80}\delta_{56}\gamma_{23}.
$$

$$
\begin{array}{rcl}
[x_{b+c+d}(u),x_{-a-2b-3c-2d}(t)] & = & x_{-a-b-2c-d}(N_{b+c+d,-a-2b-3c-2d} tu) \\[4mm] 
                                  & = & x_{-a-b-2c-d}(\epsilon_{2680}\delta_{15}\gamma_2 tu).
\end{array} \eqno{(053)}
$$

{\boldmath$b+2c:$}

$$
[x_{b+2c}(u),x_{-b-2c-d}(t)] = x_{-d}(N_{b+2c,-b-2c-d} tu) = x_{-d}(\epsilon_3 tu).\eqno{(054)}
$$

$$
[x_{b+2c}(u),x_{-a-b-2c}(t)] = x_{-a}(N_{b+2c,-a-b-2c} tu) = x_{-a}(-\delta_{1234}\gamma_1 tu).\eqno{(055)}
$$

$$
[x_{b+2c}(u),x_{-a-2b-2c}(t)] = x_{-a-b}(N_{b+2c,-a-2b-2c} tu) = x_{-a-b}(-\delta_{1234}\gamma_2 tu).\eqno{(056)}
$$

$$
\begin{array}{rcl}
[x_{b+2c}(u),x_{-a-2b-3c-d}(t)] & = & x_{-a-b-c-d}(N_{b+2c,-a-2b-3c-d} tu) \\[4mm]
                                & = & x_{-a-b-c-d}(\epsilon_{58}\delta_{1245}\gamma_2 tu).
\end{array} \eqno{(057)}
$$

$$
\begin{array}{rcl}
[x_{b+2c}(u),x_{-a-2b-4c-2d}(t)] & = & x_{-a-b-2c-2d}(N_{b+2c,-a-2b-4c-2d} tu) \\[4mm]
                                 & = & x_{-a-b-2c-2d}(-\epsilon_{6780}\delta_{1256}\gamma_2 tu).
\end{array}\eqno{(058)}
$$

$$
\begin{array}{rcl}
[x_{b+2c}(u),x_{-a-3b-4c-2d}(t)] & = & x_{-a-2b-2c-2d}(N_{b+2c,-a-3b-4c-2d} tu) \\[4mm]
                                 & = & x_{-a-2b-2c-2d}(-\epsilon_{90}\delta_{1256}\gamma_3 tu).
\end{array}\eqno{(059)}
$$

{\boldmath$b+2c+d:$}

$$
[x_{b+2c+d}(u),x_{-a-b-2c-d}(t)] = x_{-a}(N_{b+2c+d,-a-b-2c-d} tu) = x_{-a}(-2\epsilon_{36}\delta_{1234}\gamma_1 tu).\eqno{(060)}
$$

$$
\begin{array}{rcl}
[x_{b+2c+d}(u),x_{-a-2b-4c-2d}(t)] & = & x_{-a-b-2c-d}(N_{b+2c+d,-a-2b-4c-2d} tu)\times\\[4mm]
                                   &   & x_{-a}(-(1/2) N_{b+2c+d,-a-2b-4c-2d} N_{b+2c+d,-a-b-2c-d} tu^2) \\[4mm]
                                   & = & x_{-a-b-2c-d}(\epsilon_{3680}\delta_{1256}\gamma_2 tu)\ x_{-a}(\epsilon_{80}\delta_{3456}\gamma_{12} tu^2).
\end{array} \eqno{(061)}
$$
Explanation: 
$$
-(1/2)  N_{b+2c+d,-a-2b-4c-2d} N_{b+2c+d,-a-b-2c-d}  = -(1/2) (\epsilon_{3680}\delta_{1256}\gamma_2) (-2\epsilon_{36}\delta_{1234}\gamma_1) 
$$
$$
= \epsilon_{80}\delta_{3456}\gamma_{12}.
$$

$$
[x_{b+2c+d}(u),x_{-b-2c-2d}(t)] = x_{-d}(N_{b+2c+d,-b-2c-2d} tu) = x_{-d}(\epsilon_4 tu).\eqno{(062)}
$$

$$
\begin{array}{rcl}
[x_{b+2c+d}(u),x_{-a-2b-2c-d}(t)] & = & x_{-a-b}(N_{b+2c+d,-a-2b-2c-d} tu) \\[4mm]
                                  & = & x_{-a-b}(-2\epsilon_{38}\delta_{1234}\gamma_2 tu).
\end{array} \eqno{(063)}
$$

$$
\begin{array}{rcl}
[x_{b+2c+d}(u),x_{-a-3b-4c-2d}(t)] & = & x_{-a-2b-2c-d}(N_{b+2c+d,-a-3b-4c-2d} tu)\times\\[4mm]
                                   & & x_{-a-b}(-(1/2) N_{b+2c+d,-a-3b-4c-2d} N_{b+2c+d,-a-2b-2c-d} tu^2) \\[4mm]
                                   & = & x_{-a-2b-2c-d}(\epsilon_{30}\delta_{1256}\gamma_3 tu)\ x_{-a-b}(\epsilon_{80}\delta_{3456}\gamma_{23} tu^2).
\end{array} \eqno{(064)}
$$
Explanation: 
$$
-(1/2)  N_{b+2c+d,-a-3b-4c-2d} N_{b+2c+d,-a-2b-2c-d}  = -(1/2) (\epsilon_{30}\delta_{1256}\gamma_3) (-2\epsilon_{38}\delta_{1234}\gamma_2) 
$$
$$
= \epsilon_{80}\delta_{3456}\gamma_{23}.
$$

$$
\begin{array}{rcl}
[x_{b+2c+d}(u),x_{-a-2b-3c-d}(t)] & = & x_{-a-b-c}(N_{b+2c+d,-a-2b-3c-d} tu) \\[4mm]
                                  & = & x_{-a-b-c}(-\epsilon_{38}\delta_{1245}\gamma_2 tu).
\end{array} \eqno{(065)}
$$

$$
\begin{array}{rcl}
[x_{b+2c+d}(u),x_{-a-2b-3c-2d}(t)] & = & x_{-a-b-c-d}(N_{b+2c+d,-a-2b-3c-2d} tu) \\[4mm]
                                   & = & x_{-a-b-c-d}(\epsilon_{3580}\delta_{1245}\gamma_2 tu).
\end{array} \eqno{(066)}
$$

{\boldmath$b+2c+2d:$}

$$
\begin{array}{rcl}
[x_{b+2c+2d}(u),x_{-a-b-2c-2d}(t)] & = & x_{-a}(N_{b+2c+2d,-a-b-2c-2d} tu) \\[4mm]
                                   & = & x_{-a}(-\epsilon_{3467}\delta_{1234}\gamma_1 tu).
\end{array} \eqno{(067)}
$$

$$
\begin{array}{rcl}
[x_{b+2c+2d}(u),x_{-a-2b-2c-2d}(t)] & = & x_{-a-b}(N_{b+2c+2d,-a-2b-2c-2d} tu) \\[4mm]
                                    & = & x_{-a-b}(-\epsilon_{3489}\delta_{1234}\gamma_2 tu).
\end{array} \eqno{(068)}
$$

$$
\begin{array}{rcl}
[x_{b+2c+2d}(u),x_{-a-2b-3c-2d}(t)] & = & x_{-a-b-c}(N_{b+2c+2d,-a-2b-3c-2d} tu) \\[4mm]
                                    & = & x_{-a-b-c}(-\epsilon_{3480}\delta_{1245}\gamma_2 tu).
\end{array} \eqno{(069)}
$$

$$
\begin{array}{rcl}
[x_{b+2c+2d}(u),x_{-a-2b-4c-2d}(t)] & = & x_{-a-b-2c}(N_{b+2c+2d,-a-2b-4c-2d} tu) \\[4mm]
                                    & = & x_{-a-b-2c}(-\epsilon_{3480}\delta_{1256}\gamma_2 tu).
\end{array} \eqno{(070)}
$$

$$
\begin{array}{rcl}
[x_{b+2c+2d}(u),x_{-a-3b-4c-2d}(t)] & = & x_{-a-2b-2c}(N_{b+2c+2d,-a-3b-4c-2d} tu) \\[4mm]
                                    & = & x_{-a-2b-2c}(-\epsilon_{3480}\delta_{1256}\gamma_3 tu).
\end{array} \eqno{(071)}
$$

{\boldmath$a:$}

$$
[x_a(u),x_{-a-b}(t)] = x_{-b}(N_{a,-a-b} tu) = x_{-b}(\gamma_1 tu).\eqno{(072)}
$$

$$
[x_a(u),x_{-a-b-c}(t)] = x_{-b-c}(N_{a,-a-b-c} tu) = x_{-b-c}(\delta_{13}\gamma_1 tu).\eqno{(073)}
$$

$$
[x_a(u),x_{-a-b-c-d}(t)] = x_{-b-c-d}(N_{a,-a-b-c-d} tu) = x_{-b-c-d}(\epsilon_{25}\delta_{13}\gamma_1 tu).\eqno{(074)}
$$

$$
[x_a(u),x_{-a-b-2c}(t)] = x_{-b-2c}(N_{a,-a-b-2c} tu) = x_{-b-2c}(\delta_{1234}\gamma_1 tu).\eqno{(075)}
$$

$$
[x_a(u),x_{-a-b-2c-d}(t)] = x_{-b-2c-d}(N_{a,-a-b-2c-d} tu) = x_{-b-2c-d}(\epsilon_{36}\delta_{1234}\gamma_1 tu).\eqno{(076)}
$$

$$
\begin{array}{rcl}
[x_a(u),x_{-a-b-2c-2d}(t)] & = & x_{-b-2c-2d}(N_{a,-a-b-2c-2d} tu) \\[4mm]
                           & = & x_{-b-2c-2d}(\epsilon_{3467}\delta_{1234}\gamma_1 tu).
\end{array} \eqno{(077)}
$$

$$
\begin{array}{rcl}
[x_a(u),x_{-2a-3b-4c-2d}(t)] & = & x_{-a-3b-4c-2d}(N_{a,-2a-3b-4c-2d} tu) \\[4mm]
                             & = & x_{-a-3b-4c-2d}(-\alpha_1 tu).
\end{array} \eqno{(078)}
$$

{\boldmath$a+b:$}

$$
[x_{a+b}(u),x_{-a-b-c}(t)] = x_{-c}(N_{a+b,-a-b-c} tu) = x_{-c}(\delta_3 tu).\eqno{(079)}
$$

$$
[x_{a+b}(u),x_{-a-b-c-d}(t)] = x_{-c-d}(N_{a+b,-a-b-c-d} tu) = x_{-c-d}(\epsilon_{15}\delta_3 tu).\eqno{(080)}
$$

$$
[x_{a+b}(u),x_{-a-2b-2c}(t)] = x_{-b-2c}(N_{a+b,-a-2b-2c} tu) = x_{-b-2c}(\delta_{1234}\gamma_2 tu).\eqno{(081)}
$$

$$
\begin{array}{rcl}
[x_{a+b}(u),x_{-a-2b-2c-d}(t)] & = & x_{-b-2c-d}(N_{a+b,-a-2b-2c-d} tu) \\[4mm]
                               & = & x_{-b-2c-d}(\epsilon_{38}\delta_{1234}\gamma_2 tu).
\end{array}%
 \eqno{(136)}
$$

\newpage

\section{Special Case}

\subsection{The Table of the Structure Constant}

\textbf{Theorem 2.} \textit{Let the values of all structure constants corresponding to extraspecial pairs be positive.
Then the values of arbitrary structure constants are equal to the values indicated in the following tables 5 and 6.}
\vskip5mm

\rotatebox{90}{
%
}

\newpage

\subsection{List of Formulas}

\vskip5mm
\begin{center}
\textit{Positive roots}
\end{center}

{\boldmath$d:$}

$$
[x_d(u),x_c(t)] = x_{c+d}(tu). \eqno{(01)}
$$

$$
[x_d(u),x_{b+c}(t)] = x_{b+c+d}(tu). \eqno{(02)}
$$

$$
[x_d(u),x_{b+2c}(t)] = x_{b+2c+d}(tu)x_{b+2c+2d}(-tu^2). \eqno{(03)}
$$

$$
[x_d(u),x_{b+2c+d}(t)] = x_{b+2c+2d}(2tu).\eqno{(04)}
$$

$$
[x_d(u),x_{a+b+c}(t)] = x_{a+b+c+d}(tu).\eqno{(05)}
$$

$$
[x_d(u),x_{a+b+2c}(t)] = x_{a+b+2c+d}(tu)\ x_{a+b+2c+2d}(-tu^2). \eqno{(06)}
$$

$$
[x_d(u),x_{a+b+2c+d}(t)] = x_{a+b+2c+2d}(2tu). \eqno{(07)}
$$

$$
[x_d(u),x_{a+2b+2c}(t)] = x_{a+2b+2c+d}(tu)\ x_{a+2b+2c+2d}(-tu^2). \eqno{(08)}
$$

$$
[x_d(u),x_{a+2b+2c+d}(t)] = x_{a+2b+2c+2d}(2tu).\eqno{(09)}
$$

$$
[x_d(u),x_{a+2b+3c+d}(t)] = x_{a+2b+3c+2d}(tu).\eqno{(10)}
$$

{\boldmath$c:$}

$$
[x_c(u),x_{b}(t)] = x_{b+c}(tu)\ x_{b+2c}(-tu^2). \eqno{(11)}
$$

$$
[x_c(u),x_{b+c}(t)] = x_{b+2c}(2tu). \eqno{(12)}
$$

$$
[x_c(u),x_{b+c+d}(t)] = x_{b+2c+d}(tu).\eqno{(13)}
$$

$$
[x_c(u),x_{a+b}(t)] = x_{a+b+c}(tu)\ x_{a+b+2c}(-tu^2). \eqno{(14)}
$$

$$
[x_c(u),x_{a+b+c}(t)] = x_{a+b+2c}(2tu). \eqno{(15)}
$$

$$
[x_c(u),x_{a+b+c+d}(t)] = x_{a+b+2c+d}(tu).\eqno{(16)}
$$

$$
[x_c(u),x_{a+2b+2c+d}(t)] = x_{a+2b+3c+d}(tu). \eqno{(17)}
$$

$$
[x_c(u),x_{a+2b+2c+2d}(t)] = x_{a+2b+3c+2d}(tu)\ x_{a+2b+4c+2d}(-tu^2). \eqno{(18)}
$$

$$
[x_c(u),x_{a+2b+3c+2d}(t)] = x_{a+2b+4c+2d}(2tu).\eqno{(19)}
$$

{\boldmath$c+d:$}

$$
[x_{c+d}(u),x_{b}(t)] = x_{b+c+d}(tu)\ x_{b+2c+2d}(-tu^2). \eqno{(20)}
$$

$$
[x_{c+d}(u),x_{b+c+d}(t)] = x_{b+2c+2d}(2tu). \eqno{(21)}
$$

$$
[x_{c+d}(u),x_{b+c}(t)] = x_{b+2c+d}(tu). \eqno{(22)}
$$

$$
[x_{c+d}(u),x_{a+b}(t)] = x_{a+b+c+d}(tu)\ x_{a+b+2c+2d}(-tu^2). \eqno{(23)}
$$

$$
[x_{c+d}(u),x_{a+b+c+d}(t)] = x_{a+b+2c+2d}(2tu). \eqno{(24)}
$$

$$
[x_{c+d}(u),x_{a+b+c}(t)] = x_{a+b+2c+d}(tu).\eqno{(25)}
$$

$$
[x_{c+d}(u),x_{a+2b+2c}(t)] = x_{a+2b+3c+d}(-tu)\ x_{a+2b+4c+2d}(-tu^2). \eqno{(26)}
$$

$$
[x_{c+d}(u),x_{a+2b+3c+d}(t)] = x_{a+2b+4c+2d}(-2tu). \eqno{(27)}
$$

$$
[x_{c+d}(u),x_{a+2b+2c+d}(t)] = x_{a+2b+3c+2d}(-tu). \eqno{(28)}
$$

{\boldmath$b:$}

$$
[x_b(u),x_a(t)] = x_{a+b}(tu).\eqno{(29)}
$$

$$
[x_b(u),x_{a+b+2c}(t)] = x_{a+2b+2c}(tu).\eqno{(30)}
$$

$$
[x_b(u),x_{a+b+2c+d}(t)] = x_{a+2b+2c+d}(tu).\eqno{(31)}
$$

$$
[x_b(u),x_{a+b+2c+2d}(t)] = x_{a+2b+2c+2d}(tu). \eqno{(32)}
$$

$$
[x_b(u),x_{a+2b+4c+2d}(t)] = x_{a+3b+4c+2d}(tu).\eqno{(33)}
$$

{\boldmath$b+c:$}

$$
[x_{b+c}(u),x_a(t)] = x_{a+b+c}(tu)\ x_{a+2b+2c}(tu^2). \eqno{(34)}
$$

$$
[x_{b+c}(u),x_{a+b+c}(t)] = x_{a+2b+2c}(-2tu). \eqno{(35)}
$$

$$
[x_{b+c}(u),x_{a+b+c+d}(t)] = x_{a+2b+2c+d}(-tu).\eqno{(36)}
$$

$$
[x_{b+c}(u),x_{a+b+2c+d}(t)] = x_{a+2b+3c+d}(tu). \eqno{(37)}
$$

$$
[x_{b+c}(u),x_{a+b+2c+2d}(t)] = x_{a+2b+3c+2d}(tu)\ x_{a+3b+4c+2d}(tu^2). \eqno{(38)}
$$

$$
[x_{b+c}(u),x_{a+2b+3c+2d}(t)] = x_{a+b}(-2tu).\eqno{(39)}
$$

{\boldmath$b+c+d:$}

$$
[x_{b+c+d}(u),x_{a}(t)] = x_{a+b+c+d}(tu)\ x_{a+2b+2c+2d}(tu^2). \eqno{(40)}
$$

$$
[x_{b+c+d}(u),x_{a+b+c+d}(t)] = x_{a+2b+2c+2d}(-2tu).\eqno{(41)}
$$

$$
[x_{b+c+d}(u),x_{a+b+c}(t)] = x_{a+2b+2c+d}(-tu).\eqno{(42)}
$$

$$
[x_{b+c+d}(u),x_{a+b+2c}(t)] = x_{a+2b+3c+d}(-tu)\ x_{a+3b+4c+2d}(tu^2). \eqno{(43)}
$$

$$
[x_{b+c+d}(u),x_{a+2b+3c+d}(t)] = x_{a+3b+4c+2d}(2tu). \eqno{(44)}
$$

$$
[x_{b+c+d}(u),x_{a+b+2c+d}(t)] = x_{a+2b+3c+2d}(-tu).\eqno{(45)}
$$

{\boldmath$b+2c:$}

$$
[x_{b+2c}(u),x_{a}(t)] = x_{a+b+2c}(tu). \eqno{(46)}
$$

$$
[x_{b+2c}(u),x_{a+b}(t)] = x_{a+2b+2c}(tu). \eqno{(47)}
$$

$$
[x_{b+2c}(u),x_{a+b+c+d}(t)] = x_{a+2b+3c+d}(-tu).\eqno{(48)}
$$

$$
[x_{b+2c}(u),x_{a+b+2c+2d}(t)] = x_{a+2b+4c+2d}(tu). \eqno{(49)}
$$

$$
[x_{b+2c}(u),x_{a+2b+2c+2d}(t)] =x_{a+3b+4c+2d}(tu).\eqno{(50)}
$$

{\boldmath$b+2c+d:$}

$$
[x_{b+2c+d}(u),x_{a}(t)] = x_{a+b+2c+d}(tu)\ x_{a+2b+4c+2d}(tu^2). \eqno{(51)}
$$

$$
[x_{b+2c+d}(u),x_{a+b+2c+d}(t)] = x_{a+2b+4c+2d}(-2tu).\eqno{(52)}
$$

$$
[x_{b+2c+d}(u),x_{a+b}(t)] = x_{a+2b+2c+d}(tu)\ x_{a+3b+4c+2d}(tu^2). \eqno{(53)}
$$

$$
[x_{b+2c+d}(u),x_{a+2b+2c+d}(t)] = x_{a+3b+4c+2d}(-2tu). \eqno{(54)}
$$

$$
[x_{b+2c+d}(u),x_{a+b+c}(t)] = x_{a+2b+3c+d}(tu).\eqno{(55)}
$$

$$
[x_{b+2c+d}(u),x_{a+b+c+d}(t)] = x_{a+2b+3c+2d}(-tu). \eqno{(56)}
$$

{\boldmath$b+2c+2d:$}

$$
[x_{b+2c+2d}(u),x_a(t)] = x_{a+b+2c+2d}(tu).\eqno{(57)}
$$

$$
[x_{b+2c+2d}(u),x_{a+b}(t)] = x_{a+2b+2c+2d}(tu).\eqno{(58)}
$$

$$
[x_{b+2c+2d}(u),x_{a+b+c}(t)] = x_{a+2b+3c+2d}(tu). \eqno{(59)}
$$

$$
[x_{b+2c+2d}(u),x_{a+b+2c}(t)] = x_{a+2b+4c+2d}(tu). \eqno{(60)}
$$

$$
[x_{b+2c+2d}(u),x_{a+2b+2c}(t)] = x_{a+3b+4c+2d}(tu).\eqno{(61)}
$$

{\boldmath$a:$}

$$
[x_a(u),x_{a+3b+4c+2d}(t)] = x_{2a+3b+4c+2d}(tu). \eqno{(62)}
$$

{\boldmath$a+b:$}

$$
[x_{a+b}(u),x_{a+2b+4c+2d}(t)] = x_{2a+3b+4c+2d}(-tu). \eqno{(63)}
$$

{\boldmath$a+b+c:$}

$$
[x_{a+b+c}(u),x_{a+2b+3c+2d}(t)] = x_{2a+3b+4c+2d}(2tu). \eqno{(64)}
$$

{\boldmath$a+b+c+d:$}

$$
[x_{a+b+c+d}(u),x_{a+2b+3c+d}(t)] = x_{2a+3b+4c+2d}(-2tu). \eqno{(65)}
$$

{\boldmath$a+b+2c:$}

$$
[x_{a+b+2c}(u),x_{a+2b+2c+2d}(t)] = x_{2a+3b+4c+2d}(-tu). \eqno{(66)}
$$

{\boldmath$a+b+2c+d:$}

$$
[x_{a+b+2c+d}(u),x_{a+2b+2c+d}(t)] = x_{2a+3b+4c+2d}(2tu). \eqno{(67)}
$$

{\boldmath$a+b+2c+2d:$}

$$
[x_{a+b+2c+2d}(u),x_{a+2b+2c}(t)] = x_{2a+3b+4c+2d}(-tu). \eqno{(68)}
$$

\vskip5mm
\begin{center}
\textit{Negative roots}
\end{center}

{\boldmath$-d:$}

$$
[x_{-d}(u),x_{-c}(t)] = x_{-c-d}(-tu). \eqno{(01)}
$$

$$
[x_{-d}(u),x_{-b-c}(t)] = x_{-b-c-d}(-tu). \eqno{(02)}
$$

$$
[x_{-d}(u),x_{-b-2c}(t)] = x_{-b-2c-d}(-tu)x_{-b-2c-d}(-tu). \eqno{(03)}
$$

$$
[x_{-d}(u),x_{-b-2c-d}(t)] = x_{-b-2c-2d}(-2tu).\eqno{(04)}
$$

$$
[x_{-d}(u),x_{-a-b-c}(t)] = x_{-a-b-c-d}(-tu). \eqno{(05)}
$$

$$
[x_{-d}(u),x_{-a-b-2c}(t)] = x_{-a-b-2c-d}(-tu)\ x_{-a-b-2c-2d}(-tu^2). \eqno{(06)}
$$

$$
[x_{-d}(u),x_{-a-b-2c-d}(t)] = x_{-a-b-2c-2d}(-2 tu). \eqno{(07)}
$$

$$
[x_{-d}(u),x_{-a-2b-2c}(t)] = x_{-a-2b-2c-d}(-tu)\ x_{-a-2b-2c-2d}(-tu^2). \eqno{(08)}
$$

$$
[x_{-d}(u),x_{-a-2b-2c-d}(t)] = x_{-a-2b-2c-2d}(-2tu). \eqno{(09)}
$$

$$
[x_{-d}(u),x_{-a-2b-3c-d}(t)] = x_{-a-2b-3c-2d}(- tu). \eqno{(10)}
$$

{\boldmath$-c:$}

$$
[x_{-c}(u),x_{-b}(t)] = x_{-b-c}(-tu)\ x_{-b-2c}(-tu^2). \eqno{(11)}
$$

$$
[x_{-c}(u),x_{-b-c}(t)] = x_{-b-2c}(-2tu). \eqno{(12)}
$$

$$
[x_{-c}(u),x_{-b-c-d}(t)] = x_{-b-2c-d}(-tu). \eqno{(13)}
$$

$$
[x_{-c}(u),x_{-a-b}(t)] = x_{-a-b-c}(-tu)\ x_{-a-b-2c}(-tu^2). \eqno{(14)}
$$

$$
[x_{-c}(u),x_{-a-b-c}(t)] = x_{-a-b-2c}(-2tu). \eqno{(15)}
$$

$$
[x_{-c}(u),x_{-a-b-c-d}(t)] = x_{-a-b-2c-d}(-tu). \eqno{(16)}
$$

$$
[x_{-c}(u),x_{-a-2b-2c-d}(t)] = x_{-a-2b-3c-d}(-tu). \eqno{(17)}
$$

$$
[x_{-c}(u),x_{-a-2b-2c-2d}(t)] = x_{-a-2b-3c-2d}(-tu)\ x_{-a-2b-4c-2d}(-tu^2). \eqno{(18)}
$$

$$
[x_{-c}(u),x_{-a-2b-3c-2d}(t)] = x_{-a-2b-4c-2d}(-2tu). \eqno{(19)}
$$

{\boldmath$-c-d:$}

$$
[x_{-c-d}(u),x_{-b}(t)] = x_{-b-c-d}(-tu)\ x_{b-2c-2d}(-tu^2). \eqno{(20)}
$$

$$
[x_{-c-d}(u),x_{-b-c-d}(t)] = x_{-b-2c-2d}(-2tu). \eqno{(21)}
$$

$$
[x_{-c-d}(u),x_{-b-c}(t)] = x_{-b-2c-d}(-tu). \eqno{(22)}
$$

$$
[x_{-c-d}(u),x_{-a-b}(t)] = x_{-a-b-c-d}(-tu)\ x_{-a-b-2c-2d}(-tu^2). \eqno{(23)}
$$

$$
[x_{-c-d}(u),x_{-a-b-c-d}(t)] = x_{-a-b-2c-2d}(-2tu). \eqno{(24)}
$$

$$
[x_{-c-d}(u),x_{-a-b-c}(t)] = x_{-a-b-2c-d}(-tu). \eqno{(25)}
$$

$$
[x_{-c-d}(u),x_{-a-2b-2c}(t)] = x_{-a-2b-3c-d}(tu)\ x_{-a-2b-4c-2d}(-tu^2). \eqno{(26)}
$$

$$
[x_{-c-d}(u),x_{-a-2b-3c-d}(t)] = x_{-a-2b-4c-2d}(2tu). \eqno{(27)}
$$

$$
[x_{-c-d}(u),x_{-a-2b-2c-d}(t)] = x_{-a-2b-3c-2d}(tu). \eqno{(28)}
$$

{\boldmath$-b:$}

$$
[x_{-b}(u),x_{-a}(t)] = x_{-a-b}(-tu). \eqno{(29)}
$$

$$
[x_{-b}(u),x_{-a-b-2c}(t)] = x_{-a-2b-2c}(-tu). \eqno{(30)}
$$

$$
[x_{-b}(u),x_{-a-b-2c-d}(t)] = x_{-a-2b-2c-d}(-tu). \eqno{(31)}
$$

$$
[x_{-b}(u),x_{-a-b-2c-2d}(t)] = x_{-a-2b-2c-2d}(-tu). \eqno{(32)}
$$

$$
[x_{-b}(u),x_{-a-2b-4c-2d}(t)] = x_{-a-3b-4c-2d}(-tu).\eqno{(33)}
$$

{\boldmath$-b-c:$}

$$
[x_{-b-c}(u),x_{-a}(t)] = x_{-a-b-c}(-tu)\ x_{-a-2b-2c}(tu^2). \eqno{(34)}
$$

$$
[x_{-b-c}(u),x_{-a-b-c}(t)] = x_{-a-2b-2c}(2tu). \eqno{(35)}
$$

$$
[x_{-b-c}(u),x_{-a-b-c-d}(t)] = x_{-a-2b-2c-d}(tu). \eqno{(36)}
$$

$$
[x_{-b-c}(u),x_{-a-b-2c-d}(t)] = x_{-a-2b-3c-d}(-tu). \eqno{(37)}
$$

$$
[x_{-b-c}(u),x_{-a-b-2c-2d}(t)] = x_{-a-2b-3c-2d}(-tu)\ x_{-a-3b-4c-2d}(tu^2). \eqno{(38)}
$$

$$
[x_{-b-c}(u),x_{-a-2b-3c-2d}(t)] = x_{-a-3b-4c-2d}(2tu). \eqno{(39)}
$$

{\boldmath$-b-c-d:$}

$$
[x_{-b-c-d}(u),x_{-a}(t)] = x_{-a-b-c-d}(-tu)\ x_{-a-2b-2c-2d}(tu^2). \eqno{(40)}
$$

$$
[x_{-b-c-d}(u),x_{-a-b-c-d}(t)] = x_{-a-2b-2c-2d}(2tu). \eqno{(41)}
$$

$$
[x_{-b-c-d}(u),x_{-a-b-c}(t)] = x_{-a-2b-2c-2d}(tu). \eqno{(42)}
$$

$$
[x_{-b-c-d}(u),x_{-a-b-2c}(t)] = x_{-a-2b-3c-d}(tu)\ x_{-a-3b-4c-2d}(tu^2). \eqno{(43)}
$$

$$
[x_{-b-c-d}(u),x_{-a-2b-3c-d}(t)] = x_{-a-3b-4c-2d}(-2tu). \eqno{(44)}
$$

$$
[x_{-b-c-d}(u),x_{-a-b-2c-d}(t)] = x_{-a-2b-3c-2d}(tu). \eqno{(45)}
$$

{\boldmath$-b-2c:$}

$$
[x_{-b-2c}(u),x_{-a}(t)] = x_{-a-b-2c}(-tu).\eqno{(46)}
$$

$$
[x_{-b-2c}(u),x_{-a-b}(t)] = x_{-a-2b-2c}(-tu). \eqno{(47)}
$$

$$
[x_{-b-2c}(u),x_{-a-b-c-d}(t)] = x_{-a-2b-3c-d}(tu). \eqno{(48)}
$$

$$
[x_{-b-2c}(u),x_{-a-b-2c-2d}(t)] = x_{-a-2b-4c-2d}(-tu). \eqno{(49)}
$$

$$
[x_{-b-2c}(u),x_{-a-2b-2c-2d}(t)] = x_{-a-3b-4c-2d}(-tu).\eqno{(50)}
$$

{\boldmath$-b-2c-d:$}

$$
[x_{-b-2c-d}(u),x_{-a}(t)] = x_{-a-b-2c-d}(-tu)\ x_{-a-2b-4c-2d}(tu^2).\eqno{(51)}
$$

$$
[x_{-b-2c-d}(u),x_{-a-b-2c-d}(t)] = x_{-a-2b-4c-2d}(2tu). \eqno{(52)}
$$

$$
[x_{-b-2c-d}(u),x_{-a-b}(t)] = x_{-a-2b-2c-d}(-tu)\ x_{-a-3b-4c-2d}(tu^2). \eqno{(53)}
$$

$$
[x_{-b-2c-d}(u),x_{-a-2b-2c-d}(t)] = x_{-a-3b-4c-2d}(2tu). \eqno{(54)}
$$

$$
[x_{-b-2c-d}(u),x_{-a-b-c}(t)] = x_{-a-2b-3c-d}(-tu). \eqno{(55)}
$$

$$
[x_{-b-2c-d}(u),x_{-a-b-c-d}(t)] = x_{-a-2b-3c-2d}(tu). \eqno{(56)}
$$

{\boldmath$-b-2c-2d:$}

$$
[x_{-b-2c-2d}(u),x_{-a}(t)] = x_{-a-b-2c-2d}(-tu). \eqno{(57)}
$$

$$
[x_{-b-2c-2d}(u),x_{-a-b}(t)] = x_{-a-2b-2c-2d}(-tu).\eqno{(58)}
$$

$$
[x_{-b-2c-2d}(u),x_{-a-b-c}(t)] = x_{-a-2b-3c-2d}(-tu). \eqno{(59)}
$$

$$
[x_{-b-2c-2d}(u),x_{-a-b-2c}(t)] = x_{-a-2b-4c-2d}(-tu). \eqno{(60)}
$$

$$
[x_{-b-2c-2d}(u),x_{-a-2b-2c}(t)] = x_{-a-3b-4c-2d}(-tu). \eqno{(61)}
$$

{\boldmath$-a:$}

$$
[x_{-a}(u),x_{-a-3b-4c-2d}(t)] = x_{-2a-3b-4c-2d}(-tu). \eqno{(62)}
$$

{\boldmath$-a-b:$}

$$
[x_{-a-b}(u),x_{-a-2b-4c-2d}(t)] = x_{-2a-3b-4c-2d}(tu). \eqno{(63)}
$$

{\boldmath$-a-b-c:$}

$$
[x_{-a-b-c}(u),x_{-a-2b-3c-2d}(t)] = x_{-2a-3b-4c-2d}(-2tu). \eqno{(64)}
$$

{\boldmath$-a-b-c-d:$}

$$
[x_{-a-b-c-d}(u),x_{-a-2b-3c-d}(t)] = x_{-2a-3b-4c-2d}(2tu). \eqno{(65)}
$$

{\boldmath$-a-b-2c:$}

$$
[x_{-a-b-2c}(u),x_{-a-2b-2c-2d}(t)] = x_{-2a-3b-4c-2d}(tu). \eqno{(66)}
$$

{\boldmath$-a-b-2c-d:$}

$$
[x_{-a-b-2c-d}(u),x_{-a-2b-2c-d}(t)] = x_{-2a-3b-4c-2d}(-2tu). \eqno{(67)}
$$

{\boldmath$-a-b-2c-2d:$}

$$
[x_{-a-b-2c-2d}(u),x_{-a-2b-2c}(t)] = x_{-2a-3b-4c-2d}(tu). \eqno{(68)}
$$

\vskip5mm
\begin{center}
\textit{Positive and negative roots}
\end{center}

{\boldmath$d:$}

$$
[x_d(u),x_{-c-d}(t)] =  x_{-c}(-tu). \eqno{(001)}
$$

$$
[x_d(u),x_{-b-c-d}(t)] = x_{-b-c}(-tu). \eqno{(002)}
$$

$$
[x_d(u),x_{-b-2c-d}(t)] = x_{-b-2c}(-2tu). \eqno{(003)}
$$

$$
[x_d(u),x_{-b-2c-2d}(t)] = x_{-b-2c-d}(-tu). x_{-b-2c}(-tu^2). \eqno{(004)}
$$

$$
[x_d(u),x_{-a-b-c-d}(t)] = x_{-a-b-c}(-tu). \eqno{(005)}
$$

$$
[x_d(u),x_{-a-b-2c-d}(t)] = x_{-a-b-2c}(-2tu). \eqno{(006)}
$$

$$
[x_d(u),x_{-a-b-2c-2d}(t)] = x_{-a-b-2c-d}(-tu)\ x_{-a-b-2c}(-tu^2). \eqno{(007)}
$$

$$
[x_d(u),x_{-a-2b-2c-d}(t)] = x_{-a-2b-2c}(-2tu). \eqno{(008)}
$$

$$
[x_d(u),x_{-a-2b-2c-2d}(t)] = x_{-a-2b-2c-d}(-tu)\ x_{-a-2b-2c}(-tu^2). \eqno{(009)}
$$

$$
[x_d(u),x_{-a-2b-3c-2d}(t)] = x_{-a-2b-3c-d}(-tu). \eqno{(010)}
$$

{\boldmath$c:$}

$$
[x_c(u),x_{-c-d}(t)] = x_{-d}(tu).\eqno{(011)}
$$

$$
[x_c(u),x_{-b-c}(t)] = x_{-b}(-2tu).\eqno{(012)}
$$

$$
[x_c(u),x_{-b-2c}(t)] = x_{-b-c}(-tu)\ x_{-b}(-tu^2). \eqno{(013)}
$$

$$
[x_c(u),x_{-b-2c-d}(t)] = x_{-b-c-d}(-tu). \eqno{(014)}
$$

$$
[x_c(u),x_{-a-b-c}(t)] = x_{-a-b}(-2tu). \eqno{(015)}
$$

$$
[x_c(u),x_{-a-b-2c}(t)] = x_{-a-b-c}(-tu)\ x_{-a-b}(-tu^2). \eqno{(016)}
$$

$$
[x_c(u),x_{-a-b-2c-d}(t)] = x_{-a-b-c-d}(-tu). \eqno{(017)}
$$

$$
[x_c(u),x_{-a-2b-3c-d}(t)] = x_{-a-2b-2c-d}(-tu).\eqno{(018)}
$$

$$
[x_c(u),x_{-a-2b-3c-2d}(t)] = x_{-a-2b-2c-2d}(-2tu). \eqno{(019)}
$$

$$
[x_c(u),x_{-a-2b-4c-2d}(t)] = x_{-a-2b-3c-2d}(-tu)\ x_{-a-2b-2c-2d}(-tu^2). \eqno{(020)}
$$

{\boldmath$c+d:$}

$$
[x_{c+d}(u),x_{-b-c-d}(t)] = x_{-b}(-2tu). \eqno{(021)}
$$

$$
[x_{c+d}(u),x_{-b-2c-2d}(t)] = x_{-b-c-d}(-tu)\ x_{-b}(-tu^2). \eqno{(022)}
$$

$$
[x_{c+d}(u),x_{-b-2c-d}(t)] = x_{-b-c}(-tu). \eqno{(023)}
$$

$$
[x_{c+d}(u),x_{-a-b-c-d}(t)] = x_{-a-b}(-2tu). \eqno{(024)}
$$

$$
[x_{c+d}(u),x_{-a-b-2c-2d}(t)] = x_{-a-b-c-d}(-tu)\ x_{-a-b}(-tu^2). \eqno{(025)}
$$

$$
[x_{c+d}(u),x_{-a-b-2c-d}(t)] = x_{-a-b-c}(-tu). \eqno{(026)}
$$

$$
[x_{c+d}(u),x_{-a-2b-3c-d}(t)] = x_{-a-2b-2c}(2tu). \eqno{(027)}
$$

$$
[x_{c+d}(u),x_{-a-2b-4c-2d}(t)] = x_{-a-2b-3c-d}(tu)\ x_{-a-2b-2c}(-tu^2). \eqno{(028)}
$$

$$
[x_{c+d}(u),x_{-a-2b-3c-2d}(t)] = x_{-a-2b-2c-d}(tu). \eqno{(029)}
$$

{\boldmath$b:$}

$$
[x_b(u),x_{-b-c}(t)] = x_{-c}(tu). \eqno{(030)}
$$

$$
[x_b(u),x_{-a-b}(t)] = x_{-a}(-tu). \eqno{(031)}
$$

$$
[x_b(u),x_{-b-c-d}(t)] = x_{-c-d}(tu). \eqno{(032)}
$$

$$
[x_b(u),x_{-a-2b-2c}(t)] = x_{-a-b-2c}(-tu). \eqno{(033)}
$$

$$
[x_b(u),x_{-a-2b-2c-d}(t)] = x_{-a-b-2c-d}(-tu). \eqno{(034)}
$$

$$
[x_b(u),x_{-a-2b-2c-2d}(t)] = x_{-a-b-2c-2d}(-tu). \eqno{(035)}
$$

$$
[x_b(u),x_{-a-3b-4c-2d}(t)] = x_{-a-2b-4c-2d}(-tu). \eqno{(036)}
$$

{\boldmath$b+c:$}

$$
[x_{b+c}(u),x_{-b-c-d}(t)] = x_{-d}(tu). \eqno{(037)}
$$

$$
[x_{b+c}(u),x_{-b-2c}(t)] = x_{-c}(tu). \eqno{(038)}
$$

$$
[x_{b+c}(u),x_{-a-b-c}(t)] = x_{-a}(-2tu). \eqno{(039)}
$$

$$
[x_{b+c}(u),x_{-a-2b-2c}(t)] = x_{-a-b-c}(tu)\ x_{-a}(tu^2). \eqno{(040)}
$$

$$
[x_{b+c}(u),x_{-b-2c-d}(t)] = x_{-c-d}(tu). \eqno{(041)}
$$

$$
[x_{b+c}(u),x_{-a-2b-2c-d}(t)] = x_{-a-b-c-d}(tu). \eqno{(042)}
$$

$$
[x_{b+c}(u),x_{-a-2b-3c-d}(t)] = x_{-a-b-2c-d}(-tu). \eqno{(043)}
$$

$$
[x_{b+c}(u),x_{-a-2b-3c-2d}(t)] = x_{-a-b-2c-2d}(-2tu). \eqno{(044)}
$$

$$
[x_{b+c}(u),x_{-a-3b-4c-2d}(t)] = x_{-a-2b-3c-2d}(tu)\ x_{-a-b-2c-2d}(tu^2). \eqno{(045)}
$$

{\boldmath$b+c+d:$}

$$
[x_{b+c+d}(u),x_{-a-b-c-d}(t)] = x_{-d}(-2tu). \eqno{(046)}
$$

$$
[x_{b+c+d}(u),x_{-a-2b-2c-2d}(t)] = x_{-a-b-c-d}(tu)\ x_{-a}(tu^2). \eqno{(047)}
$$

$$
[x_{b+c+d}(u),x_{-b-2c-d}(t)] = x_{-c}(tu). \eqno{(048)}
$$

$$
[x_{b+c+d}(u),x_{-b-2c-2d}(t)] = x_{-c-d}(tu). \eqno{(049)}
$$

$$
[x_{b+c+d}(u),x_{-a-2b-2c-d}(t)] = x_{-a-b-c}(tu). \eqno{(050)}
$$

$$
[x_{b+c+d}(u),x_{-a-2b-3c-d}(t)] = x_{-a-b-2c}(2tu). \eqno{(051)}
$$

$$
[x_{b+c+d}(u),x_{-a-3b-4c-2d}(t)] = x_{-a-2b-3c-d}(-tu)\ x_{-a-b-2c}(tu^2). \eqno{(052)}
$$

$$
[x_{b+c+d}(u),x_{-a-2b-3c-2d}(t)] = x_{-a-b-2c-d}(tu). \eqno{(053)}
$$

{\boldmath$b+2c:$}

$$
[x_{b+2c}(u),x_{-b-2c-d}(t)] = x_{-d}(tu). \eqno{(054)}
$$

$$
[x_{b+2c}(u),x_{-a-b-2c}(t)] = x_{-a}(-tu). \eqno{(055)}
$$

$$
[x_{b+2c}(u),x_{-a-2b-2c}(t)] = x_{-a-b}(-tu). \eqno{(056)}
$$

$$
[x_{b+2c}(u),x_{-a-2b-3c-d}(t)] = x_{-a-b-c-d}(tu). \eqno{(057)}
$$

$$
[x_{b+2c}(u),x_{-a-2b-4c-2d}(t)] = x_{-a-b-2c-2d}(-tu). \eqno{(058)}
$$

$$
[x_{b+2c}(u),x_{-a-3b-4c-2d}(t)] = x_{-a-2b-2c-2d}(-tu). \eqno{(059)}
$$

{\boldmath$b+2c+d:$}

$$
[x_{b+2c+d}(u),x_{-a-b-2c-d}(t)] = x_{-a}(-2tu). \eqno{(060)}
$$

$$
[x_{b+2c+d}(u),x_{-a-2b-4c-2d}(t)] = x_{-a-b-2c-d}(tu)\ x_{-a}(tu^2). \eqno{(061)}
$$

$$
[x_{b+2c+d}(u),x_{-b-2c-2d}(t)] = x_{-d}(tu). \eqno{(062)}
$$

$$
[x_{b+2c+d}(u),x_{-a-2b-2c-d}(t)] = x_{-a-b}(-2tu). \eqno{(063)}
$$

$$
[x_{b+2c+d}(u),x_{-a-3b-4c-2d}(t)] = x_{-a-2b-2c-d}(tu)\ x_{-a-b}(tu^2). \eqno{(064)}
$$

$$
[x_{b+2c+d}(u),x_{-a-2b-3c-d}(t)] = x_{-a-b-c}(-tu). \eqno{(065)}
$$

$$
[x_{b+2c+d}(u),x_{-a-2b-3c-2d}(t)] = x_{-a-b-c-d}(tu). \eqno{(066)}
$$

{\boldmath$b+2c+2d:$}

$$
[x_{b+2c+2d}(u),x_{-a-b-2c-2d}(t)] = x_{-a}(-tu). \eqno{(067)}
$$

$$
[x_{b+2c+2d}(u),x_{-a-2b-2c-2d}(t)] = x_{-a-b}(-tu). \eqno{(068)}
$$

$$
[x_{b+2c+2d}(u),x_{-a-2b-3c-2d}(t)] = x_{-a-b-c}(-tu). \eqno{(069)}
$$

$$
[x_{b+2c+2d}(u),x_{-a-2b-4c-2d}(t)] = x_{-a-b-2c}(-tu). \eqno{(070)}
$$

$$
[x_{b+2c+2d}(u),x_{-a-3b-4c-2d}(t)] = x_{-a-2b-2c}(-tu). \eqno{(071)}
$$

{\boldmath$a:$}

$$
[x_a(u),x_{-a-b}(t)] = x_{-b}(tu). \eqno{(072)}
$$

$$
[x_a(u),x_{-a-b-c}(t)] = x_{-b-c}(tu). \eqno{(073)}
$$

$$
[x_a(u),x_{-a-b-c-d}(t)] = x_{-b-c-d}(tu). \eqno{(074)}
$$

$$
[x_a(u),x_{-a-b-2c}(t)] = x_{-b-2c}(tu). \eqno{(075)}
$$

$$
[x_a(u),x_{-a-b-2c-d}(t)] = x_{-b-2c-d}(tu). \eqno{(076)}
$$

$$
[x_a(u),x_{-a-b-2c-2d}(t)] = x_{-b-2c-2d}(tu). \eqno{(077)}
$$

$$
[x_a(u),x_{-2a-3b-4c-2d}(t)] = x_{-a-3b-4c-2d}(-tu). \eqno{(078)}
$$

{\boldmath$a+b:$}

$$
[x_{a+b}(u),x_{-a-b-c}(t)] = x_{-c}(tu). \eqno{(079)}
$$

$$
[x_{a+b}(u),x_{-a-b-c-d}(t)] = x_{-c-d}(tu). \eqno{(080)}
$$

$$
[x_{a+b}(u),x_{-a-2b-2c}(t)] = x_{-b-2c}(tu). \eqno{(081)}
$$

$$
[x_{a+b}(u),x_{-a-2b-2c-d}(t)] = x_{-b-2c-d}(tu). \eqno{(082)}
$$

$$
[x_{a+b}(u),x_{-a-2b-2c-2d}(t)] = x_{-b-2c-2d}(tu). \eqno{(083)}
$$

$$
[x_{a+b}(u),x_{-2a-3b-4c-2d}(t)] = x_{-a-2b-4c-2d}(tu). \eqno{(084)}
$$

{\boldmath$a+b+c:$}

$$
[x_{a+b+c}(u),x_{-a-b-c-d}(t)] = x_{-d}(tu). \eqno{(085)}
$$

$$
[x_{a+b+c}(u),x_{-a-b-2c}(t)] = x_{-c}(tu). \eqno{(086)}
$$

$$
[x_{a+b+c}(u),x_{-a-b-2c-d}(t)] = x_{-c-d}(tu). \eqno{(087)}
$$

$$
[x_{a+b+c}(u),x_{-a-2b-2c}(t)] = x_{-b-c}(-tu). \eqno{(088)}
$$

$$
[x_{a+b+c}(u),x_{-a-2b-2c-d}(t)] = x_{-b-c-d}(-tu). \eqno{(089)}
$$

$$
[x_{a+b+c}(u),x_{-a-2b-3c-d}(t)] = x_{-b-2c-d}(tu). \eqno{(090)}
$$

$$
[x_{a+b+c}(u),x_{-a-2b-3c-2d}(t)] = x_{-b-2c-2d}(2tu). \eqno{(091)}
$$

$$
[x_{a+b+c}(u),x_{-2a-3b-4c-2d}(t)] = x_{-a-2b-3c-2d}(-tu)\ x_{-b-2c-2d}(tu^2). \eqno{(092)}
$$

{\boldmath$a+b+c+d:$}

$$
[x_{a+b+c+d}(u),x_{-a-b-2c-d}(t)] = x_{-c}(tu). \eqno{(093)}
$$

$$
[x_{a+b+c+d}(u),x_{-a-b-2c-2d}(t)] = x_{-c-d}(tu). \eqno{(094)}
$$

$$
[x_{a+b+c+d}(u),x_{-a-2b-2c-d}(t)] = x_{-b-c}(-tu). \eqno{(095)}
$$

$$
[x_{a+b+c+d}(u),x_{-a-2b-2c-2d}(t)] = x_{-b-c-d}(-tu). \eqno{(096)}
$$

$$
[x_{a+b+c+d}(u),x_{-a-2b-3c-d}(t)] = x_{-b-2c-d}(-2tu). \eqno{(097)}
$$

$$
[x_{a+b+c+d}(u),x_{-2a-3b-4c-2d}(t)] = x_{-a-2b-3c-d}(tu)\ x_{-b-2c}(tu^2). \eqno{(098)}
$$

$$
[x_{a+b+c+d}(u),x_{-a-2b-3c-2d}(t)] = x_{-b-2c-2d}(-tu). \eqno{(099)}
$$

{\boldmath$a+b+2c:$}

$$
[x_{a+b+2c}(u),x_{-a-b-2c-d}(t)] = x_{-d}(-tu). \eqno{(100)}
$$

$$
[x_{a+b+2c}(u),x_{-a-2b-2c}(t)] = x_{-b}(-tu). \eqno{(101)}
$$

$$
[x_{a+b+2c}(u),x_{-a-2b-3c-d}(t)] = x_{-b-c-d}(-tu). \eqno{(102)}
$$

$$
[x_{a+b+2c}(u),x_{-a-2b-4c-2d}(t)] = x_{-b-2c-2d}(tu). \eqno{(103)}
$$

$$
[x_{a+b+2c}(u),x_{-2a-3b-4c-2d}(t)] = x_{-a-2b-2c-2d}(tu). \eqno{(104)}
$$

{\boldmath$a+b+2c+d:$}

$$
[x_{a+b+2c+d}(u),x_{-a-b-2c-2d}(t)] = x_{-d}(tu). \eqno{(105)}
$$

$$
[x_{a+b+2c+d}(u),x_{-a-2b-2c-d}(t)] = x_{-b}(2tu).\eqno{(106)}
$$

$$
[x_{a+b+2c+d}(u),x_{-2a-3b-4c-2d}(t)] = x_{-a-2b-2c-d}(-tu)\ x_{-b}(tu^2). \eqno{(107)}
$$

$$
[x_{a+b+2c+d}(u),x_{-a-2b-3c-d}(t)] = x_{-b-c}(tu). \eqno{(108)}
$$

$$
[x_{a+b+2c+d}(u),x_{-a-2b-3c-2d}(t)] = x_{-b-c-d}(-tu). \eqno{(109)}
$$

$$
[x_{a+b+2c+d}(u),x_{-a-2b-4c-2d}(t)] = x_{-b-2c-d}(-tu). \eqno{(110)}
$$

{\boldmath$a+b+2c+2d:$}

$$
[x_{a+b+2c+2d}(u),x_{-a-2b-2c-2d}(t)] = x_{-b}(tu). \eqno{(111)}
$$

$$
[x_{a+b+2c+2d}(u),x_{-a-2b-3c-2d}(t)] = x_{-b-c}(tu). \eqno{(112)}
$$

$$
[x_{a+b+2c+2d}(u),x_{-a-2b-4c-2d}(t)] = x_{-b-2c}(tu). \eqno{(113)}
$$

$$
[x_{a+b+2c+2d}(u),x_{-2a-3b-4c-2d}(t)] = x_{-a-2b-2c}(tu). \eqno{(114)}
$$

{\boldmath$a+2b+2c:$}

$$
[x_{a+2b+2c}(u),x_{-a-2b-2c-d}(t)] = x_{-d}(tu). \eqno{(115)}
$$

$$
[x_{a+2b+2c}(u),x_{-a-2b-3c-d}(t)] = x_{-c-d}(-tu). \eqno{(116)}
$$

$$
[x_{a+2b+2c}(u),x_{-a-3b-4c-2d}(t)] = x_{-b-2c-2d}(tu). \eqno{(117)}
$$

$$
[x_{a+2b+2c}(u),x_{-2a-3b-4c-2d}(t)] = x_{-a-b-2c-2d}(-tu). \eqno{(118)}
$$

{\boldmath$a+2b+2c+d:$}

$$
[x_{a+2b+2c+d}(u),x_{-a-2b-2c-2d}(t)] = x_{-d}(tu). \eqno{(119)}
$$

$$
[x_{a+2b+2c+d}(u),x_{-a-2b-3c-d}(t)] = x_{-c}(tu). \eqno{(120)}
$$

$$
[x_{a+2b+2c+d}(u),x_{-a-2b-3c-2d}(t)] = x_{-c-d}(-tu). \eqno{(121)}
$$

$$
[x_{a+2b+2c+d}(u),x_{-a-3b-4c-2d}(t)] = x_{-b-2c-d}(-tu). \eqno{(122)}
$$

$$
[x_{a+2b+2c+d}(u),x_{-2a-3b-4c-2d}(t)] = x_{-a-b-2c-d}(tu). \eqno{(123)}
$$

{\boldmath$a+2b+2c+2d:$}

$$
[x_{a+2b+2c+2d}(u),x_{-a-2b-3c-2d}(t)] = x_{-c}(tu). \eqno{(124)}
$$

$$
[x_{a+2b+2c+2d}(u),x_{-a-3b-4c-2d}(t)] = x_{-b-2c}(tu). \eqno{(125)}
$$

$$
[x_{a+2b+2c+2d}(u),x_{-2a-3b-4c-2d}(t)] = x_{-a-b-2c}(-tu). \eqno{(126)}
$$

{\boldmath$a+2b+3c+d:$}

$$
[x_{a+2b+3c+d}(u),x_{-a-2b-3c-2d}(t)] = x_{-d}(tu). \eqno{(127)}
$$

$$
[x_{a+2b+3c+d}(u),x_{-a-2b-4c-2d}(t)] = x_{-c-d}(-tu). \eqno{(128)}
$$

$$
[x_{a+2b+3c+d}(u),x_{-a-3b-4c-2d}(t)] = x_{-b-c-d}(-tu). \eqno{(129)}
$$

$$
[x_{a+2b+3c+d}(u),x_{-2a-3b-4c-2d}(t)] = x_{-a-b-c-d}(-tu). \eqno{(130)}
$$

{\boldmath$a+2b+3c+2d:$}

$$
[x_{a+2b+3c+2d}(u),x_{-a-2b-4c-2d}(t)] = x_{-c}(tu). \eqno{(131)}
$$

$$
[x_{a+2b+3c+2d}(u),x_{-a-3b-4c-2d}(t)] = x_{-b-c}(-tu). \eqno{(132)}
$$

$$
[x_{a+2b+3c+2d}(u),x_{-2a-3b-4c-2d}(t)] = x_{-a-b-c}(tu). \eqno{(133)}
$$

{\boldmath$a+2b+4c+2d:$}

$$
[x_{a+2b+4c+2d}(u),x_{-a-3b-4c-2d}(t)] = x_{-b}(tu). \eqno{(134)}
$$

$$
[x_{a+2b+4c+2d}(u),x_{-2a-3b-4c-2d}(t)] = x_{-a-b}(-tu). \eqno{(135)}
$$

{\boldmath$a+3b+4c+2d:$}

$$
[x_{a+3b+4c+2d}(u),x_{-2a-3b-4c-2d}(t)] = x_{-a}(tu). \eqno{(136)}
$$

\newpage

\section{Numbers $K_{r,s}^\Delta$: Calculation and Applications}

\subsection{Basic Definitions}

In what follows $\Phi$ is a reduced indecomposable root system
in the Euclidean space $V$ with scalar product
$(\,,\,),$ $\Pi(\Phi)$ is its subsystem of fundamental roots;
$N_{r,s},$ $r,s\in\Phi,$ are structure constants of the Lie algebra of type $\Phi,$
$R$ an associative-commutative ring with identity,
$E\Phi(R)$ is an elementary Chevalley group of type $\Phi$ over $R$.

{\bf Definition 1.}
Let $r,s\in\Phi,$ and $\emptyset\ne\Delta\subseteq\Pi(\Phi)$.
If there exist $q_1,\ldots,q_t\in\Delta$ such that for any integer
$i$, $1\leqslant i\leqslant t$, the sum $s+q_1+\cdots+q_i$ is a root,
then the decomposition $r=s+q_1+\cdots+q_t$ will be called an $s$-decomposition of the root $r$
with respect to $\Delta$ (an $s$-decomposition if $\Delta=\Pi(\Phi)$).

From Definition 1 it follows that if there is an
$s$-decomposition of the root $r$ with respect to $\Delta$, then ${\rm ht}(r)-{\rm ht}(s)>0$
and therefore $r,s\in\Phi^+$, or $r,s\in\Phi^-$.
Let $r=s+q_1+\cdots+q_t$ and $r=s+q'_1+\cdots+q'_{t'}$ be two $s$-decompositions of the root $r$ with respect to $\Delta.$ Since representation of $r$ in terms of roots from  $\Pi(\Phi)$  is unique, we claim that  $t=t'$ and
exists a permutation $\pi$ such that $q_i=q'_{\pi(i)}$, $i=1, \ldots, t$.
In this regard, we say that two $s$-decompositions of $r$ with respect to
$\Delta$ are equal if and only if $\pi$ is the  identity permutation.

For fixed $r,s \in \Phi$ and $\Delta \subseteq \Pi(\Phi)$ we define
the set $Q_{r,s}^\Delta$.
We put $Q_{s,s}^\Delta=\{0\}$ and $Q_{r,s}^\Delta=\emptyset$ if $r\ne s$
and there are no $s$-expansions of the root $r$ with respect to $\Delta.$
Otherwise, the set $Q_{r,s}^\Delta$ consists of all distinct tuples
$(q_1,\ldots,q_t)$, where $q_1,\ldots,q_t\in\Delta$, such that the sum
$s+q_1 +\cdots + q_t$ is an $s$-decomposition of the root $r$ with respect to $\Delta$.

{\bf Definition 2.}
Let $r,s\in\Phi$ and $\emptyset \ne \Delta\subseteq\Pi(\Phi)$.
We put $K_{s,s}^\Delta=1$ and $K_{r,s}^\Delta=0$ if $Q_{r,s}^\Delta=\emptyset.$ In the other cases we put
$$
K_{r,s}^\Delta=\sum_{(q_1,\ldots,q_t)\in Q_{r,s}^\Delta}
N_{s,q_1}N_{s+q_1,q_2}\ldots N_{s+q_1+\cdots+q_{t-1},q_t}. \eqno{(1)}
$$

The concepts listed above, their properties and examples were introduced in \cite{EKL2023}. 

On the set of positive roots, we define the directed graph $G(\Phi^+,\Delta)$ as follows.
The vertices of the graph $G(\Phi^+,\Delta)$ are the roots of the set $\Phi^+.$ Two vertices $s$ and $r,$ where
$s,r\in\Phi^+,$ is connected by an edge directed from $s$ to $r$ if and only 
if there is a fundamental root $q\in\Delta$ such that $r=s+q.$

Let us denote the adjacency matrix of the graph $G(\Phi^+,\Delta)$ by $T_+(\Delta)=(t_{sr})$ that is a matrix of order $|\Phi^+|\times |\Phi^+|$ the element $t_{sr}$ of which is equal to $1$ if there is the edge $(s,r)$ in the graph $G(\Phi^+,\Delta)$ and is equal to zero otherwise.
We also define the matrix $N_+(\Delta)=(n_{sr})$ of order $|\Phi^+|\times|\Phi^+|$ whose element $n_{sr}$ is equal to the constant $N_{s ,q}$ if exists
a simple root $q\in\Delta$ such that $r=s+q,$ and is equal to zero otherwise. The graph $G(\Phi^-,\Delta),$ matrices $T_-(\Delta)$ and $N_-(\Delta)$ are defined in a similar way. The following theorems hold
\medskip

\textbf{Theorem 3.} \textit{Let $r,s\in \Phi^+$ ($r,s\in\Phi^-$) and $r\ne s$, $\Delta\subseteq \Pi(\Phi).$ Then the number of different 
$s$\,-decompositions of the root $r$ with respect to  $\Delta$ is equal to the element at the intersection of the row
$s$ and column $r$ of the matrix $(T_+(\Delta))^k$ (respectively $(T_-(\Delta))^k$), where $k={\rm ht}(r)-{\rm ht}(s).$}
\medskip

\textbf{Theorem 4.} \textit{Let $r,s\in \Phi^+$ ($r,s\in\Phi^-$) and $r\ne s,$ $\Delta\subseteq \Pi(\Phi).$ Then the number $K_{r,s}^{\Delta}$
is equal to the element at the intersection of row $s$ and column $r,$ of the matrix $(N_+(\Delta))^k$ (respectively $(N_-(\Delta))^k$), 
where $k={\rm ht}(r)-{\rm ht}(s).$}

\subsection{Graphs, Matrices, Tables}

Let $\Delta=\{ a,b,c,d\}.$ Below we present the directed weighted graph $G(F_4^-,\Delta),$ its adjacency matrix and
matrix of weights, using them we calculate the table of the number of paths and the table of numbers $K_{r,s}^{\Delta}.$ 
\bigskip

\newpage

\thicklines

\vskip24cm

\newpage

\centerline{\textbf{Adjacency Matrix of the Graph $G(\Phi^-,\Delta)$} }

$$
\left(%
}

\newpage

\centerline{\textbf{Weighted Graph Adjacency Matrix $G(\Phi^-,\Delta)$}}

$$
\left(%
\right)
$$

Further Next, using Theorem~4 we calculate the numbers 
$K_{r,s}^{\Delta}$ for negative $r$ and $s,$ and then,
using Theorem~5 below and for positive $r$ and $s.$
\medskip

{\bf Theorem 5.} 
{\sl For any $r,s\in\Phi^+$ and arbitrary $\Delta\subseteq\Pi(\Phi)$ we have the
equality
$$
K_{r,s}^\Delta=(-1)^{{\rm ht}(r)-{\rm ht}(s)} \frac{(r,r)}{(s,s)} K_{-s,-r}^\Delta.
$$
}

\rotatebox{90}{
%
}

\newpage

Let $\Delta=\{b,c,d\}.$ Removing from the graph shown in Picture~7, all edges marked by the root $a,$ we obtain the graph
$G(F_4^-,\{b,c,d\}).$ Below are its adjacency matrix and weight matrix, from which the table of the number of paths and tables of numbers $K_{r,s}^{\Delta}$ are calculated.
\bigskip

\centerline{\textbf{Graph Adjacency Matrix $G(\Phi^-,\{b,c,d\})$}}

$$
\left(%
}

\newpage

\centerline{\textbf{Weighted Graph Adjacency Matrix $G(\Phi^-,\{b,c,d\})$}}

$$
\left(%
}

\newpage

\subsection{Aplication to Calculate of Commutators in Chevalley Groups}

Let us fix the following notation from the theory of commutators. By definition, we assume $[B,\, _{0}A]=B$ and
$$
[B,\, _{w}A]=[[B,\, _{w-1}A],A]
$$
for $w\geqslant1.$ The following statements hold \cite{EKL2023}.
\bigskip

{\bf Theorem 6.}
{\sl Let $R$ be an associative\,-commutative ring with identity, $\Phi$ be a reduced indecomposable root system,
$s_1,\ldots,s_l\in\Phi^+$ is different roots of height $m,$ $\Delta=\{r_1,\ldots,r_k\}\subseteq\Pi(\Phi).$ Set
$A=x_{r_1}(1)\ldots x_{r_k}(1)$, $B=x_{s_1}(t_1)\ldots x_{s_l}(t_l),$ where $t_1,\ldots, t_l \in R.$
Then for any natural number $w$ we have
$$
 [B,\,_wA]=\prod_{ {\rm ht}(r)=m+w}
    x_r\left(\sum_{j=1}^l t_j K_{r,s_j}^\Delta\right)\ldots\ .
$$
Here the three dots denote the product of root elements corresponding to the root of a height greater than
$m+w.$}

{\bf Theorem 7.}
{\sl Let $R$ be an associative\,-commutative ring with identity, $J=\langle t_1,\ldots, t_l\rangle$ $\subset R$
is a quasiregular ideal, and, $J^2=0;$ also let $\Phi$ be a reduced indecomposable root system,
$s_1,\ldots,s_l\in\Phi^-$ be different roots of height $m,$ $\Delta=\{r_1,\ldots,r_k\}\subseteq\Pi(\Phi).$
Suppose $A=x_{r_1}(1)\ldots x_{r_k}(1)$, $B=x_{s_1}(t_1)\ldots x_{s_l}(t_l).$
Then for any natural number $w$ such that $m+w<0,$ we have
$$
 [B,\,_wA]=\prod_{{\rm ht}(r)=m+w}
    x_r\left(\sum_{j=1}^l t_j K_{r,s_j}^\Delta\right)\ldots\ .
$$
Here the three dots denotes the product of diagonal elements and root elements corresponding to the roots
heights greater than $m+w.$
}

\bigskip

Sergey G. Kolesnikov

Siberian Federal University

email: skolesnikov@sfu-kras.ru

\medskip

Anna I. Polovinkina

Siberian Federal University

email: aiolovinkina@sfu-kras.ru


\medskip
\begin{thebibliography}{1}

\bibitem{Bur72}
{\sl Bourbaki N.} Lie groups and algebras. M.: Mir, 1972. 334 p.

\bibitem{Car72}
{\sl Carter R.} Simple groups of Lie type.-Ney York: Wiley and
Sons, 1972. 458 p.

\bibitem{EKL2023}
{\sl Egorychev G.P., Kolesnikov S.G., Leontiev V.M.}
Necessary and sufficient conditions for the regularity of Sylow p-subgroups of Chevalley groups over $Z_p$ and $Z_{p^2}$. 
{\it Siberian Mathematical Journal}, 2023, Vol. 63, No. 3, p. 500--520.

\end{thebibliography}
\end{document}